\documentclass[12pt, reqno]{amsart}
\usepackage{amsmath}
\usepackage{amssymb,amsthm,amsmath}
\usepackage[numbers,sort&compress]{natbib}
\usepackage{amssymb,amsmath}
\usepackage{amsfonts}
\usepackage{mathrsfs}
\usepackage{latexsym}
\usepackage{amssymb}
\usepackage{amsthm}
\usepackage{indentfirst}
\usepackage{colortbl}
\newtheorem{theorem}{Theorem}[section]

\newtheorem{corollary}{Corollary}[section]

\newtheorem{definition}{Definition}[section]

\newtheorem{lemma}{Lemma}[section]

\newtheorem{proposition}{Proposition}[section]
\newtheorem{remark}{Remark}[section]

\numberwithin{equation}{section}
\newcommand{\norm}[1]{\Vert #1 \Vert}
\newcommand{\normm}[1]{\left\Vert #1 \right\Vert}
\newcommand{\nn}[1]{\left|  #1 \right|}
\newcommand{\lr}[1]{\left(  #1 \right) }

\DeclareMathOperator{\dive}{div}

\DeclareMathOperator{\mm}{\mathrm{d}t }
\DeclareMathOperator{\ms}{\mathrm{d}s }

\newcommand{\fracc}{\frac{\mathrm{d}}{\mathrm{d}t}}
\newcommand{\fr}{\frac{1}{2}}

\newcommand{\inn}{\int_{\Omega}}
\newcommand{\hv}{\hat{v}}

\newcommand{\htt}{\hat{\vartheta}}
\newcommand{\cv}{\check{v}}

\newcommand{\ct}{\check{\vartheta}}

\newcommand{\hj}{\hat{J}}
\newcommand{\cj}{\check{J}}

\newcommand{\ptl}{\partial_{l}}
\newcommand{\ptm}{\partial_{m}}

\newcommand{\dx}{\mathrm{d}x}

\newcommand{\da}{\dive_{A}}
\newcommand{\na}{\nabla_{A}}
\newcommand{\dha}{\dive_{\hat{A}}}
\newcommand{\nca}{\nabla_{\check{A}}}
\newcommand{\dca}{\dive_{\check{A}}}
\newcommand{\nha}{\nabla_{\hat{A}}}
\allowdisplaybreaks
\begin{document}

%
%

\title[Existence and uniqueness of compressible Navier-Stokes equations]{Local existence and uniqueness of heat conductive
compressible Navier-Stokes equations in the presence of vacuum and without initial compatibility conditions}
\thanks{$^*$Corresponding author}

\author[J. Li]{Jinkai Li}
\address[J. Li]{South China Research Center for Applied Mathematics and Interdisciplinary Studies,
School of Mathematical Sciences, South China Normal University, Guangzhou 510631, China}
\email{jklimath@m.scnu.edu.cn; jklimath@gmail.com}

\author[Y. Zheng]{Yasi Zheng$^*$}
\address[Y. Zheng]{South China Research Center for Applied Mathematics and Interdisciplinary Studies,
School of Mathematical Sciences, South China Normal University, Guangzhou 510631, China}
\email{yszheng@m.scnu.edu.cn}


\date{
August 24, 2021}

\begin{abstract}
In this paper, we investigate the initial-boundary value problem to the heat conductive
compressible Navier-Stokes equations. Local existence and uniqueness of
strong solutions is established with any such
initial data that the initial density $\rho_0$, velocity $u_0$, and temperature $\theta_0$ satisfy
$\rho_0\in W^{1,q}$, with $q\in(3,6)$, $u_0\in H^1$, and $\sqrt{\rho_0}\theta_0\in L^2$. The initial density is assumed
to be only nonnegative and thus the initial vacuum is
allowed. In addition to the necessary regularity assumptions, we do not require any initial compatibility conditions
such as those proposed in (Y. Cho and H. Kim, \emph{Existence results for viscous polytropic fluids with vacuum}, J. Differential Equations {\bf 228} (2006), no.~2, 377--411.), which although are widely used in
many previous works but put some inconvenient constraints on the initial data.
Due to the weaker regularities of the initial data and the absence of the initial
compatibility conditions,
leading to weaker regularities of the solutions compared with those in the previous works, the uniqueness of solutions
obtained in the current paper does not follow from the arguments used in the existing literatures. Our proof of the
uniqueness of solutions is based on the following new idea of two-stages argument: (i) showing that the difference
of two solutions (or part of their components) with the same initial data is controlled by some power function
of the time variable; (ii) carrying out some singular-in-time weighted energy differential inequalities fulfilling
the structure of the Gr\"onwall inequality. The existence is established in the Euler coordinates, while the
uniqueness is proved in the Lagrangian coordinates first and then transformed back to the Euler coordinates.
\end{abstract}

\maketitle


\section{Introduction}
\subsection{The compressible Navier-Stokes equations}
Let $\Omega$ be a bounded domain in $\mathbb R^3$ with suitably smooth boundary $\partial\Omega$. Consider
the following heat conductive compressible Navier-Stokes equations in $\Omega\times(0,T)$:
\begin{align}
\label{a}\rho_{t}+\dive(\rho u)=0,\\
\label{b}\rho (u_{t}+u\cdot\nabla u) -\mu\Delta u-(\mu+\lambda)\nabla\dive u+\nabla P=0,\\
\label{c}c_{v}\rho\left( \theta_{t}+ u\cdot\nabla\theta\right) -\kappa\Delta \theta+P\dive u=\mathcal{Q}(\nabla u),
\end{align}
where the unknowns $\rho\geq0, u\in\mathbb R^3,$ and $\theta\geq0$, respectively, represent the density, velocity, and absolute temperature, $p=R\rho\theta$, with positive constant $R$, is the pressure, $c_v$ is a positive constant,
constants $\mu$ and $\lambda$ are the bulk and shear viscous coefficients, respectively, satisfying the physical constraints
$$
\mu>0,\quad 2\mu+3\lambda\geq0,
$$
positive constant
$\kappa$ is the heat conductive coefficient, and
$$\mathcal Q(\nabla u)=\frac\mu2|\nabla u+(\nabla u)^T|^2+\lambda(\text{div}\,u)^2,$$
with $(\nabla u)^T$ being the transpose of $\nabla u$.

There has been many works on the mathematical studies on the compressible Navier-Stokes equations.
In the absence of vacuum, i.e., in the case that the initial density has a uniform positive lower bound, the uniqueness of solutions
was first established by Graffi \cite{G1953} for the isentropic case, and later extended by
Serrin \cite{S1959} to the general case. Local existence of solutions to the compressible Navier-Stokes equations
was first established by Nash \cite{N1962} in the Sobolev type spaces and later by Itaya \cite{Itaya1971}
in the H\"{o}lder type spaces, see also \cite{VOL1972,SOL1977}
for further developments.
Global well-posedness of strong solutions in one dimension with arbitrary large initial data
was first discovered by Kanel \cite{KANEL1968} for the isentropic case, and thereafter extended by Kazhihov--Shelukhin
\cite{KS1977} and  Kazhihov\cite{Ka1982} to the heat conductive case, see also
\cite{CHT2000,JZ2004,LL2016,ZA1998} for some related results. For the multi-dimensional case, global well-posedness
of solutions was first established by Matsumura--Nishida \cite{MN1981,MN1980,MN1982} for small perturbed
initial data around the non-vacuum equilibrium sates in some Sobolev spaces of high order and by Hoff \cite{Hoff1997}
for discontinuous initial data. For the local and global well-posedness of strong solutions in the critical spaces,
one refers to \cite{CMZ2010,DX2018,CD2015,Dan2000,Dan2001,FZZ2018} and the references therein.

In the presence of vacuum, i.e., in the case that the initial density vanishes in some region,
there has been a considerable number of works on the compressible Navier-Stokes equations
since the work of Lions \cite{Lions1998}, where the global existence of weak solutions to the isentropic
compressible Navier-Stokes equations was established with the adiabatic constant $\gamma\geq\frac{9}{5}$.
This was extended to the case that $\gamma\geq\frac{3}{2}$ by Feireisl--Novotn\'{y}--Petzeltov\'{a}
\cite{FNP2001}, and further by Jiang--Zhang \cite{JZ2001,JZ2003} to the case that $\gamma>1$
but only for the spherically symmetric or axisymmetric solutions, see Bresch--Jabin \cite{BJ2018} for some
recent developments where the more general stress tensor and pressure laws are allowed.
Global existence of weak solutions to the full
compressible Navier-Stokes equations under some structure assumptions on the viscous and heat conductive
coefficients as well the equations of states was established by Feireisl
\cite{Feireisl2004}. However, the uniqueness of weak solutions is still an open problem.

Same to the case in the absence of vacuum, one can also establish the local existence and uniqueness of
strong solutions in the presence of vacuum, if the initial data have suitably high regularities.
In fact Salvi--Stra\v{s}kraba \cite{SS1993}, Choe--Kim \cite{CK2003}, and Cho--Choe--Kim \cite{CCK2004} established the
local well-posedness of strong solutions to the isentropic compressible Navier-Stokes equations
with suitably regular initial data satisfying some compatibility condition.
For the full case, Cho--Kim \cite{CK2006} proved the local well-posedness of strong solutions for
the Cauchy problem in $\mathbb R^3$ with initial data satisfying $\rho_0-\rho_\infty\in W^{1,q}(\mathbb R^3)\cap W^{1,r}(\mathbb R^3)$, with $q\in(3,6]$, $\rho_\infty\geq0$, $(u_0,\theta_0)\in D_0^1(\mathbb R^3)\cap D^2(\mathbb R^3)$, and the following compatibility conditions
\begin{eqnarray}
-\mu\Delta u_{0}-(\mu+\lambda)\nabla\dive u_{0}+\nabla(R\rho_{0}\theta_{0})=\sqrt{\rho_{0}}g_{1}, \label{COM1}\\
\kappa\Delta\theta_{0}+\frac{\mu}{2}\nn{\nabla u_{0}+\lr{\nabla u_{0}}^{T}}^{2}+\lambda(\dive u_{0})^{2}=\sqrt{\rho_{0}}g_{2},\label{COM2}
\end{eqnarray}
for $(g_1,g_2)\in L^2(\mathbb R^3)$, where $r=2$ or $3$ if $\rho_\infty>0$ and $r=2$ if $\rho_\infty=0$.
Global existence of strong solutions of small energy but allowed to have large oscillations was first established by Huang--Li--Xin \cite{HUANGLIXIN2012} to the Cauchy problem of
the isentropic compressible Navier-Stokes equations in $\mathbb R^3$; see \cite{LIXIN2019,HUANGLI2018,WENZHU2017,LI2020} for some further developments in this direction. Different from the multi-dimensional case, in the one-dimensional case, the global well-posedness of strong solutions can be established for
arbitrary large initial data for both heat conductive and non-heat-conductive cases, see
\cite{LJ2019,LJ2020,LX2020,LX2020arxiv}. In
particular, local and global well-posedness of entropy-bounded solutions was established firstly in
\cite{LX2020,LX2020arxiv}.

It should be pointed out that the compatibility conditions (\ref{COM1})--(\ref{COM2}) or their natural amendments
play an essential role in the well-posedness theories established in \cite{SS1993,CK2003,CK2006} and, as a results,
they are accepted as standard assumptions
to get the well-posedness of the compressible Navier-Stokes equations in the presence of vacuum. Note that the compatibility
conditions (\ref{COM1})--(\ref{COM2}) ask for some restrictive constraints on the initial data in the vacuum region and
also in the neighborhood of the vacuum-nonvacuum interface. In fact, by the compatibility conditions (\ref{COM1})
and (\ref{COM2}), the initial velocity $u_0$ and
temperature $\theta_0$ are destined to obey
$$
-\mu\Delta u_{0}-(\mu+\lambda)\nabla\dive u_{0}=0\quad\mbox{and}\quad \kappa\Delta\theta_{0}+\frac{\mu}{2}\nn{\nabla u_{0}+\lr{\nabla u_{0}}^{T}}^{2}+\lambda(\dive u_{0})^{2}=0
$$
in the vacuum region, which however seem not physically relevant. From this point of view, the well-posedness theory established in
\cite{SS1993,CK2003,CK2006} does not always match the physical requirements. In particular, it does not always provide
the desired well-posedness for any suitably smooth initial data without any extra constrains.

Due to the analysis in the above paragraph, it is both mathematically and
physically important to establish an alternative well-posedness theory without requiring any initial compatibility conditions
like (\ref{COM1}) and (\ref{COM2}). The first study towards this direction was made by the first author of this paper
in \cite{LJ2017} for the inhomogeneous incompressible Navier-Stokes equations, where the local well-posedness was successfully
established without any compatibility conditions on the initial data, see Danchin--Mucha \cite{DANCHINMUCHA} for further developments aiming
to relax the smoothness of the initial density.
Similar local well-posedness theory without any initial compatibility conditions was later established for the isentropic compressible
Navier-Stokes equations by Gong--Li--Liu--Zhang \cite{GL2020} and Huang \cite{Huang2020} independently. However, for the full
compressible Navier-Stokes equations, to the best of out knowledge, the desired local well-posedness theory without any compatibility
conditions on the initial data has not been established, and only part result is available, see Lai--Xu--Zhang
\cite{LXZ2021}, where they removed (\ref{COM2}) but still required (\ref{COM1}).

The aim of this paper is to establish the desired local well-posedness theory to the full compressible
Navier-Stokes equations without any extra compatibility conditions beyond the necessary smoothness conditions on the initial data.
We also pay some attention to find some minimal regularities on the initial data to guarantee the well-posedness. In this paper, we
consider the initial-boundary value problem; however, the result and method present this paper work also for the Cauchy problem.

The initial and boundary conditions read as:
\begin{eqnarray}
\label{d}&(\rho,\rho u, \rho\theta)|_{t=0}=(\rho_{0}, \rho_{0} u_{0}, \rho_{0}\theta_{0}),\\
\label{e}&u|_{\partial\Omega}=0,\quad\theta|_{\partial\Omega}=0,
\end{eqnarray}
where $\rho_0, u_0,$ and $\theta_0$ are given functions.

It should be pointed out that the real values of $u_0$ and $\theta_0$ that we need are only in the non-vacuum region
$\Omega_+:=\{x\in\Omega|\rho_0(x)>0\}$ but not
in the vacuum region $\Omega_0:=\{x\in\Omega|\rho_0(x)=0\}$. Precisely, denote by $u_0^{\text{real}}$ and
$\theta_0^{\it{real}}$, respectively, the initial velocity and temperature in the non-vacuum region $\Omega_+$, and define
$\mathscr {S}_{\text{ext}}$ as
\begin{align*}
  \mathscr {S}_{\text{ext}}=\Big\{&(\tilde u_0,\tilde\theta_0)\Big|(\tilde u_0, \tilde\theta_0)=(u_0^\text{real},\theta_0^\text{real})
  \text{ on }\Omega_+,
  \tilde u_0\in H_0^1(\Omega), \\
&\tilde\theta_0\text{ is measurable and finitely valued a.e.\,in } \Omega
\Big\},
\end{align*}
then any $(u_0,\theta_0)\in \mathscr S_\text{ext}$ can be chosen as the ``initial"
velocity and temperature without changing the initial condition (\ref{d}).
In fact, for any $(u_0,\theta_0)\in \mathscr S_\text{ext}$, it is clear that
  \begin{equation*}
    {\rho_0}u_0=
    \left\{
    \begin{array}{cr}
    {\rho_0}u_0^{real}&\mbox{ in }\Omega_+,\\
    0&\mbox{ in }\Omega_0,
    \end{array}
    \right.\qquad  {\rho_0}\theta_0=
    \left\{
    \begin{array}{cr}
    {\rho_0}\theta_0^{real}&\mbox{ in }\Omega_+,\\
    0&\mbox{ in }\Omega_0.
    \end{array}
    \right.
  \end{equation*}
Due to the above explanation, throughout this paper, we always assume that the ``initial" velocity $u_0$ and temperature $\theta_0$ are defined
on the whole domain such that $u_0\in H_0^1(\Omega)$ and that $\theta_0$ is Lebesgue measurable and
finitely valued almost everywhere.

%

Before stating the main results, we first clarify some necessary notations
being used throughout this paper and state the definition of solutions to be established.

For $1\leq q\leq\infty$ and positive integer $m$, we use $L^q=L^q(\Omega)$ and $W^{m,q}=W^{m,q}(\Omega)$
to denote the standard Lebesgue and Sobolev spaces, respectively, and we use $H^m$
to replace $W^{m,2}$. For simplicity, we also use notations $L^q$ and
$H^m$ to denote the $N$ product spaces $(L^q)^N$ and $(H^m)^N$, respectively.
We always use $\|u\|_q$ to denote the $L^q$ norm of $u$, while the $L^2$ norm is further simplified as $\|\cdot\|$. For shortening the
expressions, we sometimes use $\|(f_1,f_2,\cdots,f_n)\|_X$ to denote the norm
$\sum_{i=1}^N\|f_i\|_X$ or its equivalence $\left(\sum_{i=1}^N\|f_i\|_X^2
\right)^{\frac12}$.

The strong solutions to be established in this paper are defined as follows.

\begin{definition}
\label{def1}
Given a positive time $ T\in(0,\infty)$ and let $q\in(3,6)$. Assume that $\theta_0$ is nonegative, Lebesgue measurable, and finitely valued a.e.\,in $\Omega$, and that
\begin{equation*}
	0\leq\rho_{0}\in W^{1,q}(\Omega),\quad u_{0}\in H_{0}^{1}(\Omega), \quad\sqrt{\rho_{0}}\theta_{0}\in L^{2}(\Omega).
\end{equation*}
A triple $(\rho, u, \theta)$ is called a strong solution to system \eqref{a}--\eqref{c} in $\Omega\times(0, T)$, subject to \eqref{d}--\eqref{e}, if it has the regularities
 \begin{align*}
&0\leq\rho\in C([0,T];L^{2})\cap L^{\infty}(0,T;W^{1,q}),\qquad\rho_{t}\in L^{\infty}(0,T;L^{2}),\\
&\sqrt\rho u\in C([0,T];L^{2}),\quad u\in L^{\infty}(0,T;H_{0}^{1})\cap L^{2}(0,T;H^{2})
\cap L^1(0,{T_0}; W^{2,q}),\\
&\sqrt{\rho}u_{t}\in L^{2}(0,T;L^{2}),\quad\sqrt{t}u\in L^{\infty}(0,T;H^{2})\cap L^{2}(0,T;W^{2,q}),\\
& \sqrt{t}u_{t}\in L^{2}(0,T;H_{0}^{1}),\quad\sqrt\rho\theta\in C([0,T];L^{2}),\quad0\leq\theta\in L^{2}(0,T;H_{0}^{1}),\\
&\sqrt{t}\theta\in L^{\infty}(0,T;H_{0}^{1})\cap	L^{2}(0,T;H^{2}),\quad\sqrt t\sqrt\rho\theta_t\in
L^2(0,{T_0}; L^2),\\
&t\theta\in L^2(0,{T_0}; W^{2,6}),\quad  t\theta_{t}\in L^{2}(0,T;H_{0}^{1}),
\end{align*}
satisfies equations \eqref{a}--\eqref{c} a.e.\,in $\Omega\times(0,T)$, and fulfills the initial condition \eqref{d}.
\end{definition}

\begin{remark}
By the regularities of $\rho$ stated in Definition \ref{def1}, it follows from the Gagliardo-Nirenberg inequality that
$\rho\in C([0, T]; C(\overline\Omega))$. Thanks to this and recalling that $\sqrt\rho u,\sqrt\rho\theta\in C([0,T];L^{2})$, it is clear that $\rho u,\rho\theta\in C([0, T]; L^2)$.
Therefore, the initial values of $\rho u$ and $\rho\theta$ are well-defined.
\end{remark}

We are now ready to state the main result of this paper.
\begin{theorem}\label{u}
Let $q\in(3,6)$. Assume that $\theta_0$ is nonegative, Lebesgue measurable, and finitely valued a.e.\,in $\Omega$, and that
\begin{eqnarray*}
	0\leq\rho_{0}\in W^{1,q}(\Omega),\quad u_{0}\in H_{0}^{1}(\Omega), \quad\sqrt{\rho_{0}}\theta_{0}\in L^{2}(\Omega).
\end{eqnarray*}
Then, there exists a positive time ${T_0}$ depending only on $R, \mu, \lambda, c_{v}, \gamma, q,$ and  $\Phi_{0}$, such that system
\eqref{a}--\eqref{c}, subject to \eqref{d}--\eqref{e}, has a unique strong solution in $\Omega\times(0,{T_0})$, where $\Phi_0:=\norm{\rho_{0}}_{\infty}+\norm{\nabla\rho_{0}}_{q}+\norm{(\sqrt{\rho_{0}}\theta_{0},\nabla u_{0})}^{2}$.
\end{theorem}

\begin{remark}
(i) No compatibility conditions as those in \cite{SS1993,CK2006,CK2003,CCK2004,LXZ2021} are required in Theorem \ref{u}. Comparing with
the result proved in \cite{LXZ2021}, where compatibility condition (\ref{COM2}) was removed but (\ref{COM1}) was still required, in Theorem \ref{u}, both
(\ref{COM1}) and (\ref{COM2}) were removed.

(ii) The arguments present in this paper with slightly modifications work also for the Cauchy problem and similar result as in Theorem \ref{u} still holds, with the assumptions on $u_0$ and $\theta_0$ replaced by
$(u_0,\sqrt{\rho_0}\theta_0)\in D_0^1\times L^2.$ Note that these assumptions are weaker than those in
\cite{CK2006,LXZ2021}. In fact, \cite{CK2006}
requires $(u_0, \theta_0)\in D_0^1\cap D^2$ while \cite{LXZ2021} requires $u_0\in D_0^1\cap D^2$ and $\theta_0\in D^1$.
\end{remark}

The key of proving the existence part of Theorem \ref{u} is to carry out some suitable
a priori estimates of the following quantity
 \begin{align}
\nonumber\Phi(t):=&\int_{0}^{t}\left( \norm{\sqrt{\rho}u_{t}}^{2}+\norm{\nabla^{2}u}^{2}+\norm{\sqrt{t}\nabla u_{t}}^{2}+\norm{\nabla\theta}^{2}+\norm{\sqrt{t}\sqrt{\rho}\theta_{t}}^{2}+\norm{\sqrt{t}\nabla^{2}\theta}^{2}\right) \mathrm{d}s \\
&+\sup_{0\leq s\leq t} \left( \norm{\rho}_{\infty}+\norm{\nabla\rho}_{q}+\norm{\sqrt{\rho}\theta}^{2}+\norm{\nabla u}^{2}+\norm{\sqrt{t}\nabla\theta}^{2}\right)+1 \label{PHI}
\end{align}
for any approximate solution $(\rho, u, \theta)$ to system \eqref{a}--\eqref{c}, subject to \eqref{d}--\eqref{e}.
The estimates for $\Phi$ have to be independent of
the compatibility conditions. Roughly speaking, the estimate of $\Phi(T)$ is achieved
based on the following conditional a priori estimate by the continuity argument: it holds that $\Phi(T)\leq C$, as long as
$T^{\frac{6-q}{4q}}\Phi^{2}(T)\leq\epsilon_{0}$, where $\epsilon_0$ and $C$ are two positive constants
independent of the compatibility conditions, see Corollary \ref{x}. With this at hand,
one can get the time weighted higher order a priori estimates. Then, one can prove that
the existence time of the approximate solutions can be chosen independent of the initial
compatibility conditions, so are the corresponding a priori estimates. These will then yield a preparing existence result
by passing the limit to the approximate solutions;
however, the regularities that $\sqrt\rho u,\sqrt\rho\theta\in C([0,T];L^{2})$ are not guaranteed in this
passage, and instead
what we have are $\rho u\in C([0,T]; L^2)$ and $\rho\theta\in C_w([0,T]; L^2)$, here $C_w$ represents
the weak continuity. As a compensation, we prove that $\|\sqrt{\rho}u\|^2(t)
\leq\|\sqrt{\rho_{0}}u_{0}\|^2+Ct$ and $\norm{\sqrt{\rho}\theta}^{2}(t)\leq\norm{\sqrt{\rho_{0}}\theta_{0}}^{2}+C\sqrt t$,
which are employed to prove the continuities with respect to time of $\sqrt\rho u$ and $\sqrt\rho\theta$ in $L^2$
in the Lagrangian coordinates first and finally transformed back to those in the Euler coordinates.

Due to the lower regularities on the initial data and the absence of the initial compatibility conditions, the regularities of the
strong solutions obtained in this paper are weaker than those required in \cite{SS1993,CK2006,CCK2004,LXZ2021} to prove the
uniqueness in their ways. We also note that even though the uniqueness was achieved in \cite{LJ2017,GL2020,Huang2020,DANCHINMUCHA}
for the inhomogeneous incompressible Navier-Stokes equations and the isentropic compressible Navier-Stokes equations
with the initial data $(\rho_0, u_0)$ satisfying some similar regularities as in this paper
and without any initial compatibility conditions, still the arguments in these works do not apply to the current paper.
The main reasons are that the regularities of the initial temperature assumed in this paper are weaker than those of the initial velocity,
and even worse that the entropy production term $\mathcal Q(\nabla u)$ has stronger nonlinearities than the convection
terms. As a result, in matter of overcoming the difficulties cased by the low regularities and lack of compatibility conditions on the initial
data, the ideas used to deal with the velocity are not sufficient to deal with the temperature.

Our strategies of proving the uniqueness are illustrated as follows.
Let $(\overline\rho, \overline u, \overline\theta)$ and $(\hat\rho, \hat u, \hat\theta)$
be two solutions with the same initial data and denote by $(\rho, u, \theta)$
their subtraction. Then, one has some differential inequalities of the form
\begin{eqnarray}
  \frac{d}{dt}\|\sqrt{\overline\rho}u\|^2+ \|\nabla u\|^2&\leq& C \|\nabla\hat u_t\|^2 \|\rho\|^2+\cdots,\label{DF1}\\
  \frac{d}{dt}\|\sqrt{\overline\rho}\theta\|^2+\|\nabla\theta\|^2&\leq& C \|\nabla^2\overline u\| \|\nabla u\|^2
  +C\|\nabla\overline\theta_t\|^2\|\rho\|^2 +\cdots,\label{DF2} \\
  \frac{d}{dt}\|\rho\|^2&\leq& C\|\nabla\overline u\|_\infty\|\rho\|^2+C\|\nabla u\|\|\rho\|,\label{DF3}
\end{eqnarray}
where all other quantities that can be dealt with relatively easier are omitted in the suspension points.
We want to derive some Gr\"onwall type structure from the above inequalities.
Recalling that one only has $t\nabla\overline\theta_t\in L^2(0,T; L^2)$, the hardest
term $\|\nabla\overline\theta_t\|^2\|\rho\|^2$ in (\ref{DF2})
has to be dealt with as $\|t\nabla\overline\theta_t\|^2\frac{\|\rho\|^2}{t^2}$ and,
as a result, one has to consider the differential inequality for $\frac{\|\rho\|^2}{t^2}$, which can be derived from (\ref{DF3}) as
$$
\frac{d}{dt}\frac{\|\rho\|^2}{t^2}+\frac{\|\rho\|^2}{t^3}\leq  C\|\nabla\overline u\|_\infty\frac{\|\rho\|^2}{t^2}+C\frac{\|\nabla u\|^2}{ t}.
$$
This motivates us to divide (\ref{DF1}) with $t$, leading to
\begin{equation}
\label{DF4}
  \frac{d}{dt}\frac{\|\sqrt{\overline\rho}u\|^2}{t}+ \frac{\|\sqrt{\overline\rho}u\|^2}{t^2}+\frac{\|\nabla u\|^2}t
   \leq  C\|\sqrt t\nabla\hat u_t\|^2 \frac{\|\rho\|^2}{t^2}+\cdots.
\end{equation}
Noticing that $\sqrt t\nabla^2\overline u\in L^\infty(0,T; L^2)$, one can derive from the above two and (\ref{DF2}) that
$$
\frac{d}{dt}\left(\|\sqrt{\overline\rho}\theta\|^2+\frac{\|\sqrt{\overline\rho}u\|^2}{t}+\frac{\|\rho\|^2}{t^2}\right)
\leq  C(\|(\sqrt t\nabla\hat u_t,t\nabla\overline\theta_t)\|^2+\|\nabla\overline u\|_\infty) \frac{\|\rho\|^2}{t^2}+\cdots,
$$
which meets the Gr\"onwall type structure. It remains to guarantee that the quantity with singular weights
$\|\sqrt{\overline\rho}\theta\|^2+\frac{\|\sqrt{\overline\rho}u\|^2}{t}+\frac{\|\rho\|^2}{t^2}$ tends to zero
when approaching the initial time. This is expected to be verified from (\ref{DF1})--(\ref{DF3}) by using the
regularities of the solutions.

It is worth to point out some technique points in the arguments explained in the above paragraph. First, in deriving
the singular $t$-weighted energy inequality (\ref{DF4}), one also encounters a term of the form
$\frac1t\int_\Omega\overline\rho\theta\text{div} udx$. To deal with this term,
we need the $L^\infty(0,T; L^3)$ bound of $\nabla\sqrt{\overline \rho}$ to match the singular weights and, as a result,
one needs the extra condition that $\nabla\sqrt{\rho_0}\in L^3$. However, this
is not assumed in Theorem \ref{u}. Second, in order to prove the uniqueness in the way as
explained in the above paragraph, one needs to
 show that the initial values of $\sqrt{\bar\rho}u$
and $\sqrt{\overline\rho}\theta$ are identically zero. However, it is a subtle issue to verify this in the Euler
coordinates, as the initial condition is
$(\overline\rho\bar u,\overline\rho\overline\theta)|_{t=0}=(\hat\rho\hat u,\hat\rho\hat\theta)|_{t=0}$. Due to the above
two technical difficulties, even though we use the ideas explained as in the previous paragraph to prove the
uniqueness, our proof of the uniqueness is actually carried out in the Lagrangian coordinates first and then
transformed back to the Euler coordinates. Concerning the first technique
point mentioned above, it turns out that, in the Lagrangian
coordinates, the term corresponding to $\int_\Omega\overline\rho\theta\text{div} udx$ reads as $\int_\Omega
\rho_0\vartheta\text{div}_{\bar A}vdy$, for which one can make use of the information $\sqrt{\rho_0}\vartheta$
to avoid the requirement $\nabla\sqrt{\rho_0}\in L^3$, where $\bar A$ is the deformation matrix of the transformation
between the Euler and Lagrangian coordinates. As for the second technique point, in the Lagrangian coordinate, the
corresponding requirement is then $\sqrt{\rho_0}\hat u=\sqrt{\rho_0}\bar u$ at the initial time of which the proof
is given in Proposition \ref{cn}. Finally, we would like to remak that the singular weights used in the Lagrangian
coordinates are actually less singular than those in the Euler coordinates. This also reflects another advantage of
using the Lagrangian coordinates to prove the uniqueness.

Throughout this paper, we use $C$, which may vary from place to place, to denote a generic constant depending only on $R, \mu, \lambda, c_{v}, \gamma, q,$ and the upper bound of $\Phi_{0}$ unless we clearly specify.


\section{A priori estimates independent of compatibility conditions}
The aim of this section is to derive some a priori estimates for the strong solutions to system \eqref{a}--\eqref{c},
subject to \eqref{d}--\eqref{e}, with initial data satisfying some compatibility conditions. We emphasize that although we
assume the initial compatibility conditions, the a priori estimate established in this section do not depend on these
conditions. This is crucial to finally establish the existence of strong solutions without any compatibility conditions.

We start with the following local well-posedness result which can be proved in the same way as in \cite{CK2006} where the compatibility
conditions are required.


\begin{proposition}\label{f}
Let $q\in(3,6]$ and assume that $(\rho_{0},u_{0},\theta_{0})$ satisfies
\begin{align}
	\nonumber0\leq\rho_{0}\in W^{1,q}(\Omega),\quad u_{0}\in H_{0}^{1}(\Omega)\cap H^{2}(\Omega),\quad 0\leq \theta_{0}\in H_{0}^{1}(\Omega)\cap H^{2}(\Omega),
\end{align}
and the compatibility conditions
\begin{align}
	\nonumber-\mu\Delta u_{0}-(\mu+\lambda)\nabla\dive u_{0}+\nabla(R\rho_{0}\theta_{0})=\sqrt{\rho_{0}}g_{1},\\
	\nonumber\kappa\Delta\theta_{0}+\frac{\mu}{2}\nn{\nabla u_{0}+\lr{\nabla u_{0}}^{T}}^{2}+\lambda(\dive u_{0})^{2}=\sqrt{\rho_{0}}g_{2},
\end{align}
for some $g_{1},g_{2}\in L^{2}(\Omega)$.

Then, there exists a positive time $T_{*}$ depending on $R$, $\mu$, $\lambda$, $c_{v}$, $\gamma$, $q$,
$\norm{\nabla^{2}u_{0}}$, $\norm{\nabla^{2}\theta_{0}}$, $\norm{g_{1}}$, and $\norm{g_{2}}$, such that system
\eqref{a}--\eqref{c}, subject to \eqref{d}--\eqref{e}, admits a unique strong solution $(\rho,u,\theta)$ in
$\Omega\times(0,T_{*})$, satisfying
\begin{eqnarray*}
	\nonumber\rho\in C([0,T_{*}];W^{1,q}),\qquad \rho_{t}\in C([0,T_{*}];L^{q}), \\
	\nonumber(u_{t},\theta_{t})\in L^{2}(0,T_{*};H_{0}^{1}),\qquad (\sqrt{\rho}u_{t},\sqrt{\rho}\theta_{t})\in L^{\infty}(0,T_{*};L^{2}),\\
\nonumber (u,\theta)\in C([0,T_{*}];H_{0}^{1}\cap H^{2})\cap L^{2}(0,T_{*};W^{2,q}).
\end{eqnarray*}
\end{proposition}

It will be shown in this section that the existence time $T_{*}$ in the above proposition can be chosen depending only on $R, \mu, \lambda, c_{v}, \gamma, q,$ and the upper bound of
\begin{align}
	\nonumber\Phi_{0}:=\norm{\rho_{0}}_{\infty}+\norm{\nabla\rho_{0}}_{q}+\norm{(\sqrt{\rho_{0}}\theta_{0},\nabla u_{0})}^{2}.
\end{align}
In particular, $T_{*}$ can be chosen independent of $\norm{\nabla^{2}u_{0}}, \norm{\nabla^{2}\theta_{0}}, \norm{g_{1}}$, and $\norm{g_{2}}$. Let $\Phi$ be the quantity given by (\ref{PHI}).
The main issue of this section is to derive the local in time estimate of $\Phi$ independent of $\norm{\nabla^{2}u_{0}}, \norm{\nabla^{2}\theta_{0}}, \norm{g_{1}}$, and $\norm{g_{2}}$, and therefore independent of the initial compatibility conditions.

In the rest of this section until the last proposition, we always assume that $(\rho, u, \theta)$ is a solution to system
\eqref{a}--\eqref{c}, subject to \eqref{d}--\eqref{e}, in $\Omega\times(0,T)$, for some positive time $T\leq1$,
satisfying the regularities in Proposition \ref{f} with $T_{*}$ there replaced by $T$. We emphasize again that
$C$, which may vary from place to place, is a generic constant depending only on
$R,\ \mu,\ \lambda,\ c_{c},\ \gamma,\ q,$ and the upper bound of $\Phi_{0}$.

\begin{proposition}\label{g}
It holds that
\begin{align}
\nonumber\int_{0}^{T}\left( \norm{\nabla u}_{\infty}+\norm{\nabla^{2}u}_{q}\right) \mm\leq CT^{\frac{6-q}{4q}}\Phi^{2}(T).
\end{align}
\end{proposition}
\begin{proof}
Applying the elliptic estimates to \eqref{b}, one obtains
$$\norm{\nabla^{2}u}_{q}\leq C(\norm{\rho u_{t}}_{q}+\norm{\rho(u\cdot\nabla)u}_{q}+\norm{\nabla P}_{q}).$$
It follows from the H\"{o}lder, Gagliardo-Nirenberg, Sobolev, and Poincar\'e inequalities that
\begin{eqnarray*}
	\norm{\rho u_{t}}_{q}
	&\leq& C\norm{\rho}_{\infty}^{\frac{1}{2}}\norm{\sqrt{\rho} u_{t}}^{\frac{6-q}{2q}}\norm{\sqrt{\rho} u_{t}}_{6}^{\frac{3q-6}{2q}}\leq C\norm{\rho}_{\infty}^{\frac{5q-6}{4q}}\norm{\sqrt{\rho }u_{t}}^{\frac{6-q}{2q}}\norm{\nabla u_{t}}^{\frac{3q-6}{2q}},\\
	\norm{\rho(u\cdot\nabla)u}_{q}&\leq&\norm{\rho}_{\infty}\norm{u}_{\infty}\norm{\nabla u}_{q}\leq C\norm{\rho}_{\infty}\norm{\nabla u}^{\frac{1}{2}}\norm{\nabla^{2} u}^{\frac{3}{2}},\\
	\norm{\nabla P}_{q}&\leq& C(\norm{\nabla\rho}_{q}\norm{\theta}_{\infty}+\norm{\rho}_{\infty}\norm{\nabla\theta}_{q})\leq C(\norm{\nabla\rho}_{q}+\norm{\rho}_{\infty})\norm{\nabla\theta}_{q}\\
	&\leq& C(\norm{\nabla\rho}_{q}+\norm{\rho}_{\infty})\norm{\nabla\theta}^{\frac{6-q}{2q}}\norm{\nabla^{2}\theta}^{\frac{3q-6}{2q}}.
\end{eqnarray*}
Integrating the above estimates from 0 and $T$, one gets
\begin{eqnarray*}
\int_{0}^{T}\norm{\rho u_{t}}_{q}\mm
&\leq& C\int_{0}^{T}\norm{\rho}_{\infty}^{\frac{5q-6}{4q}}\norm{\sqrt{\rho }u_{t}}^{\frac{6-q}{2q}}\norm{\nabla u_{t}}^{\frac{3q-6}{2q}}\mm\\
&\leq& C\Phi^{\frac{5q-6}{4q}}(T)\left(\int_{0}^{T}\norm{\sqrt{\rho}u_{t}}^{2} \mm\right)^{\frac{6-q}{4q}} \lr{\int_{0}^{T}\norm{\sqrt{t}\nabla u_{t}}^{2}\mm}^{\frac{3q-6}{4q}}T^{\frac{6-q}{4q}}\\
&\leq& CT^{\frac{6-q}{4q}}\Phi^{\frac{7q-6}{4q}}(T),
\end{eqnarray*}
and
\begin{eqnarray*}
	\int_{0}^{T}\norm{\rho(u\cdot\nabla)u}_{q}\mm&\leq& C\int_{0}^{T}\norm{\rho}_{\infty}\norm{\nabla u}^{\frac{1}{2}}\norm{\nabla^{2} u}^{\frac{3}{2}}\mm\\
	&\leq& \Phi^{\frac{5}{4}}(T)\left(\int_{0}^{T} \norm{\nabla^{2} u}^{2}\mm\right) ^{\frac{3}{4}}T^{\frac{1}{4}} \leq CT^{\frac{1}{4}}\Phi^{2}(T),
\end{eqnarray*}
as well as
\begin{eqnarray*}
\int_{0}^{T}\norm{\nabla P}_{q}\mm&\leq&
C\int_{0}^{T}(\norm{\nabla\rho}_{q}+\norm{\rho}_{\infty})\norm{\nabla\theta}^{\frac{6-q}{2q}}
\norm{\nabla^{2}\theta}^{\frac{3q-6}{2q}}\mm\\
&\leq& C\Phi(T)\left( \int_{0}^{T}\norm{\nabla\theta}^{2}\mm\right) ^{\frac{6-q}{4q}}\left( \int_{0}^{T}\norm{\sqrt{t}\nabla^{2}\theta}^{2}\mm\right)^{\frac{3q-6}{4q}} T^{\frac{6-q}{4q}}\\
&\leq& CT^{\frac{6-q}{4q}}\Phi^{\frac{3}{2}}(T).
\end{eqnarray*}
Therefore, we show
\begin{equation*}
	\int_{0}^{T}\norm{\nabla^{2}u}_{q}\mm\leq C\left( T^{\frac{6-q}{4q}}\Phi^{\frac{7q-6}{4q}}(T)+T^{\frac{1}{4}}\Phi^{2}(T)+T^{\frac{6-q}{4q}}\Phi^{\frac{3}{2}}(T)\right)
	\leq CT^{\frac{6-q}{4q}}\Phi^{2}(T),
\end{equation*}
where $\frac{1}{4}\geq\frac{6-q}{4q}$ for $q\in(3,6)$, $T\leq 1$, and $\Phi(T)\geq1$ were used. Thanks to this, it follows from the
Sobolev and the Poincar\'e inequalities that
 $$
 \int_{0}^{T}\norm{\nabla u}_{\infty}\mm
 \leq\int_{0}^{T}\norm{\nabla^{2}u}_{q}\mm
 \leq CT^{\frac{6-q}{4q}}\Phi^{2}(T).
 $$
The proof is complete.
\end{proof}
 \begin{proposition}\label{h}
 It holds that
 \begin{align*}
 	&\sup_{0\leq t\leq T}\norm{\rho}_{\infty}\leq \norm{\rho_{0}}_{\infty}\exp\left\{CT^{\frac{6-q}{4q}}\Phi^{2}(T)\right\},\\
 	&\sup_{0\leq t\leq T}\norm{\rho}_{W^{1,q}} \leq C\left( 1+T^{\frac{6-q}{4q}}\Phi^{2}(T)\right) \exp\left\{CT^{\frac{6-q}{4q}}\Phi^{2}(T)\right\}.
 \end{align*}
  \end{proposition}
\begin{proof}
For any given $x\in\Omega$ and $s\in[0,T]$, define $U(x,t;s)$ as
\begin{equation*}
	\left\{
	\begin{aligned}
	&\frac{\mathrm{d}}{\mathrm{d}t}U(x,t;s)=u(U(x,t;s),t),\qquad \forall\ t\in[0,T],\\
	&U(x,s;s)=x.
	\end{aligned}
	\right.
\end{equation*}
 Note that $u\in L^{1}(0,T;W^{1,\infty})$ guaranteed by Proposition \ref{g}, $U$ is well-defined. Then, it follows from \eqref{a} that
$$\frac{\mathrm{d}}{\mathrm{d}t}\rho(U(x,t;s),t)=-\dive u(U(x,t;s),t)\rho(U(x,t;s),t),\qquad\forall\ t\in(0,T).$$
Solving the above ordinary differential equation yields
\begin{align*}
\rho(U(x,t;s),t)=\rho(x,s)e^{-\int_{s}^{t}\dive u(U(x,\tau;s),\tau)\mathrm{d}\tau },\qquad \forall\ s,t\in[0,T].
\end{align*}
Choosing $t=0$ in the above, one gets
\begin{align*}
\rho(x,s)=\rho_{0}(U(x,0;s))e^{-\int_{0}^{s}\dive u(U(x,\tau;s),\tau)\mathrm{d}\tau }, \qquad \forall\ (x,s)\in\Omega\times[0,T].
\end{align*}
Therefore, by Proposition \ref{g}, it holds that
\begin{align*}
\sup_{0\leq t\leq T}\norm{\rho}_{\infty}(t)\leq
\norm{\rho_{0}}_{\infty}e^{\int_{0}^{T}\norm{\dive u}_{\infty}(\tau)\mathrm{d}\tau }\leq \norm{\rho_{0}}_{\infty}e^{CT^{\frac{6-q}{4q}}\Phi^{2}(T)},
\end{align*}
proving the first conclusion.

Multiplying \eqref{a} with $\rho^{q-1}$ and integrating over $\Omega$, one deduces by integration by parts that
\begin{align}
	\nonumber\frac{\mathrm{d}}{\mathrm{d}t}\norm{\rho}_{q}^{q}\leq C\norm{\nabla u}_{\infty}\norm{\rho}_{q}^{q}.
\end{align}
Applying the operator $\nabla$ to \eqref{a} and multiplying the resultant with $|\nabla\rho|^{q-2}\nabla\rho$, it follows from
integration by parts that
\begin{eqnarray}
	\nonumber\frac{\mathrm{d}}{\mathrm{d}t}\norm{\nabla\rho}_{q}^{q}&\leq& C\int_{\Omega}\left( |\nabla u||\nabla\rho|^{q}+|\rho||\nabla\rho|^{q-1}|\nabla^{2}u|\right) \mathrm{d}x\\
	\nonumber&\leq& C(\norm{\nabla u}_{\infty}\norm{\nabla\rho}_{q}^{q}+\norm{\rho}_{\infty}\norm{\nabla\rho}_{q}^{q-1}\norm{\nabla^{2}u}_{q}).
\end{eqnarray}
Hence,
\begin{align}
\frac{\mathrm{d}}{\mathrm{d}t}\norm{\rho}_{W^{1,q}}=\frac{\mathrm{d}}{\mathrm{d}t}\left( \norm{\rho}_{q}+\norm{\nabla\rho}_{q}\right)\leq C(\norm{\nabla u}_{\infty}\norm{\rho}_{W^{1,q}}+\|\rho\|_\infty\norm{\nabla^{2}u}_{q}),
\end{align}
from which, by the Gr\"onwall inequality, using the first conclusion, and by Proposition \ref{g}, it follows that
\begin{eqnarray*}
\norm{\rho}_{W^{1,q}}&\leq&\left( \norm{\rho_{0}}_{W^{1,q}}+\int_{0}^{T}\|\rho\|_\infty\norm{\nabla^{2}u}_{q}\mm \right) e^{C\int_{0}^{T}\norm{\nabla u}_{\infty}\mm }\\
&\leq& C\left( 1+T^{\frac{6-q}{4q}}\Phi^{2}(T)\right)e^{CT^{\frac{6-q}{4q}}\Phi^{2}(T)},
\end{eqnarray*}
proving the second conclusion.
\end{proof}

As a direct corollary of Proposition \ref{h}, one obtains:

\begin{corollary}\label{d1}
	There is a sufficiently small positive constant $\epsilon_{0}\leq1$ depending only on $R, \mu, \lambda, c_{v}, \gamma, q,$ and $\Phi_{0}$, such that
	\begin{align*}
	\sup_{0\leq t\leq T}\norm{\rho}_{\infty}\leq2\norm{\rho_{0}}_{\infty},\qquad\sup_{0\leq t\leq T}\norm{\rho}_{W^{1,q}}\leq C,
	\end{align*}
	as long as
	\begin{align}
	\label{d2}
    T^{\frac{6-q}{4q}}\Phi^{2}(T)\leq\epsilon_{0}.
	\end{align}
\end{corollary}

Under the assumption (\ref{d2}) and since $\Phi(T)\geq1$, it is easy to check that the following relations hold:
\begin{eqnarray}
\label{dj}T\Phi^{3}(T)&=&\left(T^{\frac{6-q}{4q}}\Phi^{2}(T) \right)^{\frac{3}{2}}T^{\frac{11q-18}{8q}} \leq\epsilon_{0}^{\frac{3}{2}}\leq1,\\
\label{dk}T^{\frac{1}{2}}\Phi^{2}(T)&=&\left(T^{\frac{6-q}{4q}}\Phi^{2}(T) \right)T^{\frac{3q-6}{4q}}\leq\epsilon_{0}\leq1,\\
\label{dl}T^{\frac{1}{4}}\Phi^{\frac{3}{2}}(T)&\leq&\left(T^{\frac{6-q}{4q}}\Phi^{2}(T) \right)^\frac34T^{\frac{7q-18}{16q}}\leq\epsilon_{0}^\frac34\leq1.
\end{eqnarray}

These will be frequently used in the rest of this section.

\begin{proposition}\label{p}
Let $\epsilon_{0}$ be the number stated in Corollary \ref{d1} and assume that \eqref{d2} holds. Then, the following estimate holds
	\begin{align}
		\nonumber\sup_{0\leq t\leq T}\norm{\sqrt{\rho}\theta}^{2}+\int_{0}^{T}\norm{\nabla\theta}^{2}\mm\leq C.
	\end{align}
\end{proposition}

\begin{proof}
Multiply \eqref{c} with $\theta$ and integrate it over $\Omega$ to get
\begin{align}
\nonumber\frac{c_{v}}{2}\frac{\mathrm{d}}{\mathrm{d}t}\norm{\sqrt{\rho}\theta}^{2}+\kappa\norm{\nabla\theta}^{2}=-\int_{\Omega}\dive uP\theta\mathrm{d}x+\int_{\Omega}\mathcal{Q}(\nabla u) \theta\mathrm{d}x.
\end{align}
The terms on the right-hand side are estimated by the H\"{o}lder, Sobolev, and Young inequalities as
\begin{eqnarray}
\nonumber\int_{\Omega}\dive uP\theta\mathrm{d}x&\leq& R\int_{\Omega}\rho|\theta|^{2}|\nabla u|\mathrm{d}x\leq C\norm{\nabla u}_{\infty}\norm{\sqrt{\rho}\theta}^{2},\\
\nonumber\int_{\Omega}\mathcal{Q}(\nabla u)\theta\mathrm{d}x
&\leq& C\norm{\nabla u}\norm{\nabla u}_{3}\norm{\theta}_{6}\leq C\norm{\nabla u}^{\frac{3}{2}}\norm{\nabla^{2} u}^{\frac{1}{2}}\norm{\nabla\theta}\\
\nonumber&\leq&\frac{\kappa}{2}\norm{\nabla\theta}^{2}+C\norm{\nabla u}^{3}\norm{\nabla^{2}u}.
\end{eqnarray}
Therefore, it follows that
\begin{align}
\nonumber c_{v}\frac{\mathrm{d}}{\mathrm{d}t}\norm{\sqrt{\rho}\theta}^{2}+\kappa\norm{\nabla\theta}^{2}\leq C(\norm{\nabla u}_{\infty}\norm{\sqrt{\rho}\theta}^{2}+\norm{\nabla u}^{3}\norm{\nabla^{2}u}),
\end{align}
from which, by the Gr\"onwall inequality, it follows from the H\"{o}lder inequality, Proposition \ref{g}, and Corollary \ref{d1} that
\begin{eqnarray}
\nonumber&&\sup_{0\leq t\leq T}c_{v}\norm{\sqrt{\rho}\theta}^{2}+\kappa\int_{0}^{T}\norm{\nabla\theta}^{2}\mm\\
\nonumber&\leq&\left(c_{v} \norm{\sqrt{\rho_{0}}\theta_{0}}^{2}+C\int_{0}^{T}\norm{\nabla u}^{3}\norm{\nabla^{2}u}\mm\right)e^{C\int_{0}^{T}\norm{\nabla u}_{\infty}\mm }\\
\nonumber& \leq& C\left( 1+T^{\frac{1}{2}}\Phi^{2}(T)\right)e^{CT^{\frac{6-q}{4q}}\Phi^{2}(T) }\leq C,
\end{eqnarray}
where in the last step (\ref{d2}) and (\ref{dl}) were used. This completes the proof.
\end{proof}

\begin{proposition}\label{s}
	Under the assumptions of Proposition \ref{p}, it holds that
\begin{align}
\nonumber\sup_{0\leq t\leq T}\norm{\nabla u}^{2}+\int_{0}^{T}\left( \norm{\sqrt{\rho}u_{t}}^{2}+\norm{\nabla^{2} u}^{2}\right)\mm\leq C.
\end{align}
\end{proposition}

\begin{proof}
By Corollary \ref{d1} and the H\"{o}lder, Gagliardo-Nirenberg, and Young inequalities, one deduces
	\begin{eqnarray*}
		\nonumber\norm{\nabla^{2}u}^{2}&\leq& C\left(\norm{\rho u_{t}}^{2}+\norm{\rho(u\cdot\nabla)u}^{2}+\norm{\nabla P}^{2} \right)\\
		\nonumber&\leq &C\left(\norm{\rho}_{\infty}\norm{\sqrt{\rho} u_{t}}^{2}+\norm{\rho}_{\infty}\norm{u}_{6}^{2}\norm{\nabla u}^{2}_{3}+\norm{\nabla\rho}_{3}^{2}\norm{\theta}_{6}^{2}+\norm{\rho}_{\infty}^{2} \norm{\nabla\theta}^{2} \right) \\
		\nonumber&\leq& C\left( \norm{\sqrt\rho u_{t}}^{2}+ \norm{\nabla u}^{3}\norm{\nabla^{2}u}+ \norm{\nabla\theta}^{2} \right)\\
		 &\leq&\frac{1}{2}\norm{\nabla^{2}u}^{2}+C\left(\norm{\sqrt{\rho}u_{t}}^{2}+\norm{\nabla u}^{6}+\norm{\nabla\theta}^{2}\right)
	\end{eqnarray*}
and, thus,
\begin{equation}
  \label{o}
  \norm{\nabla^{2}u}^{2}\leq C\left(\norm{\sqrt{\rho}u_{t}}^{2}+\norm{\nabla u}^{6}+\norm{\nabla\theta}^{2}\right).
\end{equation}
Note that \eqref{c} implies
\begin{equation}
\label{EQP}
P_{t}=(\gamma-1)\bigg(\mathcal{Q} (\nabla u)+\kappa\Delta\theta-P\dive u-c_{v}\dive (\rho u\theta)\bigg).
\end{equation}
Thus, by integration by parts, one gets
	\begin{align}
		\nonumber& \int_{\Omega}\nabla Pu_{t}\mathrm{d}x =-\frac{d}{dt}\int_\Omega P\text{div}udx+\int_\Omega P_t\text{div}u dx\\
		\nonumber=&-\frac{\mathrm{d}}{\mathrm{d}t}\int_{\Omega}P\dive u\mathrm{d}x+(\gamma-1)\int_{\Omega}\dive u\bigg(\mathcal{Q}(\nabla u)+\kappa\Delta\theta-P\dive u -c_{v}\dive (\rho u\theta)\bigg)\mathrm{d}x\\
		\nonumber=&-\frac{\mathrm{d}}{\mathrm{d}t}\int_{\Omega}P\dive u\mathrm{d}x+(\gamma-1)\int_{\Omega}\dive u\mathcal{Q}(\nabla u)\mathrm{d}x-\kappa(\gamma-1)\int_{\Omega}\nabla\dive u\cdot\nabla\theta\mathrm{d}x\\
		\nonumber& -(\gamma-1)\int_{\Omega}P(\dive u)^{2}\mathrm{d}x+c_{v}(\gamma-1)\int_{\Omega}\rho\theta u\cdot\nabla\dive u\mathrm{d}x.
	\end{align}
Multiplying \eqref{b} with $u_{t}$, integrating over $\Omega$, and using the above identity, it follows
	\begin{eqnarray}
		\nonumber&&\frac{\mathrm{d}}{\mathrm{d}t}\left(\frac{\mu}{2}\norm{\nabla u}^{2}+\frac{\mu+\lambda}{2}\norm{\dive u}^{2}-\int_{\Omega}P\dive u\mathrm{d}x \right) +\norm{\sqrt{\rho}u_{t}}^{2}\\
		\nonumber &=&-\int_{\Omega}\rho(u\cdot\nabla)u\cdot u_{t}\mathrm{d}x -(\gamma-1)\int_{\Omega}\dive u\mathcal{Q}(\nabla u)\mathrm{d}x+\kappa(\gamma-1)\int_{\Omega}\nabla\dive u\cdot\nabla\theta\mathrm{d}x\\
	\label{d4}&&+(\gamma-1)\int_{\Omega}P(\dive u)^{2}\mathrm{d}x-c_{v}(\gamma-1)\int_{\Omega}\rho\theta u\cdot\nabla\dive u\mathrm{d}x:=\sum\limits_{i=1}^{5}J_{i}.
	\end{eqnarray}
By Corollary \ref{d1}, it follows from the Gagliardo-Nirenberg, Sobolev, Poincar\'e, and Young inequalities that
\begin{eqnarray*}
	|J_{1}|&\leq&\norm{\rho}_{\infty}^{\frac{1}{2}}\norm{u}_{6}\norm{\nabla u}_{3}\norm{\sqrt{\rho}u_{t}}\leq C\norm{\nabla u}^{\frac{3}{2}}\norm{\nabla^{2}u}^{\frac{1}{2}}\norm{\sqrt{\rho}u_{t}}\\
&\leq&\frac18\|\sqrt\rho u_t\|^2+\eta\|\nabla^2u\|^2+C_\eta\|\nabla u\|^6,\\
	|J_{2}|&\leq& C\int_{\Omega}\nn{\nabla u}^{3}\mathrm{d}x\leq C\norm{\nabla u}^{\frac{3}{2}}\norm{\nabla^{2} u}^{\frac{3}{2}}
\leq\eta\|\nabla^2u\|^2+C_\eta\|\nabla u\|^6,\\
	|J_{3}|&\leq& C\norm{\nabla^{2}u}\norm{\nabla\theta}\leq \eta\norm{\nabla^{2}u}^2+C_\eta\norm{\nabla\theta}^2,\\
	|J_{4}|&\leq& C \norm{\rho}_{\infty}\norm{\theta}_{6}\norm{\nabla u}\norm{\nabla u}_{3}\leq C\norm{\nabla\theta}\norm{\nabla u}^{\frac{3}{2}}\norm{\nabla^{2}u}^{\frac{1}{2}}\\
&\leq&\eta\|\nabla^2u\|_2^2+C_\eta(\|\nabla\theta\|^2+\|\nabla u\|^6),\\
	|J_{5}|&\leq& C\norm{\rho}_{\infty}^{\frac{1}{2}}\norm{u}_{\infty}\norm{\sqrt{\rho}\theta}_{2}\norm{\nabla^{2}u}\leq C\norm{\nabla u}^{\fr}\norm{\nabla^{2}u}^{\frac{3}{2}}\norm{\sqrt{\rho}\theta}\nonumber\\
&\leq&\eta\|\nabla^2u\|^2+C_\eta\|\nabla u\|^2\|\sqrt\rho\theta\|^4,
\end{eqnarray*}
for any positive number $\eta$.
Plugging the above estimates into \eqref{d4}, adding the resultant with \eqref{o} multiplied with a small positive number $\epsilon_{1}$, and choosing $\eta$ sufficiently small, one obtains
\begin{eqnarray*}
\frac{\mathrm{d}}{\mathrm{d}t}\left(\frac{\mu}{2}\norm{\nabla u}^{2}+\frac{\mu+\lambda}{2}\norm{\dive u}^{2}-\int_{\Omega}P\dive u\mathrm{d}x \right) +\norm{\sqrt{\rho}u_{t}}^{2}+\frac{\epsilon_{1}}{2}\norm{\nabla^{2}u}^{2}\\
 \leq \ C\bigg(\norm{\nabla u}^{6}+ \norm{\nabla\theta}^{2}+\|\nabla u\|^2\|\sqrt\rho\theta\|^4\bigg).
\end{eqnarray*}
Integrating the above inequality over $(0,T)$ and using Proposition \ref{p}, one deduces
\begin{eqnarray*}
&&\mu\sup_{0\leq t\leq T}\norm{\nabla u}^{2}+\int_{0}^{T}\left(\norm{\sqrt{\rho}u_{t}}^{2}+\epsilon_{1}\norm{\nabla^{2}u}^{2} \right)\mm\\
	\nonumber
	&\leq&2\sup_{0\leq t\leq T}\left|\int_\Omega P\text{div} u dx\right|+C\left(\|\nabla u_0\|^2+\left|\int_\Omega P_0\text{div}u_0dx\right|\right)\\
&&+C\int_0^T\left(\|\nabla u\|^6+\|\nabla\theta\|^2+\|\nabla u\|^2\|\sqrt\rho\theta\|^4\right)dt\\
&\leq&\frac\mu2\sup_{0\leq t\leq T}\|\nabla u\|^2+C\Big[1+\sup_{0\leq t\leq T}\left(\|\rho\|_\infty\|\sqrt\rho\theta\|^2\right)+T\Phi^3(T)\Big]\\
&\leq&\frac\mu2\sup_{0\leq t\leq T}\|\nabla u\|^2+C\Big(1+T\Phi^3(T)\Big),
\end{eqnarray*}
from which, by (\ref{dj}), the conclusion follows.
\end{proof}

\begin{proposition}\label{d6}
	Under the conditions of Proposition \ref{p}, it holds that
$$
  \sup_{0\leq t\leq T}\norm{(\sqrt t\nabla\theta,\sqrt{t}\sqrt{\rho}u_{t},\sqrt t\nabla^2u)}^{2}+ \int_{0}^{T}\norm{(\sqrt{t}\nabla u_{t},\sqrt t\nabla^2\theta)}^{2}\mm
  \leq C .
$$
\end{proposition}

\begin{proof}
By Corollary \ref{d1}, it follows from \eqref{c} and the Sobolev and Young inequalities that
\begin{eqnarray*}
	\norm{\nabla^{2}\theta}^{2}
	\nonumber &\leq& C\left(\norm{\rho\theta_{t}}^{2}+\norm{\rho(u\cdot\nabla)\theta}^{2}+\norm{\dive uP}^{2}+\norm{\mathcal{Q}(\nabla u)}^{2} \right)\\\nonumber
	&\leq& C\left(\norm{\rho}_{\infty}\norm{\sqrt\rho\theta_{t}}^{2}+\norm{\rho}_{\infty}^{2}\norm{u}_{6}^{2}\norm{\nabla\theta}_{3}^{2}+\norm{\rho}_{\infty}^{2}\norm{\theta}_{6}^{2}\norm{\nabla u}^{2}_{3}+\norm{\nabla u}_{4}^{4} \right) \\
	\nonumber
	&\leq& C\left(\norm{\sqrt\rho\theta_{t}}^{2}+\norm{\nabla u}^{2}\norm{\nabla\theta}\norm{\nabla^{2}\theta}+\norm{\nabla\theta}^{2}\norm{\nabla u}\norm{\nabla^{2}u}+\norm{\nabla u}\norm{\nabla^{2}u}^{3} \right) \\\nonumber
	&\leq& \frac{1}{2}\norm{\nabla^{2}\theta}^{2}+C\big(\norm{\sqrt{\rho}\theta_{t}}^{2}+\norm{\nabla u}^{4}\norm{\nabla\theta}^{2}\\
  &&+\norm{\nabla\theta}^{2}\norm{\nabla u}\norm{\nabla^{2}u}+\norm{\nabla u}\norm{\nabla^{2}u}^{3} \big)
\end{eqnarray*}
and, thus,
\begin{align}
	\label{i}\norm{\nabla^{2}\theta}^{2}\leq C\left(\norm{\sqrt{\rho}\theta_{t}}^{2}+\norm{\nabla u}^{4}\norm{\nabla\theta}^{2}+\norm{\nabla\theta}^{2}\norm{\nabla u}\norm{\nabla^{2}u}+\norm{\nabla u}\norm{\nabla^{2}u}^{3} \right).
\end{align}
Testing \eqref{c} with $\theta_{t}$ yields
\begin{align}	\frac{\kappa}{2}\frac{\mathrm{d}}{\mathrm{d}t}\norm{\nabla\theta}^{2}+c_{v}\norm{\sqrt{\rho}\theta_{t}}^{2}
\label{j}	=\int_{\Omega}\Big[-c_{v}\rho(u\cdot\nabla)\theta\theta_{t}-P\dive u\theta_{t}+\mathcal{Q}(\nabla u)\theta_{t}\Big]\mathrm{d}x.
\end{align}
Terms on the right-hand side of (\ref{j}) are estimated by Proposition \ref{p} and the Gagliardo-Nirenberg, Poincar\'e, and Young inequalities as follows
\begin{eqnarray}
\nonumber\left|\int_{\Omega}\rho(u\cdot\nabla)\theta\theta_{t}\mathrm{d}x\right|&\leq& \norm{\rho}_{\infty}^{\frac{1}{2}}\norm{u}_{\infty}\norm{\nabla\theta}\norm{\sqrt{\rho}\theta_{t}}\leq C\norm{\rho}_{\infty}^{\frac{1}{2}}\norm{\nabla u}^{\frac{1}{2}}\norm{\nabla^{2}u}^{\frac{1}{2}}\norm{\nabla\theta}\norm{\sqrt{\rho}\theta_{t}}\\
&\leq& \frac{1}{8}\norm{\sqrt{\rho}\theta_{t}}^{2}+C\norm{\nabla u}\norm{\nabla^{2} u}\norm{\nabla \theta}^{2},\label{k}\\
\nonumber\left|\int_{\Omega}P\dive u\theta_{t}\mathrm{d}x\right|&\leq& R\int_{\Omega}\rho|\theta||\nabla u||\theta_{t}|\mathrm{d}x \leq C\norm{\rho}_{\infty}^{\frac{1}{2}}\norm{\theta}_{6}\norm{\nabla u}_{3}\norm{\sqrt{\rho}\theta_{t}}\\
&\leq& C\norm{\rho}_{\infty}^{\frac{1}{2}}\norm{\nabla\theta}\norm{\nabla u}^{\frac{1}{2}}\norm{\nabla^{2} u}^{\frac{1}{2}}\norm{\sqrt{\rho}\theta_{t}} \nonumber\\
&\leq& \frac{c_{v}}{8}\norm{\sqrt{\rho}\theta_{t}}^{2}+C\norm{\nabla u}\norm{\nabla^{2} u}\norm{\nabla \theta}^{2},\\
\nonumber\int_{\Omega}\mathcal{Q}(\nabla u)\theta_{t}\mathrm{d}x&=&\frac{\mathrm{d}}{\mathrm{d}t}\int_{\Omega}\mathcal{Q}(\nabla u)\theta\mathrm{d}x-\int_{\Omega}\left(4\mu Du:Du_{t}+2\lambda\dive u\dive u_{t} \right) \theta\mathrm{d}x\\
\nonumber &\leq &\frac{\mathrm{d}}{\mathrm{d}t}\int_{\Omega}\mathcal{Q}(\nabla u)\theta\mathrm{d}x+C\norm{\nabla u}_{3}\norm{\nabla u_{t}}\norm{\theta}_{6} \\
&\leq& \frac{\mathrm{d}}{\mathrm{d}t}\int_{\Omega}\mathcal{Q}(\nabla u)\theta\mathrm{d}x+C\norm{\nabla u}^{\frac{1}{2}}\norm{\nabla^{2} u}^{\frac{1}{2}}\norm{\nabla u_{t}}\norm{\nabla\theta}.\label{m}
\end{eqnarray}
Plugging \eqref{k}--\eqref{m} into \eqref{j} and adding the resultant with \eqref{i} multiplied with
a small positive number $\epsilon_{2}$, one obtains
\begin{align*}
\frac{\mathrm{d}}{\mathrm{d}t}\Big( \frac{\kappa}{2}\norm{\nabla\theta}^{2}&-\int_{\Omega}\mathcal{Q}(\nabla u)\theta\mathrm{d}x\Big)+\frac{c_{v}}{2}\norm{\sqrt{\rho}\theta_{t}}^{2}+\epsilon_{2}\norm{\nabla^{2}\theta}^{2} \\
\leq& \ \ C\Big(\norm{\nabla u}^{4}\norm{\nabla\theta}^{2}+\norm{\nabla\theta}^{2}\norm{\nabla u}\norm{\nabla^{2}u}+ \norm{\nabla u}\norm{\nabla^{2}u}^{3}\\
&\ \ +\norm{\nabla u}^{\frac{1}{2}}\norm{\nabla^{2}u}^{\frac{1}{2}}\norm{\nabla u_{t}}\norm{\nabla\theta}\Big).
\end{align*}
Multiplying the above inequality with $t$ yields
\begin{align}
\nonumber\frac{\mathrm{d}}{\mathrm{d}t}\left( \frac{\kappa}{2}\right.\norm{\sqrt{t}&\nabla\theta}^{2}\left.-t\int_{\Omega}\mathcal{Q}(\nabla u)\theta\mathrm{d}x\right)+\frac{c_{v}}{2}\norm{\sqrt{t}\sqrt{\rho}\theta_{t}}^{2}+\epsilon_{2}\norm{\sqrt{t}\nabla^{2}\theta}^{2} +\int_{\Omega}\mathcal{Q}(\nabla u)\theta\mathrm{d}x\\
\nonumber\label{n}  \leq&\  \ C(\norm{\nabla u}^{4}\norm{\sqrt{t}\nabla\theta}^{2}+\norm{\sqrt{t}\nabla\theta}^{2}\norm{\nabla u}\norm{\nabla^{2}u} +\sqrt{t}\norm{\nabla u}\norm{\nabla^{2}u}^{2}\norm{\sqrt{t}\nabla^{2}u}\\
&\ +\norm{\nabla u}^{\frac{1}{2}}\norm{\nabla^{2}u}^{\frac{1}{2}}\norm{\sqrt{t}\nabla u_{t}}\norm{\sqrt{t}\nabla\theta}+\norm{\nabla\theta}^{2}).
\end{align}
It follows from the Gagliardo-Nirenberg inequality and (\ref{dj}) that
\begin{eqnarray}
\nonumber t\int_{\Omega}\mathcal{Q}(\nabla u)\theta\mathrm{d}x&\leq& Ct\int_{\Omega}|\nabla u|^{2}|\theta|\mathrm{d}x
 \leq Ct\norm{\nabla u}\norm{\nabla u}_{3}\norm{\theta}_{6}\\
&\leq& CT^{\frac{1}{4}}\norm{\nabla u}^{\frac{3}{2}}\norm{\sqrt{t}\nabla^{2}u}^{\frac{1}{2}}\norm{\sqrt{t}\nabla\theta}\leq C[T\Phi^3(T)]^\frac14\|\sqrt t\nabla\theta\|\|\sqrt t\nabla^2u\|^\frac12\nonumber\\
&\leq& \eta\sup_{0\leq t\leq T}\|\sqrt t\nabla\theta\|^2+C_\eta\sup_{0\leq t\leq T}\|\sqrt t\nabla^2u\|^2,\label{d3}
\end{eqnarray}
for any positive $\eta>0$.
Integrating \eqref{n} over $(0,T)$ and using \eqref{d3}, one deduces by Proposition \ref{p} and the Young inequality that
\begin{eqnarray*}
\nonumber&&\kappa\sup_{0\leq t\leq T}\norm{\sqrt{t}\nabla\theta}^{2}+\int_{0}^{T}\left(c_{v} \norm{\sqrt{t}\sqrt{\rho}\theta_{t}}^{2}+2\epsilon_{2}\norm{\sqrt{t}\nabla^{2}\theta}^{2}\right)\mm  \\
\nonumber&\leq& \sup_{0\leq t\leq T}\left(t\int_{\Omega}\mathcal{Q}(\nabla u)\theta\mathrm{d}x\right)+C\int_{0}^{T}\Big(\norm{\nabla u}^{4}\norm{\sqrt{t}\nabla\theta}^{2}+\norm{\sqrt{t}\nabla\theta}^{2}\norm{\nabla u}\norm{\nabla^{2}u}\\
\nonumber&&+ \sqrt{t}\norm{\nabla u}\norm{\nabla^{2}u}^{2}\norm{\sqrt{t}\nabla^{2}u}
+\norm{\nabla u}^{\frac{1}{2}}\norm{\nabla^{2}u}^{\frac{1}{2}}\norm{\sqrt{t}\nabla u_{t}}\norm{\sqrt{t}\nabla\theta}+\norm{\nabla\theta}^{2}\Big)\mm\\
&\leq&\eta\sup_{0\leq t\leq T}\|\sqrt t\nabla\theta\|^2+C_\eta \left(1+T^\frac12\Phi^\frac32(T)\right)\sup_{0\leq t\leq T}\|\sqrt t\nabla^2u\|\nonumber\\
&&+C_\eta\left(T\Phi^{3}(T)+T^{\frac{1}{2}}\Phi^{2}(T)+T^{\frac{1}{4}}\Phi^{\frac{3}{2}}(T)+1\right),
\end{eqnarray*}
for any $\eta>0$,
from which, choosing $\eta$ sufficiently small and by \eqref{dj}--\eqref{dl}, one gets
\begin{equation}
  \sup_{0\leq t\leq T}\norm{\sqrt{t}\nabla\theta}^{2}+\int_{0}^{T}\left(\norm{\sqrt{t}\sqrt{\rho}\theta_{t}}^{2}+\norm{\sqrt{t}\nabla^{2}\theta}^{2} \right)\mm\leq C\left(\sup_{0\leq t\leq T}\|\sqrt t\nabla^2u\|+1\right). \label{ESTNT}
\end{equation}

Differentiating \eqref{b} with respect to $t$ yields
\begin{align}
\rho(u_{tt}+u\cdot\nabla u_{t})&+\rho u_{t}\cdot\nabla u+\rho_{t}(u_{t}+u\cdot\nabla u)\nonumber\\
&-\mu\Delta u_{t}-(\mu+\lambda)\nabla\dive u_{t}+\nabla P_{t}=0.\label{r}
\end{align}
It follows from (\ref{EQP}) that
\begin{eqnarray}
	\nonumber&&\int_{\Omega}\nabla P_{t}u_{t}\mathrm{d}x=-\int_{\Omega}P_{t}\dive u_{t}\mathrm{d}x\\\nonumber
	&=&(\gamma-1)\int_{\Omega}\dive u_{t}\left(-\mathcal{Q}(\nabla u)-\kappa\Delta\theta+P\dive u+c_{v}\dive(\rho u\theta) \right)\mathrm{d}x,
\end{eqnarray}
Testing \eqref{r} by $u_{t}$ and utilizing the above equality yield
\begin{eqnarray*}
&&\frac{1}{2}\frac{\mathrm{d}}{\mathrm{d}t}\norm{\sqrt{\rho}u_{t}}^{2}+\mu\norm{\nabla u_{t}}^{2}
+(\mu+\lambda)\norm{\dive u_{t}}^{2}\\
&=&(\gamma-1)\int_{\Omega}\dive u_{t}\Big(\mathcal{Q}(\nabla u)+\kappa\Delta\theta-P\dive u-c_{v}
\dive(\rho u\theta)\Big)\mathrm{d}x\\
&&+\int_{\Omega}\dive(\rho u)\big(u_{t}+(u\cdot\nabla)u\big)
\cdot u_{t}\mathrm{d}x-\int_{\Omega}\rho(u_{t}\cdot\nabla)u\cdot u_{t}\mathrm{d}x.
\end{eqnarray*}
Multiplying the above identity with $t$ and integrating over $(0,T)$ lead to
\begin{eqnarray}
\nonumber&& \frac12\sup_{0\leq t\leq T}\norm{\sqrt{t}\sqrt{\rho}u_{t}}^{2}+\mu\int_{0}^{T}
\norm{\sqrt{t}\nabla u_{t}}^{2}\mm\\
\nonumber&\leq&\frac12\int_{0}^{T}\norm{\sqrt{\rho}u_{t}}^{2}\mm+C\int_{0}^{T}t\int_{\Omega}|\nabla u_{t}
||\nabla
u|^{2}\mathrm{d}x\mathrm{d}t+\kappa(\gamma-1)\int_{0}^{T}t\int_{\Omega}\text{div}u_t\Delta\theta\mathrm{d}x\mathrm{d}t \\
\nonumber&&-(\gamma-1)\int_{0}^{T}t\int_{\Omega}\text{div}u_t
P\text{div}u\mathrm{d}x\mathrm{d}t-c_v(\gamma-1)\int_{0}^{T}t\int_{\Omega}
\text{div}u_t\text{div}(\rho u\theta)\mathrm{d}x\mathrm{d}t\\
\nonumber&&-\int_{0}^{T}t\int_{\Omega}\rho u\cdot\nabla|u_{t}|^2 \mathrm{d}x\mathrm{d}t-\int_{0}^{T}t\int_{\Omega}\rho  u
\cdot\nabla((u\cdot\nabla) u\cdot u_{t})\mathrm{d}x\mathrm{d}t
\\
\label{d5} && -\int_{0}^{T}t\int_{\Omega}\rho(u_t\cdot\nabla)u\cdot u_t\mathrm{d}x\mathrm{d}t =\sum\limits_{i=1}^{8}K_{i}.
\end{eqnarray}
By the Cauchy and Gagliardo-Nirenberg inequalities, it follows from Propositions \ref{p}--\ref{d6} and Corollary \ref{d1} that $K_{1}\leq C$,
\begin{align*}
K_{2}\leq&  C\int_{0}^{T}t\norm{\nabla u_{t}}\norm{\nabla u}_{4}^{2}\mm
\leq C\int_{0}^{T}\norm{\sqrt{t}\nabla u_{t}}\norm{\nabla u}^{\frac{1}{2}}\norm{\sqrt{t}\nabla^{2}u}\norm{\nabla^{2}
u}^{\frac{1}{2}}\mm\\
\nonumber\leq&  CT^\frac14\Phi(T)\sup_{0\leq t\leq T}
\|\sqrt t\nabla^2u\|,\\
K_{3}\leq&   \frac{\mu}{2}\int_{0}^{T}\norm{\sqrt{t}\nabla u_{t}}^{2}\mm+C\int_0^T\|\sqrt t\nabla^2\theta\|^2dt,\\
K_{4}\leq&  C\int_{0}^{T}t\norm{\rho}_{\infty}\norm{\theta}_{6}\norm{\nabla u}_{3}\norm{\nabla u_{t}}\mm\\
\leq& C\int_{0}^{T}\norm{\sqrt{t}\nabla\theta}\norm{\nabla
u}^{\frac{1}{2}}\norm{\nabla^{2}u}^{\frac{1}{2}}\norm{\sqrt{t}\nabla u_{t}}\mm\leq
CT^{\frac{1}{4}}\Phi^{\frac{3}{2}}(T),\\
K_{5}\leq&  C\int_{0}^{T}t\int_{\Omega}|\nabla u_{t}|\left(|\nabla\rho||u||\theta|+\rho|\nabla
u||\theta|+\rho|u||\nabla\theta| \right)\mathrm{d}x\mm \\
\leq& C\int_{0}^{T}t\norm{\nabla
u_{t}}\left(\norm{\nabla\rho}_{3}\norm{u}_{\infty}\norm{\theta}_{6}+\norm{\rho}_{\infty}\norm{\nabla
u}_{3}\norm{\theta}_{6}+\norm{\rho}_{\infty}\norm{u}_{6 }\norm{\nabla\theta}_{3} \right)\mm\\
\leq&  C\int_0^T(\norm{\sqrt{t}\nabla u_{t}}\norm{\nabla
u}^{\frac{1}{2}}\norm{\nabla^{2}u}^{\frac{1}{2}}\norm{\sqrt{t}\nabla\theta}dt \\
&+C\int_{0}^{T} \norm{\sqrt{t}\nabla u_{t}}\norm{\nabla
u}\norm{\sqrt{t}\nabla\theta}^{\frac{1}{2}}\norm{\sqrt{t}\nabla^{2}\theta}^{\frac{1}{2}}dt
\leq CT^{\frac{1}{4}}\Phi^{\frac{3}{2}}(T),\\
K_{6}\leq& C \int_{0}^{T}t\norm{\rho}_{\infty}^{\frac{1}{2}}\norm{u}_{6}\norm{\sqrt{\rho} u_{t}}_{3}\norm{\nabla
u_{t}}\mm
\leq C\int_{0}^{T}t\norm{\rho}_{\infty}^{\frac{3}{4}}\norm{\nabla u}\norm{\sqrt{\rho}u_{t}}^{\frac{1}{2}}\norm{\nabla
u_{t}}^{\frac{3}{2}}\mm\\
\leq& C\int_{0}^{T}t^\frac14\norm{\nabla u}\norm{\sqrt{\rho}u_{t}}^{\frac{1}{2}}\norm{\sqrt{t}\nabla
u_{t}}^{\frac{3}{2}}\mm
\leq  CT^{\frac{1}{4}}\Phi^{\frac{3}{2}}(T),\\
K_{7}\leq& C \int_{0}^{T}t\int_{\Omega}\rho|u|\left(|\nabla u|^{2}|u_{t}|+|u||\nabla^{2}u||u_{t}|+|u||\nabla u||\nabla
u_{t}| \right) \mm\\
\nonumber\leq&  C\int_{0}^{T}t\norm{\rho}_{\infty}\left(\norm{u}_{6}\norm{\nabla
u}_{3}^{2}\norm{u_{t}}_{6}+\norm{u}_{6}^{2}\norm{\nabla^{2}u}\norm{u_{t}}_{6}+\norm{u}_{6}^{2}\norm{\nabla
u}_{6}\norm{\nabla u_{t}} \right)\mm\\
\leq& C\int_{0}^{T}\sqrt{t}\norm{\nabla u}^{2}\norm{\nabla^{2}u}\norm{\sqrt{t}\nabla u_{t}}\mm\leq
CT^{\frac{1}{2}}\Phi^{2}(T),\\
K_{8}\leq& C \int_{0}^{T}t\norm{\rho}_{\infty}^{\frac{1}{2}}\norm{\nabla u}\norm{\sqrt{\rho}u_{t}}_{3}\norm{
u_{t}}_{6}\mm\leq C\int_{0}^{T}t\norm{\rho}_{\infty}^{\frac{3}{4}}\norm{\nabla
u}\norm{\sqrt{\rho}u_{t}}^{\frac{1}{2}}\norm{\nabla u_{t}}^{\frac{3}{2}}\mm\\
	\leq& C\int_{0}^{T}t^\frac14\norm{\nabla u}\norm{\sqrt{\rho}u_{t}}^{\frac{1}{2}}\norm{\sqrt{t}\nabla
u_{t}}^{\frac{3}{2}}\mm\leq CT^{\frac{1}{4}}\Phi^{\frac{3}{2}}(T).
\end{align*}
Plugging the above estimates into \eqref{d5}, it follows from (\ref{dk}), (\ref{dl}), and (\ref{ESTNT})
that
\begin{align*}
&\sup_{0\leq t\leq T}\norm{\sqrt{t}\sqrt{\rho}u_{t}}^{2}+\mu\int_{0}^{T}\norm{\sqrt{t}\nabla u_{t}}^{2}\mm\nonumber\\
\leq&CT^\frac14\Phi(T)\sup_{0\leq t\leq T}
\|\sqrt t\nabla^2u\|+C \int_0^T\|\sqrt t\nabla^2\theta\|^2dt +  C \left( 1+T^{\frac{1}{4}}\Phi^{\frac{3}{2}}(T)+T^{\frac{1}{2}}\Phi^{2}(T)\right)\nonumber\\
\leq&C \left(\sup_{0\leq t\leq T}
\|\sqrt t\nabla^2u\|+ \int_0^T\|\sqrt t\nabla^2\theta\|^2dt +1\right)
\leq C \left(\sup_{0\leq t\leq T}
\|\sqrt t\nabla^2u\|+1\right),
\end{align*}
that is
\begin{equation}
  \sup_{0\leq t\leq T}\norm{\sqrt{t}\sqrt{\rho}u_{t}}^{2}+\mu\int_{0}^{T}\norm{\sqrt{t}\nabla u_{t}}^{2}\mm
  \leq C \left(\sup_{0\leq t\leq T}
\|\sqrt t\nabla^2u\|+1\right).\label{ESTRUT}
\end{equation}

Recalling \eqref{o}, it follows from (\ref{ESTNT}), (\ref{ESTRUT}), (\ref{dj}), and the Young inequality that
\begin{eqnarray*}
\nonumber\sup_{0\leq t\leq T}\norm{\sqrt{t}\nabla^{2}u}^{2}&\leq& C\sup_{0\leq t\leq T}\left(\norm{\sqrt{t}\sqrt{\rho}u_{t}}^{2}+t\norm{\nabla u}^{6}+ \norm{\sqrt{t}\nabla \theta}^{2}  \right),\\
\nonumber&\leq& C\sup_{0\leq t\leq T}\left(\norm{\sqrt{t}\sqrt{\rho}u_{t}}^{2}+\norm{\sqrt{t}\nabla \theta}^{2}+T\Phi^3(T) \right)\\
&\leq& C \left(\sup_{0\leq t\leq T}
\|\sqrt t\nabla^2u\|+1\right)\leq\frac12\sup_{0\leq t\leq T}
\|\sqrt t\nabla^2u\|^2+C
\end{eqnarray*}
and, thus, $\sup_{0\leq t\leq T}\norm{\sqrt{t}\nabla^{2}u}^{2}\leq C.$ With the aid of this, the conclusion follows easily from (\ref{ESTNT}) and (\ref{ESTRUT}).
\end{proof}

As a direct corollary of Corollary \ref{d1} and Propositions \ref{p}--\ref{d6}, we have the following corollary.

\begin{corollary}\label{x}
Under the assumptions of Proposition \ref{p}, it holds that
$$
\Phi(T)+\sup_{0\leq t\leq T}\|(\sqrt t\nabla^2u,\sqrt t\sqrt\rho u_t)\|^2\leq C,
$$
that is
\begin{eqnarray*}
\nonumber \sup_{0\leq s\leq t} \left( \norm{\rho}_{\infty}+\norm{\nabla\rho}_{q}+\norm{\left(\sqrt{\rho}\theta,\nabla u,\sqrt{t}\nabla\theta,\sqrt t\nabla^2u,\sqrt t\sqrt\rho u_t\right)}^2\right) \\
+\int_{0}^{t} \norm{\left(\sqrt{\rho}u_{t},\nabla^{2}u,\sqrt{t}\nabla u_{t},\nabla\theta,\sqrt{t}\sqrt{\rho}\theta_{t},\sqrt{t}\nabla^{2}\theta\right)}^{2} \mathrm{d}s\leq C.
\end{eqnarray*}
\end{corollary}

\begin{proposition}
\label{prop2.9}
Under the assumptions of Proposition \ref{p}, it holds that
 \begin{align}
	\nonumber\sup_{0\leq t\leq {T_0}} \left\|\left(\rho_{t},t\nabla^{2}\theta,t\sqrt{\rho}\theta_{t}\right)\right\|^{2} +\int_{0}^{{T_0}}\left(\norm{t\nabla\theta_{t}}^{2}+\|t\nabla^2\theta\|_6^2
+\norm{\sqrt{t}\nabla^{2}u}^{2}_{q}\right)\mm\leq C.
\end{align}
\end{proposition}
\begin{proof}
Using\eqref{a}, it follows from the H\"older and Sobolev inequalities and Corollary \ref{x} that
	\begin{eqnarray}
	\norm{\rho_{t}}&=&\norm{u\cdot\nabla\rho+\dive u\rho}
	\leq\norm{u}_{6}\norm{\nabla\rho}_{3}+\norm{\nabla u}\norm{\rho}_{\infty}\nonumber\\
	&\leq& C(\norm{\nabla\rho}_{3}+\norm{\rho}_{\infty})\norm{\nabla u}\leq C.\label{INEQ1}
	\end{eqnarray}
Differentiating \eqref{c} with respect to $t$ yields
	\begin{align}
		\nonumber c_{v}\rho (\theta_{tt}&+u\cdot\nabla\theta_{t})+c_{v}\rho_{t}(\theta_{t}+u\cdot\nabla\theta)+c_{v}\rho u_{t}\cdot\nabla\theta-\kappa\Delta\theta_{t}\\
		\nonumber&+P_{t}\dive u+P\dive u_{t}= 4\mu Du:Du_{t}+2\lambda\dive u\dive u_{t}.
	\end{align}
Multiply the above equality with $\theta_{t}$ and integrating over $\Omega$, one gets
\begin{eqnarray}
&&\frac{c_{v}}{2}\frac{\mathrm{d}}{\mathrm{d}t}\norm{\sqrt{\rho}\theta_{t}}^{2}+\kappa\norm{\nabla\theta_{t}}^{2}\nonumber\\
&=&-c_{v}\int_{\Omega}\rho_{t}|\theta_{t}|^{2}\mathrm{d}x-c_{v}\int_{\Omega}\rho_{t} u\cdot\nabla\theta \theta_{t}\mathrm{d}x-c_{v}\int_{\Omega}\rho u_{t}\cdot\nabla \theta \theta_{t}\mathrm{d}x-\int_{\Omega}P_{t}\dive u\theta_{t}\mathrm{d}x\nonumber \\
&&-\int_{\Omega}P\dive u_{t}\theta_{t}\mathrm{d}x+\int_{\Omega}\left(4\mu Du:Du_{t}+2\lambda\dive u\dive u_{t} \right)\theta_{t}\mathrm{d}x:=\sum_{i=1}^{6}L_{i}.\label{d7}
\end{eqnarray}
Terms on the right-hand side are estimated by \eqref{a}, Corollary \ref{x}, (\ref{INEQ1}), and the H\"{o}lder, Gagliardo-Nerenberg, and Young inequalities as follows
\begin{eqnarray*}
|L_{1}|&=&c_v\left|\int_\Omega\text{div}(\rho u)|\theta_t|^2dx\right|\leq c_{v}\int_{\Omega}\rho|u||\theta_{t}||\nabla\theta_{t}|\mathrm{d}x\\
&\leq& C\norm{\rho}_{\infty}^{\frac{1}{2}}\norm{u}_{6}\norm{\sqrt{\rho}\theta_{t}}_{3}\norm{\nabla\theta_{t}} \leq  C\norm{\rho}_{\infty}^{\frac{3}{4}}\norm{\nabla u}\norm{\sqrt{\rho}\theta_{t}}^{\frac{1}{2}}\norm{\nabla\theta_{t}}^{\frac{3}{2}}\\
&\leq&\frac{\kappa}{8}\norm{\nabla\theta_{t}}^{2}+C \norm{\sqrt{\rho}\theta_{t}}^{2},\\
|L_{2}|&=&c_v\left|\int_\Omega\text{div}(\rho u) u\cdot\nabla\theta\theta_tdx\right|\\
&\leq&  c_{v}\int_{\Omega}\rho|u|\left(|\nabla u||\nabla\theta||\theta_{t}|+|u||\nabla^{2}\theta||\theta_{t}|+|u||\nabla\theta||\nabla\theta_{t}| \right) \mathrm{d}x\\
	&\leq& C\norm{\rho}_{\infty}\left(\norm{u}_{\infty}\norm{\nabla u}_{3}\norm{\nabla\theta}\norm{\theta_{t}}_{6}+\norm{u}_{6}^{2}\norm{\nabla^{2}\theta}\norm{\theta_{t}}_{6}+\norm{u}_{6}^{2}\norm{\nabla\theta}_{6}\norm{\nabla\theta_{t}} \right)\\
	&\leq& C\norm{\rho}_{\infty}\left( \norm{\nabla u}\norm{\nabla^{2}u}\norm{\nabla\theta}\norm{\nabla\theta_{t}}+\norm{\nabla u}^{2}\norm{\nabla^{2}\theta}\norm{\nabla\theta_{t}}\right) \\
	&\leq&\frac{\kappa}{8}\norm{\nabla\theta_{t}}^{2}+C\left(\norm{\nabla^{2}u }^{2}\norm{\nabla\theta}^{2}+\norm{\nabla^{2}\theta}^{2} \right),\\
	 |L_{3}|&\leq& c_{v}\int_{\Omega}\rho|u_{t}||\nabla\theta||\theta_{t}|\mathrm{d}x \leq C\norm{\rho}^{\frac{1}{2}}_\infty\norm{\sqrt{\rho}u_{t}}_{3}\norm{\nabla\theta}\norm{\theta_{t}}_{6}\\
	&\leq& C\norm{\rho}^{\frac{3}{4}}_\infty\norm{\sqrt{\rho}u_{t}}^{\frac{1}{2}}\norm{\nabla u_{t}}^{\frac{1}{2}}\norm{\nabla \theta}\norm{\nabla\theta_{t}}
 \leq \frac{\kappa}{8}\norm{\nabla\theta_{t}}^{2}+C\norm{\sqrt{\rho}u_{t}}\norm{\nabla u_{t}}\norm{\nabla\theta}^{2},\\
	|L_{4}|&\leq& R\int_{\Omega}\left(|\rho_{t}|\theta|\nabla u||\theta_{t}|+\rho|\theta_{t}|^{2}|\nabla u| \right)\mathrm{d}x\\
	&\leq& C\left( \norm{\rho_{t}}\norm{\theta}_{6}\norm{\nabla u}_{6}\norm{\theta_{t}}_{6}+\norm{\rho}_{\infty}^{\frac{1}{2}}\norm{\nabla u}\norm{\sqrt{\rho}\theta_{t}}_{3}\norm{\theta_{t}}_{6}\right) \\
	&\leq &C\left(  \norm{\nabla\theta}\norm{\nabla^{2}u}\norm{\nabla\theta_{t}}+ \norm{\sqrt{\rho}\theta_{t}}^{\frac{1}{2}}\norm{\nabla\theta_{t}}^{\frac{3}{2}}\right)\\
	&\leq&\frac{\kappa}{8}\norm{\nabla\theta_{t}}^{2}+C\left(  \norm{\nabla\theta}^{2}\norm{\nabla^{2}u}^{2}+ \norm{\sqrt{\rho}\theta_{t}}^{2}\right) ,\\
	|L_{5}|&\leq& R\int_{\Omega}\rho \theta |\nabla u_{t}||\theta_{t}|\mathrm{d}x\leq C\norm{\rho}_{\infty}^{\frac{1}{2}}\norm{\sqrt{\rho}\theta}_{3}\norm{\nabla u_{t}}\norm{\theta_{t}}_{6}\\
	&\leq& C\norm{\rho}_{\infty}^{\frac{3}{4}}\norm{\sqrt{\rho}\theta}^{\frac{1}{2}}\norm{\nabla\theta}^{\frac{1}{2}}\norm{\nabla u_{t}}\norm{\nabla\theta_{t}} \leq \frac{\kappa}{8}\norm{\nabla\theta_{t}}^{2}+C \norm{\nabla\theta}\norm{\nabla u_{t}}^{2},\\
	|L_{6}|&\leq& C\int_{\Omega}|\nabla u||\nabla u_{t}||\theta_{t}|\mathrm{d}x\leq C\norm{\nabla u}_{3}\norm{\nabla u_{t}}\norm{\theta_{t}}_{6}\\
	\nonumber&\leq& C\norm{\nabla u}^{\frac{1}{2}}\norm{\nabla ^{2}u}^{\frac{1}{2}}\norm{\nabla u_{t}}\norm{\nabla\theta_{t}} \leq \frac{\kappa}{8}\norm{\nabla\theta_{t}}^{2}+C \norm{\nabla^{2}u}\norm{\nabla u_{t}}^{2}.
	\end{eqnarray*}
Substituting the above estimates into \eqref{d7} yields
	\begin{eqnarray*}
	&&{\frac{c_{v}}{2}}\frac{\mathrm{d}}{\mm}\norm{\sqrt{\rho}\theta_{t}}^{2}+\frac{\kappa}{4}\norm{\nabla\theta_{t}}^{2}\\
	& \leq &C\left(  \norm{\sqrt{\rho}\theta_{t}}^{2}+ \norm{\nabla^{2}u}^{2}\norm{\nabla\theta}^{2}+ \norm{\nabla^{2}\theta}^{2}+\norm{\sqrt{\rho}u_{t}}\norm{\nabla u_{t}}\norm{\nabla\theta}^{2}\right.\\
	&&+ \left.\norm{\nabla\theta}\norm{\nabla u_{t}}^{2}+ \norm{\nabla^{2}u}\norm{\nabla u_{t}}^{2} \right).
	\end{eqnarray*}
Multiplying the above inequality by $t^{2}$ and by Corollary \ref{x}, it follows
\begin{eqnarray}
\frac{c_{v}}{2}\frac{\mathrm{d}}{\mm}\norm{t\sqrt{\rho}\theta_{t}}^{2}+\frac{\kappa}{4}\norm{t\nabla\theta_{t}}^{2}
\leq C\left( \norm{\sqrt{t}\sqrt{\rho}\theta_{t}}^{2}+\norm{\sqrt{t}\nabla^{2}\theta}^{2}+
	 \norm{\sqrt{t}\nabla u_{t}} +1\right).\label{t}
\end{eqnarray}
Integrating \eqref{t} over $(0,T)$ and using Corollary \ref{x} yield
\begin{equation}
\sup_{0\leq t\leq {T_0}} \norm{t\sqrt{\rho}\theta_{t}}^{2}+ \int_{0}^{{T_0}}\norm{t\nabla\theta_{t}}^{2}\mm\leq C.
\label{INEQ2}
\end{equation}
Recalling \eqref{i}, it follows from (\ref{INEQ2}) and Corollary \ref{x} that 	
 \begin{eqnarray}
	\sup_{0\leq t\leq T}\norm{t\nabla^{2}\theta}^{2}
	&	\leq& C\sup_{0\leq t\leq T}\left(\norm{t\sqrt{\rho}\theta_{t}}^{2}+ t\|\nabla u\|^4 \norm{\sqrt{t}\nabla\theta}^{2}+\sqrt t\|\nabla u\|\norm{\sqrt{t}\nabla^{2}u}^{3}\right)\nonumber\\
&&+\sup_{0\leq t\leq T}\sqrt t\norm{\sqrt{t}\nabla\theta}^{2}\|\nabla u\|\norm{\sqrt{t}\nabla^{2}u} \leq  C. \label{INEQ3}
	\end{eqnarray}
Applying the elliptic estimates to \eqref{b} and by Corollary \ref{x}, one obtains from the H\"older, Sobolev, and Poincar\'e
inequalities that
	\begin{eqnarray*}
	\norm{\nabla^{2}u }_{q}
	&\leq& C\left(\norm{\rho u_{t}}_{q}+\norm{\rho(u\cdot\nabla)u}_{q}+\norm{\nabla P}_{q} \right)\\
	&\leq &C\left( \norm{\rho}_{\infty}\norm{u_{t}}_{6}+\norm{\rho}_{\infty}\norm{u}_{\infty}\norm{\nabla u}_{6}+\norm{\nabla\rho}_{q}\norm{\theta}_{\infty}+\norm{\rho}_{\infty}\norm{\nabla\theta}_{6}\right)\\
	&\leq& C\left(\norm{\nabla u_{t}}+\norm{\nabla^{2}u}^{2} +\norm{\nabla^{2}\theta} \right)
		\end{eqnarray*}
and, thus, by (\ref{INEQ3}) and Corollary \ref{x}, one gets
\begin{eqnarray}
\int_{0}^{{T_0}}\norm{\sqrt{t}\nabla^{2} u}_{q}^{2}\mm
 &\leq& C\int_{0}^{{T_0}}\left(\norm{\sqrt{t}\nabla u_{t}}^{2}+\norm{\sqrt{t}\nabla^{2}u}^{2}\norm{\nabla^{2}u}^{2}+\norm{\sqrt{t}\nabla^{2}\theta}^{2}  \right)\mm\nonumber\\
 &\leq &C.
 \label{INEQ4}
\end{eqnarray}
Finally, applying the elliptic estimates to \eqref{c} and using the Sobolev and Poincar\'e inequalities, one deduces by Corollary \ref{x} that
 \begin{eqnarray}
 	\nonumber\norm{\nabla^{2}\theta}_{6}^{2}&\leq& C\lr{\norm{\rho\theta_{t}}_{6}^{2}+\norm{\rho\theta\dive u}_{6}^{2}+\norm{\mathcal{Q}(\nabla u)}_{6}^{2}}\\
 	\nonumber&\leq& C\lr{ \norm{\rho}_{\infty}^{2}\norm{\theta_{t}}_{6}^{2}+\norm{\rho}_{\infty}^{2}\norm{\theta}^{2}_{6}\norm{\nabla u}^{2}_{\infty}+\norm{\nabla u}_{\infty}^{2}\norm{\nabla u}_{6}^{2}}\\
 	\nonumber&\leq & C\lr{\norm{\nabla\theta_{t}}^{2}+\norm{\nabla\theta}^{2}\norm{\nabla^{2}u}_{q}^{2}+\norm{\nabla^{2}u}^{2}
 \norm{\nabla^{2}u}_{q}^{2}}.
 	\end{eqnarray}
Hence, it follows from (\ref{INEQ2}), (\ref{INEQ4}), and Corollary \ref{x} that
\begin{eqnarray}
\int_{0}^{{T_0}}\norm{t\nabla^{2}\theta}_{6}^{2}\mm&\leq& C\int_{0}^{{T_0}}\lr{\norm{t\nabla\theta_{t}}^{2}+\norm{\sqrt{t}\nabla\theta}^{2}\norm{\sqrt{t}\nabla^{2}u}_{q}^{2}}\mm
\nonumber\\
 &&+C\int_{0}^{{T_0}}\norm{\sqrt{t}\nabla^{2}u}^{2}\norm{\sqrt{t}\nabla^{2}u}_{q}^{2}\mm\leq C.\label{INEQ5}
\end{eqnarray}
Combining (\ref{INEQ1}) with (\ref{INEQ2})--(\ref{INEQ5}) yields the conclusion.
\end{proof}

As the end of this section, we prove in the next proposition that the existence time ${T_0}$ depends only on $R, \mu, \lambda, c_{v},
\gamma, q$, and the upper bound of $\Phi_{0}$, but is independent of the quantities
$\norm{\nabla^{2}u_0}, \norm{\nabla^{2}\theta_0}, \norm{g_{1}}$, and $\norm{g_{2}}$.
\begin{proposition}\label{v}
Let $q\in(3,6)$ and assume that $(\rho_{0},u_{0},\theta_{0})$ satisfies
\begin{align}
	\nonumber\underline\rho\leq \rho_{0}\in W^{1,q}(\Omega),\quad u_{0}\in H_{0}^{1}(\Omega)\cap H^{2}(\Omega),\quad 0\leq \theta_{0}\in H_{0}^{1}(\Omega)\cap H^{2}(\Omega),
\end{align}
for some positive number $\underline\rho$.

Then, there exist two positive constants $T_0$ and C depending only on $R, \mu, \lambda, c_{v}, \gamma, q,$ and the upper bound of $\Phi_{0}$, such that system \eqref{a}--\eqref{c}, subject to \eqref{d}--\eqref{e}, admits a unique solution $(\rho, u, \theta)$, in $\Omega\times(0,T_0)$, satisfying
$$
\sup_{0\leq t\leq T_0}\left( \|\rho\|_\infty+\norm{\rho}_{W^{1,q}}+\norm{(\rho_{t},\sqrt{\rho}\theta,\nabla u)}^2 \right) +\int_{0}^{T_0} \norm{(\nabla\theta,\sqrt{\rho}u_{t},\nabla^{2}u)}^2 \mm\leq C
$$
and
\begin{eqnarray*}
\sup_{0\leq t\leq T_0}\norm{(\sqrt{t}\nabla\theta,t\sqrt{\rho}\theta_{t},t\nabla^{2}\theta,\sqrt{t}\sqrt{\rho}
u_{t},\sqrt{t}\nabla^{2}u)}^{2}+\int_0^{T_0}\left(\norm{\sqrt{t}\nabla^{2}u}_{q}^{2}+\|t\nabla^2\theta\|_6^2\right)dt \\
 +\int_{0}^{T_0} \norm{(\sqrt{t}\sqrt{\rho}\theta_{t},\sqrt{t}\nabla^{2}\theta,t\nabla\theta_{t},\sqrt{t}\nabla u_{t})}^{2} \mm  \leq C.
\end{eqnarray*}
\end{proposition}

\begin{proof}
By Proposition \ref{f}, there is a unique local strong solution $(\rho, u, \theta)$ on $\Omega\times(0,T_*)$
satisfying the regularities stated
in Proposition \ref{f}. By applying Proposition \ref{f} inductively, one can extend the local solution uniquely to the
maximal time of existence $T_\text{max}$. Then, the following holds
\begin{equation}
  \sup_{T_*\leq t<T_\text{max}}\left(\left\|\frac1\rho\right\|_\infty+\|\rho\|_{W^{1,q}}+\|u\|_{H^2}+\|\theta\|_{H^2}\right)
  =\infty. \label{EXT1}
\end{equation}
For any $T\in(0,T_\text{max})$, $(\rho, u, \theta)$ satisfies the regularities in Proposition \ref{f} with $T_*$ there
replaced by $T$. Let $\epsilon_0$ be the constant stated
in Corollary \ref{d1}, $\Phi$ the function given by (\ref{PHI}), and
set
$$
T_0=\sup\left\{T\in(0,T_\text{max})~\Big|~T^\frac{6-q}{4q}\Phi^2(T)\leq\epsilon_0\right\}.
$$

Claim: $T_0<T_\text{max}$. Assume by contradiction that $T_0=T_\text{max}$. Then, by definition
\begin{equation}
T^\frac{6-q}{4q}\Phi^2(T)\leq\epsilon_0,\quad\forall T\in(0,T_\text{max}).\label{EXT2}
\end{equation}
Since $\Phi(T)\geq1$, it follows from (\ref{EXT2}) that
$T_\text{max}\leq\epsilon_0^\frac{4q}{6-q}.$
Thanks to (\ref{EXT2}), it follows from Corollary \ref{x} and Proposition \ref{prop2.9} that
$$
\Phi(T)+\sup_{0\leq t\leq T}\|(\sqrt t\nabla^2u,t\nabla^2\theta)\|^2\leq C,\quad\forall T\in(0,T_\text{max}),
$$
for a positive constant $C$ independent of $T\in(0,T_\text{max})$. By following the arguments in Proposition \ref{h}, one
deduces by Proposition \ref{g} and (\ref{EXT2}) that
$$
  \inf_{\Omega\times(0, T)}\rho\geq\underline\rho e^{-\int_0^T\|\text{div}u\|_\infty dt}
  \geq\underline\rho e^{-C\epsilon_0},\quad\forall T\in(0,T_\text{max}),
$$
for a positive constant $C$ independent of $T\in(0,T_\text{max})$. Combining the above two yields
$$
  \sup_{T_*\leq t<T_\text{max}}\left(\left\|\frac1\rho\right\|_\infty+\|\rho\|_{W^{1,q}}+\|u\|_{H^2}+\|\theta\|_{H^2}\right)
  \leq C,
$$
which contradicts to (\ref{EXT1}). This contradiction proves the claim.

Since $T_0<T_\text{max}$ and noticing that $\Phi(T)$ is continuous on $[0,T_\text{max})$, one gets by the definition of
$T_0$ that
\begin{equation}
\label{EXT3}
T_0^\frac{6-q}{4q}\Phi^2(T_0)=\epsilon_0.
\end{equation}
Thanks to this and recalling that $\Phi(T)\geq1$, it follows from Corollary \ref{x} that
$T_0\leq\epsilon_0^\frac{4q}{6-q}$ and $\Phi(T_0)\leq C_0$ for a positive constant $C_0$ depending only on $R, \mu, \lambda, c_{v}, \gamma, q,$ and the upper bound of $\Phi_{0}$. Therefore, it follows from (\ref{EXT3}) that
$T_0\geq\left(\frac{\epsilon_0}{C_0^2}\right)^\frac{4q}{6-q}.$ The corresponding estimates follow
from Corollary \ref{x} and Proposition \ref{prop2.9}. This completes the proof.
\end{proof}

\section{A preparing existence result}
In this section, we prove the following existence result, which is a preparation of proving the existence part of Theorem \ref{u}.
Note that the uniqueness is not included here.

\begin{proposition}
\label{existence1}
Assume that all the conditions of Theorem \ref{u} hold. Denote
$$
\Phi_{0}:=\norm{\rho_{0}}_{\infty}+\norm{\nabla\rho_{0}}_{q}+\norm{(\sqrt{\rho_{0}}\theta_{0},\nabla u_{0})}^{2}.
$$

(i) Then, there exists
a positive time ${T_0}$ depending only on
$R, \mu, \lambda, c_{v}, \gamma, q,$ and $\Phi_{0}$,
such that system \eqref{a}--\eqref{c}, subject to \eqref{d}--\eqref{e}, in $\Omega\times(0,{T_0}),$
admits a solution $(\rho, u, \theta)$, which satisfies all the properties stated in Definition \ref{def1} except that the regularities $\sqrt\rho u,\sqrt\rho\theta\in C([0,{T_0}];L^{2})$
are replaced by
$$\rho u\in C([0,{T_0}];L^{2})\quad\mbox{ and }\quad \rho\theta\in C_w([0,{T_0}];L^{2}),$$
where $C_w$ represents the weak continuity.

(ii) Moreover, for any $t\in(0,{T_0})$, it holds that
\begin{eqnarray*}
  \|\sqrt{\rho}u\|^2(t) \leq\|\sqrt{\rho_{0}}u_{0}\|^2+Ct,\quad \norm{\sqrt{\rho}\theta}^{2}(t)\leq\norm{\sqrt{\rho_{0}}\theta_{0}}^{2}+C\sqrt t.
\end{eqnarray*}
\end{proposition}

\begin{proof}
 (i) \textbf{Step 1. Construction of the initial data.} Choose $\{u_{0n}\}_{n=1}^\infty\subseteq H_0^1\cap H^2$ such that $u_{0n}\rightarrow u_0$ in $H^1$ as $n\rightarrow\infty$. Set $\rho_{0n}=\rho_0+\frac{1}{n^2}$. Then, it is clear that
 \begin{equation}
 \label{APID1}
 \|\rho_{0n}\|_\infty+\|\nabla\rho_{0n}\|_q+\|\nabla u_{0n}\|^2\leq \|\rho_0\|_\infty+\|\nabla\rho_0\|_q+\|\nabla u_0\|^2+\frac12
 \end{equation}
 for large $n$. Put
  \begin{equation*}
 \overline{\theta}_{0n}(x)=
 \left\{
 \begin{aligned}
 &0,\qquad&&x\in\left\{x\in\Omega~\Big|~\rho_{0}(x)<\tfrac{1}{n}\right\},\\
 &\theta_{0},&&x\in\left\{x\in\Omega~\Big|~\rho_{0}(x)\geq\tfrac{1}{n}\right\},
 \end{aligned}
 \right.
 \end{equation*}
 and take $\theta_{0n}\geq0$ such that
 \begin{equation}
   \label{APID2}
   \|\theta_{0n}-\overline\theta_{0n}\|\leq\frac1n.
 \end{equation}
 Note that such $\theta_{0n}$ exists. For example, one can take $\theta_{0n}=j_{\varepsilon_n}*\widetilde\theta_{0n}$ for sufficiently small
 positive $\varepsilon_n$, where
 $j_\varepsilon$ is the standard mollifier and $\widetilde\theta_{0n}$ is the zero extension of $\overline\theta_{0n}$ on $\mathbb R^3$, that
 is, $\widetilde\theta_{0n}=\overline\theta_{0n}$ on $\Omega$ and $\widetilde\theta_{0n}=0$ on $\mathbb R^3\setminus\Omega$.

 We want to show
 \begin{eqnarray}
   \|\rho_{0n}\|_\infty+\|\nabla\rho_{0n}\|_q+\|\nabla u_{0n}\|^2+\|\sqrt{\rho_{0n}}\theta_{0n}\|^2\leq \Phi_0+1\label{ESTAPID}
 \end{eqnarray}
 for large $n$, and
 \begin{equation}
      \int_\Omega\rho_{0n}\theta_{0n}\chi dx \rightarrow\int_\Omega\rho_0\theta_0\chi dx\quad\mbox{as }n\rightarrow\infty,\quad\forall\chi \in L^2(\Omega). \label{COVAPID}
 \end{equation}

 The quantity $\|\sqrt{\rho_{0n}}\theta_{0n}\|$ is estimated as follows. By the elementary inequality
 $\sqrt{a+b}\leq\sqrt a+\sqrt b$ for $a,b\geq0$ and recalling the definition of $\overline\theta_{0n}$, it follows that
 \begin{eqnarray*}
 \sqrt{\rho_{0n}}\theta_{0n}&=&\sqrt{\rho_0+\frac1{n^2}}(\theta_{0n}-\overline\theta_{0n})+\sqrt{\rho_0+\frac1{n^2}}\overline\theta_{0n}\\
 &\leq&\left(\sqrt{\rho_0}+\frac1n\right)|\theta_0n-\overline\theta_{0n}|+\sqrt{\rho_0}\overline\theta_{0n}+\frac{\overline\theta_{0n}}{n}\\
 &\leq&\left(\sqrt{\rho_0}+\frac1n\right)|\theta_0n-\overline\theta_{0n}|+\sqrt{\rho_0} \theta_{0}+\frac{\overline\theta_{0n}}{n}
 \end{eqnarray*}
 and, thus, recalling (\ref{APID2}), one gets
 \begin{eqnarray*}
   \|\sqrt{\rho_{0n}}\theta_{0n}\|\leq \left(\|\rho_0\|_\infty^\frac12+\frac1{n}\right)\frac1n+\|\sqrt{\rho_0} \theta_{0}\|+
   \frac{\|\overline\theta_{0n}\|}{n}.
 \end{eqnarray*}
 With the aid of the above and noticing that
  \begin{eqnarray}
  \|\overline\theta_{0n}\|&=&\left(\int_{\Omega\cap\left\{x|\rho_0(x)\geq\frac1n\right\}} \theta_0^2 dx\right)^\frac12 \leq
    \sqrt n\left(\int_{\Omega\cap\{x|\rho_0\geq\frac1n\}}\rho_0\theta_0^2dx\right)^\frac12\nonumber\\
    &\leq&\sqrt n\left(\int_{\Omega}\rho_0\theta_0^2dx\right)^\frac12=\sqrt n \|\sqrt{\rho_0} \theta_{0}\|,
   \label{APID3}
 \end{eqnarray}
 one obtains
  \begin{eqnarray}
   \|\sqrt{\rho_{0n}}\theta_{0n}\|
   \leq\left(\|\rho_0\|_\infty^\frac12+\frac1{n}\right)\frac1n+\left(1+\frac{1}{\sqrt n}\right)\|\sqrt{\rho_0} \theta_{0}\|.
   \label{APID-ESTHETA}
 \end{eqnarray}
 This implies
 \begin{eqnarray}
   \|\sqrt{\rho_{0n}}\theta_{0n}\|^2\leq\|\sqrt{\rho_0}\theta_0\|^2+\frac12,\quad\mbox{for large }n. \label{APID4}
 \end{eqnarray}
 Combining (\ref{APID1}) with (\ref{APID4}) leads to (\ref{ESTAPID}).

Thanks to (\ref{APID2}) and (\ref{APID3}), it follows for any $\chi\in L^2(\Omega)$ that
\begin{eqnarray*}
  &&\left|\int_\Omega\rho_{0n}\theta_{0n}\chi dx-\int_\Omega\rho_0\theta_0\chi dx\right|\\
  &=&\left|\int_\Omega\left[\rho_{0n}(\theta_{0n}-\overline\theta_{0n})+(\rho_{0n}-\rho_0)\overline\theta_{0n}+\rho_0(\overline\theta_{0n}-\theta_0)
  \right]\chi dx\right|\\
  &\leq&\|\rho_{0n}\|_\infty\|\theta_{0n}-\overline\theta_{0n}\|\|\chi\|+\frac1{n^2}\|\overline\theta_{0n}\|\|\chi\|+\int_{\Omega\cap\{x|\rho_0
  (x)<\frac1n\}}\rho_0\theta_0|\chi|dx\\
  &\leq&\frac{\|\chi\|}{n}\left(\|\rho_0\|_\infty+\frac1{n^2}\right)+\frac{\|\chi\|}{n^{\frac32}}\|\sqrt{\rho_0}\theta_0\|+\frac{\|\chi\|}
  {\sqrt n}\|\sqrt{\rho_0}\theta_0\|,
\end{eqnarray*}
which implies (\ref{COVAPID}).

\textbf{Step 2. Approximate solutions and convergence.} Thanks to (\ref{ESTAPID}) and Proposition \ref{v}, there are two positive
constants ${T_0}$ and $C$ independent of $n$ such that system \eqref{a}-\eqref{c}, subject to \eqref{d}-\eqref{e}, admits a unique solution $(\rho_{n},u_{n},\theta_{n})$, in $\Omega\times(0,{T_0})$, and the following a priori estimates hold
\begin{equation}
  \label{aa}
  \left.
  \begin{aligned}
    &\int_{0}^{{T_0}} \norm{(\nabla\theta_n,\sqrt{\rho_n}\partial_tu_n,\nabla^{2}u_n)}^2dt\leq C,\\
    &\sup_{0\leq t\leq {T_0}}  \left(\|\rho_n\|_\infty+
    \norm{\rho_n}_{W^{1,q}}+\norm{(\partial_t\rho_n,\sqrt{\rho_n}\theta_n,\nabla u_n)}^2\right)\leq C,\\
    &\int_0^{{T_0}}\left(\norm{(\sqrt{t}\sqrt{\rho_n}\partial_t\theta_n,\sqrt{t}\nabla^{2}\theta_n,\sqrt{t}\nabla \partial_tu_n)}^{2}+\norm{\sqrt{t}\nabla^{2}u_n}_{q}^{2}\right)\mm\leq C,\\
    &\int_0^{{T_0}} t^2(\|\nabla\partial_t\theta_n\|^{2}+\|\nabla^2\theta_n\|_6^2)\mm\leq C,\\
    &\sup_{0\leq t\leq {T_0}}  \norm{(\sqrt{t}\nabla\theta_n,t\sqrt{\rho_n}\partial_t\theta_n,t\nabla^{2}\theta_n,\sqrt{t}\sqrt{\rho_n}\partial_tu_n,\sqrt{t}\nabla^{2}u_n)}^{2}  \leq C,
  \end{aligned}
  \right\}
\end{equation}
for large $n$. Then, by
the Banach-Alaoglu theorem and using the Cantor's diagonal arguments, there is a subsequence, still denoted by $(\rho_{n},u_{n},\theta_{n})$, and $(\rho,u,\theta)$ satisfying
\begin{align}
&\rho\in \label{ab}L^{\infty}(0,{T_0};W^{1,q}),\quad\rho_{t}\in L^{\infty}(0,{T_0};L^{2}),\\
&	\label{ac}\theta\in L^{2}(0,{T_0};H_{0}^{1}),\quad u\in L^{\infty}(0,{T_0};H_{0}^{1})\cap L^{2}((0,{T_0};H^{2}),\\
&\label{ad}\sqrt{t}\nabla\theta,\sqrt t\nabla^2u\in L^{\infty}(0,{T_0};L^{2}),\quad t\nabla^{2}\theta\in L^{\infty}(0,{T_0};L^{2})
\cap L^2(0,{T_0}; L^6),\\
&\label{ae}\sqrt{t}\nabla^{2}\theta\in L^{2}(0,{T_0};L^{2}),\quad t\nabla\theta_{t}\in L^{2}(0,{T_0};L^{2}),\\
&\label{af}\sqrt{t}\nabla u_{t}\in L^{2}(0,{T_0};L^{2}),\quad\sqrt{t}\nabla^{2}u\in L^{2}(0,{T_0};L^{q}),
\end{align}
such that
\begin{align}
\label{ag}&\rho_{n}\stackrel{*}\rightharpoonup\rho,\qquad\mathrm{in}\ L^{\infty}(0,{T_0};W^{1,q}),\\
\label{ah}&\partial_{t}\rho_{n}\stackrel{*}\rightharpoonup\rho_{t},\qquad\mathrm{in}\ L^{\infty}(0,{T_0};L^{2}),\\
\label{ai}&u_{n}\stackrel{*}\rightharpoonup u,\qquad\mathrm{in}\ L^{\infty}(0,{T_0};H_{0}^{1}),\\
\label{aj}&u_{n}\rightharpoonup u,\quad\mathrm{in}\ L^{2}(0,{T_0};H^{2}),\\
\label{ak}&\partial_{t}u_{n}\rightharpoonup u_{t},\qquad\mathrm{in}\ L^{2}(\delta,{T_0};H_{0}^{1}),\\
\label{al}&\theta_{n}\stackrel{*}\rightharpoonup\theta,\qquad \mathrm{in}\ L^{\infty}(\delta,{T_0};H_{0}^{1}),\\
\label{am}&\theta_{n}\rightharpoonup\theta,\qquad \mathrm{in}\ L^{2}(\delta,{T_0};W^{2,6}),\\
\label{an}&\partial_{t}\theta_{n}\rightharpoonup\theta_{t},\qquad \mathrm{in}\ L^{2}(\delta,{T_0};H^{1}_{0}),
\end{align}
for any $\delta\in(0,{T_0})$. Moreover, since $W^{1,q}\hookrightarrow\hookrightarrow C(\overline{\Omega})$ for $q\in(3,6)$,
and $H^2\hookrightarrow\hookrightarrow H^1\hookrightarrow\hookrightarrow
L^2$, it follows from the Aubin-Lions lemma and (\ref{ag})--(\ref{an}) that
\begin{align}
 \label{ao}&\rho_{n}\rightarrow\rho, &&\mathrm{in}\  C([0,{T_0}] ;C(\overline{\Omega})),\\
 \label{ap}&u_{n}\rightarrow u,&&\mathrm{in}\  C([\delta,{T_0}];L^{2}(\Omega))\cap  L^{2}(\delta,{T_0};H_{0}^{1}(\Omega)),\\
 \label{aq}&\theta_{n}\rightarrow\theta,&&\mathrm{in}\ C([\delta,{T_0}];L^{2}(\Omega))\cap  L^{2}(\delta,{T_0};H_{0}^{1}(\Omega)).
\end{align}
Due to the convergence (\ref{ak}), (\ref{an}), and \eqref{ao}--\eqref{aq}, one has
the following convergence of the nonlinear terms
\begin{align}
\label{3.20-1}&(\rho_{n}u_{n},\sqrt{\rho_n}u_n,\rho_n\theta_n,\sqrt{\rho_n}\theta_n)\rightarrow (\rho u,\sqrt\rho u,\rho\theta,\sqrt\rho\theta) \quad\mbox{ in } C([\delta,T_0]; L^{2}),\\
\label{3.20-2}&(\rho_{n}\partial_{t}u_{n},\sqrt{\rho_n}\partial_tu_n,\rho_n\partial_t\theta_n,\sqrt{\rho_n}\partial_t\theta_n)
\rightharpoonup(\rho u_{t},\sqrt\rho u_t,\rho\theta_t,\sqrt\rho\theta_t) \quad\mbox{ in } L^{2}(\Omega\times(\delta,{T_0})),\\
\label{3.20-3}&\rho_{n}(u_{n}\cdot\nabla)u_{n}\rightarrow\rho(u\cdot\nabla)u, \quad \rho_{n}(u_{n}\cdot\nabla)\theta_{n}\rightarrow\rho(u\cdot\nabla)\theta,\quad \mbox{ in } L^{1}(\Omega\times(\delta,{T_0})), \\
\label{3.20-4}&\rho_{n}\theta_n\text{div}u_n\rightarrow \rho\theta\text{div}u,\quad\mathcal Q(\nabla u_n)\rightarrow\mathcal Q(\nabla u),\quad \mbox{ in }L^{1}(\Omega\times(\delta,{T_0})),
\end{align}
for any $\delta\in(0,{T_0})$.
By the weakly lower semi-continuity of norms, it follows from \eqref{aa}, \eqref{3.20-1}, \eqref{3.20-2} that
\begin{eqnarray*}
  \int_\delta^{T_0}(\|\sqrt\rho u_t\|^2+\|\sqrt t\sqrt\rho\theta_t\|^2)dt&\leq&\varliminf_{n\rightarrow\infty}\int_\delta^{T_0}(\|\sqrt{\rho_n}
  \partial_tu_n\|^2+\|\sqrt t\sqrt{\rho_n}\partial_t\theta_n\|^2)dt\leq C, \\
  \|\sqrt\rho\theta\|(t)&=&\lim_{n\rightarrow\infty}\|\sqrt{\rho_n}\theta_n\|\leq C,
\end{eqnarray*}
for any $\delta, t\in(0,{T_0})$ and for a positive constant $C$ independent of $\delta$ and $t$. Therefore,
\begin{equation}
  \label{3.20-5}
  \sqrt\rho\theta\in L^\infty(0,{T_0}; L^2)\quad\mbox{and}\quad\sqrt\rho u_t,\sqrt t\sqrt\rho\theta_t
  \in L^2(0,{T_0}; L^2).
\end{equation}
The regularity $u\in L^1(0,{T_0}; W^{2,q})$ can be proved in the same as in Proposition \ref{g}.

\textbf{Step 3. The existence.} Thanks to the convergence \eqref{ag}--\eqref{3.20-4}, one can take the limit as $n\rightarrow\infty$ to the equations of $(\rho_n, \theta_n, u_0)$ to show that $(\rho, u, \theta)$ satisfies equations
(\ref{a})--(\ref{c}) in the sense of distribution. Due to the regularities (\ref{ab})--(\ref{af}), one can further
show that $(\rho, u, \theta)$ satisfies (\ref{a})--(\ref{c}), a.e.\,in $\Omega\times(0,{T_0})$. The initial
condition $\rho|_{t=0}=\rho_0$ is guaranteed by (\ref{ao}) by recalling that $\rho_n|_{t=0}=\rho_0+\frac1{n^2}.$

To complete the proof of (i), one still needs to show the regularities
$\rho u\in C([0,{T_0}]; L^2)$ and $\rho\theta\in C_w([0,{T_0}]; L^2)$, as well as
the initial condition $(\rho u, \rho\theta)|_{t=0}=(\rho_0 u_0, \rho_0\theta_0).$ To this end,
noticing that $\rho u, \rho\theta\in C((0,{T_0}];L^2)$ guaranteed by (\ref{3.20-1}), it suffices to
show
\begin{eqnarray}
  \rho u\rightarrow\rho_0u_0\quad\mbox{ in }L^2, \quad\mbox{ as }t\rightarrow0,\label{ICU}\\
  \rho\theta\rightharpoonup\rho_0\theta_0\quad\mbox{ in }L^2,\quad\mbox{ as }t\rightarrow0. \label{ICT}
\end{eqnarray}

We first verify (\ref{ICU}). By (\ref{aa}), it follows from the Gagliardo-Nirenberg and H\"older inequalities that
\begin{eqnarray}
\nonumber\int_{0}^{{T_0}}\norm{\partial_{t}(\rho_n u_n)}^{2}\mm&\leq& 2\int_0^{T_0}\left(\|\partial_t\rho_nu_n\|^2+\|\rho_n\partial_tu_n\|^2\right)dt\\
&\leq& C\int_0^{T_0}\left(\|\partial_t\rho_n\|^2\|u_n\|_\infty^2+\|\rho_n\|_\infty\|\sqrt{\rho_n}\partial_tu_n\|^2\right) dt\nonumber\\
\label{ar}&\leq& C\int_{0}^{{T_0}} \|\nabla u_n\|\norm{\nabla^{2}u_n}\mm+C \leq C
\end{eqnarray}
for large $n$.
Thanks to this, it follows from the Newton-Leibnitz formula, the Minkowski and H\"older inequalities that
\begin{eqnarray}
	\nonumber&&\norm{\rho u(\cdot,t)-\rho_{0} u_{0}}\\
	\nonumber&\leq& \norm{\rho u-\rho_{n}u_{n}}(t)+\norm{\rho_{n}u_{n}-\rho_{0n}u_{0n}}(t)+\norm{\rho_{0n}u_{0n}-\rho_{0n}u_{0}}+\norm{\rho_{0n}u_{0}-\rho_{0}u_{0}}\\
	\nonumber&\leq&\norm{\rho u-\rho_{n}u_{n}}(t)+\int_{0}^{t}\norm{\partial_{t}(\rho_{n}u_{n})}\mathrm{d}\tau+\norm{\rho_{0n}(u_{0n}-u_{0})}
+\frac{C}{n^2}\norm{u_{0}}\\
	\label{dm}&\leq&\norm{\rho u-\rho_{n}u_{n}}(t)+C\sqrt{t}+\norm{\rho_{0n}}_{\infty}\norm{u_{0n}-u_{0}}+\frac{C}{n^2}\norm{u_{0}}
\end{eqnarray}
for large $n$, from which, recalling (\ref{3.20-1}) and $u_{0n}\rightarrow u_0$ in $H^1$ as $n\rightarrow\infty$, one gets
by taking $n\rightarrow\infty$ that
$\|\rho u-\rho_0u_0\|(t)\leq C\sqrt t$, proving (\ref{ICU}).

Then, we verify (\ref{ICT}). Since $\rho\theta\in L^\infty(0,{T_0}; L^2)$ and $C_c^\infty(\Omega)$ is dense
in $L^2$, it suffices to verify
\begin{equation}
  \label{ICT'}
  \left(\int_\Omega\rho\theta\phi dx\right)(t)\rightarrow\int_\Omega\rho_0\theta_0\phi dx\quad\mbox{as }t\rightarrow0,\quad\forall \phi\in C_c^\infty(\Omega).
\end{equation}
Rewrite the equation for $\theta_n$ as
$$
c_v\left[\partial_t(\rho_n\theta_n)+\text{div}(\rho_n\theta_nu_n)\right]+R\rho_n\theta_n\text{div}u_n-\kappa\Delta\theta_n=
\mathcal Q(\nabla u_n).
$$
Multiplying the above equation with $\phi\in C_c^\infty(\Omega)$ and integrating over $\Omega\times(0,t)$ yield
\begin{eqnarray*}
&&c_v\left[\left(\int_\Omega\rho_n\theta_n\phi dx\right)(t)-\int_{\Omega}\rho_{0n}\theta_{0n}\phi\mathrm{d}x \right]\\
&=&c_{v}\int_{0}^{t}\int_{\Omega}\rho_{n}\theta_{n}u_{n}\cdot\nabla\phi\mathrm{d}xd\tau+\kappa\int_{0}^t
\int_{\Omega}\Delta\theta_{n}\phi\mathrm{d}xd\tau\\
&&-R\int_{0}^t\int_{\Omega}\rho_{n}\theta_{n}\dive u_{n}\phi\mathrm{d}xd\tau+\int_{0}^t\int_{\Omega}\mathcal{Q}(\nabla u)\phi\mathrm{d}xd\tau=:\sum\limits_{i=1}^{4}M_{i}.
\end{eqnarray*}
Terms on the right-hand side are estimated by integration by parts, the H\"older inequality, and (\ref{aa}) as follows:
\begin{eqnarray}
\nonumber |M_1|&\leq& c_v
\int_{0}^t\|\rho_n\|^\frac12_\infty\|\sqrt{\rho_n}\theta_n\|\|u_{n}\|_6\norm{\nabla\phi}_3 d\tau\leq Ct,\\
	\nonumber |M_{2}|&\leq&\kappa\int_{0}^t\int_{\Omega}|\nabla\theta_{n}||\nabla\phi|\mathrm{d}xd\tau\leq C\left(\int_{0}^t\norm{\nabla\theta_{n}}^2d\tau\right)^\frac12\sqrt t\leq C\sqrt t,\\
	\nonumber |M_{3}|&\leq& C\int_{0}^t\norm{\rho_{n}}_{\infty}^{\frac{1}{2}}\norm{\sqrt{\rho}\theta_{n}}\norm{\nabla u_{n}} \norm{\phi}_\infty d\tau \leq C t,\\
	\nonumber |M_4|&\leq& C\int_{0}^t\int_{\Omega}|\nabla u_{n}|^{2}|\phi|\mathrm{d}xd\tau\leq C\int_{0}^t\norm{\nabla u_{n}}^2\norm{\phi}_\infty d\tau\leq Ct,
\end{eqnarray}
for large $n$.
Therefore, for large $n$, it follows
$$
\left|\left(\int_\Omega\rho_n\theta_n\phi dx\right)(t)-\int_{\Omega}\rho_{0n}\theta_{0n}\phi\mathrm{d}x\right|\leq C\sqrt t,\quad\forall t\in(0,{T_0}),
$$
for any $\phi\in C_c^\infty(\Omega)$ and for a positive constant $C$ independent of $n$. Thanks to this and recalling
(\ref{COVAPID}) and (\ref{3.20-1}), one gets by taking $n\rightarrow\infty$ that
$$
\left|\left(\int_\Omega\rho \theta \phi dx\right)(t)-\int_{\Omega}\rho_{0}\theta_{0}\phi\mathrm{d}x\right|\leq C\sqrt t,\quad\forall t\in(0,{T_0}),\quad\forall \phi\in C_c^\infty(\Omega),
$$
verifying (\ref{ICT}).

(ii) Multiplying equation \eqref{b} for $(\rho_n, u_n, \theta_n)$ with $u_{n}$ and integrating over $\Omega$, one gets by integration by
parts that
$$
 \fracc\int_{\Omega}\frac{\rho_{n}}{2}|u_{n}|^{2}\dx+\mu\inn|\nabla u_{n}|^{2}\dx+(\mu+\lambda)\inn|\dive u_{n}|^{2}\dx-\int_{\Omega} P_{n}\text{div} u_{n}\dx=0.
$$
Integrating the above with respect to $t$, by the H\"older inequality, and using (\ref{aa}), one deduces for large $n$ that
\begin{eqnarray*}
  \|\sqrt{\rho_n}u_n\|^2(t)&\leq&\|\sqrt{\rho_{0n}}u_{0n}\|^2+2\int_0^t\int_\Omega P_n\text{div}u_n dxd\tau\\
  &\leq&\|\sqrt{\rho_{0n}}u_{0n}\|^2+2R\int_0^t\|\rho_n\|_\infty^\frac12\|\sqrt{\rho_n}\theta_n\|\|\nabla u_n\|d\tau\\
  &\leq&\|\sqrt{\rho_{0n}}u_{0n}\|^2+Ct,\qquad\forall t\in(0,{T_0}),
\end{eqnarray*}
for a positive constant $C$ independent of $n$. Thanks to the above, recalling (\ref{3.20-1}) and noticing that $\sqrt{\rho_{0n}}u_{0n}\rightarrow\sqrt{\rho_0}u_0$ in $L^2$ as $n\rightarrow\infty$, one gets by taking $n\rightarrow\infty$
that
\begin{equation}
  \|\sqrt{\rho}u\|^2(t) \leq\|\sqrt{\rho_{0}}u_{0}\|^2+Ct,\qquad\forall t\in(0,{T_0}).\label{3.35}
\end{equation}
Multiplying equation \eqref{c} for $(\rho_n, u_n, \theta_n)$ with $\theta_{n}$ and integrating over $\Omega$, one gets by integration by parts,
the Sobolev embedding inequality, and (\ref{aa}) that
\begin{eqnarray*}
	\frac{c_{v}}{2}\fracc\norm{\sqrt{\rho_{n}}\theta_{n}}^{2} +\kappa\norm{\nabla\theta_{n}}^{2}
	&=&-\inn\dive u_{n}P_{n}\theta_{n}\dx+\inn\mathcal{Q}(\nabla u_{n})\theta_{n}\dx\\
	&\leq& C\left(\norm{\rho_{n}}_\infty^{\fr}\norm{\sqrt{\rho_{n}}\theta_{n}} \norm{\nabla u_{n}}_{3}+\norm{\nabla u_{n}}\norm{\nabla u_{n}}_{3} \right)\norm{\theta_{n}}_{6} \\
&\leq&\frac{\kappa}{2}\norm{\nabla\theta_{n}}^{2}+ C \norm{\nabla^{2}u_{n}},
\end{eqnarray*}
for large $n$, from which, integrating with respect to $t$, using \eqref{aa} again, and by the H\"older inequality, one obtains
\begin{eqnarray*}
\norm{\sqrt{\rho_{n}}\theta_{n}}^{2}(t)&\leq&\norm{\sqrt{\rho_{0n}}\theta_{0n}}^{2}+C\int_{0}^{t} \norm{\nabla^{2}u_{n}} \mathrm{d}\tau\\
	&\leq&\norm{\sqrt{\rho_{0n}}\theta_{0n}}^{2}+C\sqrt t,\qquad\forall t\in(0,{T_0}).
\end{eqnarray*}
Thanks to this and recalling (\ref{APID-ESTHETA}) and (\ref{3.20-1}), one can take $n\rightarrow\infty$ to get
\begin{align*}
	\norm{\sqrt{\rho}\theta}^{2}(t)\leq\norm{\sqrt{\rho_{0}}\theta_{0}}^{2}+C\sqrt t,\qquad\forall t\in(0,{T_0}).
\end{align*}
Combining this with (\ref{3.35}), the conclusion follows.
\end{proof}

 \section{Proof of Theorem \ref{u}}

This section is devoted to the proof of Theorem \ref{u}. As already explained in the introduction that the existence of strong solutions, which
enjoy all the regularities stated in Definition \ref{def1} except that $\sqrt\rho u, \sqrt\rho\theta\in C([0,T_0]; L^2)$,
is proved directly in the Euler coordinates, but the regularities $\sqrt\rho u, \sqrt\rho\theta\in C([0,T_0]; L^2)$ and the uniqueness
are proved in the Lagrangian coordinates first and later transformed back to the Euler coordinates.

 \subsection{Lagrangian coordinates and some properties}
 Given a velocity field $u\in L^1(0,{T_0}; C^1(\overline\Omega))$ satisfying $u|_{\partial\Omega}=0$ and let
$x=\varphi(y,t)$
 be the corresponding coordinates transform, governed by the velocity field $u$, between the Euler coordinates $(x,t)$ and the Lagrangian
 coordinates $(y,t)$, that is,
 \begin{equation}\label{LG1}
 \left\{
 \begin{aligned}
 &\partial_{t}\varphi(y,t)=u(\varphi(y,t),t),\qquad\forall\ t\in[0,{T_0}],\\
 &\varphi(y,0)=y.
 \end{aligned}
 \right.
 \end{equation}
By the classic theory for ODEs, $\varphi$ is well-defined and $\varphi: \Omega\times[0,{T_0}]\rightarrow\Omega$. Moreover,
by the unique solvability of ODEs, for each $t\in[0,{T_0}]$, $\varphi(\cdot,t): \Omega\rightarrow\Omega$
is bijective. Denote by
$y=\psi(x,t)$ the inverse mapping of $x=\varphi(y,t)$ with respect to $y$, which satisfies
 \begin{equation}\label{LG1-1}
 \left\{
 \begin{aligned}
 &\partial_t\psi(x,t)+(u(x,t)\cdot\nabla)\psi(x,t)=0,\\
 &\psi(x,t)|_{t=0}=x.
 \end{aligned}
 \right.
 \end{equation}
Set
\begin{eqnarray}
  &&A(y,t)  =(a_{ij}(y,t))_{3\times3},\quad  a_{ij}(y,t)=\partial_i\psi_j(x,t)|_{x=\varphi(y,t)},\label{EQA} \\
  &&J(y,t)=\text{det } A(y,t)=\text{det }\nabla\psi(x,t)|_{x=\varphi(y,t)},\label{EQJ}\\
  &&B(y,t)=(b_{ij}(y,t))_{3\times3},\quad b_{ij}(y,t)=\partial_i\varphi_j(y,t).\label{EQB}
\end{eqnarray}
Then, one can check
\begin{equation}
\label{LG5}
  \left\{
  \begin{array}{l}
  \partial_tA(y,t)=-\nabla u(\varphi(y,t),t)A(y,t),\\
  A(y,t)|_{t=0}=I
  \end{array}
  \right.
\end{equation}
and
\begin{equation}
\label{LG4}
  \left\{
  \begin{array}{l}
  \partial_tB(y,t)=B(y,t)\nabla u(\varphi(y,t),t),\\
  B(y,t)|_{t=0}=I,
  \end{array}
  \right.
\end{equation}
here we set $\nabla u=(\partial_i u_j)_{3\times3}$.
Recalling the definition of $J$, one derives from (\ref{LG5}) that
\begin{equation}\label{LG6}
  \left\{
  \begin{array}{l}
  \partial_t J(y,t)=-\text{div} u(\varphi(y,t),t)J(y,t),\\
  J(y,t)|_{t=0}=1,
  \end{array}
  \right.
\end{equation}

Some properties of the mapping $\varphi$ are stated in the next two propositions whose proofs are postponed
in the Appendix.

\begin{proposition}
  \label{cs}
Given $u\in L^\infty(0,{T_0}; H_0^1)\cap L^1(0,{T_0}; W^{2,q})$, with $q\in(3,6)$, and let $\varphi$,
$\psi$, $A$, $B$, and $J$ be defined as before.

Then, $J>0$ on $\Omega\times(0,{T_0})$ and the following hold:
\begin{align}
&\sup_{0\leq t\leq{T_0}}\left(\left\|\left(\frac1J,J,A,B\right)\right\|_\infty+
 \|(J_t, A_t, B_t)\|+\|(\nabla J, \nabla A, \nabla B)\|_q\right)\leq C, \label{CS2}\\
&\|\nabla[g(\varphi(\cdot,t))]\|_{W^{1,\alpha}} \simeq \|\nabla g\|_{W^{1,\alpha}}
  \simeq\|\nabla[g(\psi(\cdot,t))]\|_{W^{1,\alpha}},\quad\forall\alpha\in[1,q],\label{CS3}\\
&\|\nabla[g(\varphi(\cdot,t))]\|_\alpha \simeq \|\nabla g\|_\alpha\simeq\|\nabla[g(\psi(\cdot,t))]\|_\alpha,\quad\forall\alpha\in[1,\infty],
\label{CS4}\\
&\|g(\varphi(\cdot ,t))\|_\alpha \simeq \|g\|_\alpha\simeq\|g(\psi(\cdot,t))\|_\alpha,\quad\forall\alpha\in[1,\infty], \label{CS5}
\end{align}
for any function $g$ such that all the relevant quantities are finite, here we denote $\mathscr Q_1\simeq \mathscr Q_2$ means $\frac{\mathscr Q_1}{\overline C}\leq\mathscr Q_2\leq\overline C\mathscr Q_1$ for a positive
constant $\overline C$ depending only on $\Omega, \alpha, q, T_0,$ and $\|u\|_{L^\infty(0,T_0; H^1(\Omega))\cap L^1(0,T_0; L^1(0,T_0; W^{2,q}(\Omega))}$.
\end{proposition}

\begin{proposition}
  \label{CONTINUITY}
Under the assumptions as in Proposition \ref{cs}, the following hold:

(i) $h(\varphi(\cdot,t),t)\in C([0,{T_0}]; L^2)$ if $h\in C([0,{T_0}]; L^2)$;

(ii) $h(\varphi(\cdot,t),t)\in C_w([0,{T_0}]; L^2)$ if $h\in C_w([0,{T_0}]; L^2)$.
\end{proposition}

\subsection{Regularities and reduced system in the Lagrangian coordinates}
Given initial data $(\rho_0, u_0, \theta_0)$ satisfying the assumptions in Theorem \ref{u}. Let $(\rho, u, \theta)$
be the solution established in Proposition \ref{existence1} and $\varphi$ the corresponding mapping defined by (\ref{LG1}).
Set
\begin{equation}\label{RVT}
\left\{
\begin{aligned}
  \varrho(y,t)=\rho(\varphi(y,t),t),\\
  v(y,t)=u(\varphi(y,t),t),\\
  \vartheta(y,t)=\theta(\varphi(y,t),t).
\end{aligned}
\right.
\end{equation}

As direct corollaries of Proposition \ref{cs} and Proposition \ref{CONTINUITY} and recalling the regularities of $(\rho, u, \theta)$ in Proposition \ref{existence1}, one has:
\begin{equation}
  \label{REGRVT1}
  \left\{
\begin{array}{l}
  \varrho\in L^\infty(0,{T_0}; W^{1,q}),\\
  v\in L^\infty(0,{T_0}; H^1)\cap L^2(0,{T_0}; H^2)\cap L^1(0,{T_0}; W^{2,q}),\\
  \sqrt tv\in L^\infty(0,{T_0}; H^2)\cap L^2(0,{T_0}; W^{2,q}), \\
  \sqrt\varrho\vartheta\in L^\infty(0,{T_0}; L^2),\quad \vartheta\in L^2(0,{T_0}; H^1), \\
  \sqrt t\vartheta\in L^\infty(0,{T_0}; H^1)\cap L^2(0,{T_0}; H^2),\\
  \varrho v\in C([0,{T_0}]; L^2),\quad\varrho\vartheta\in C_w([0,{T_0}]; L^2).
\end{array}
\right.
\end{equation}

Direct calculations show
\begin{eqnarray*}
  &&\partial_t\varrho(y,t)=(\partial_t\rho+u\cdot\nabla\rho)(\varphi(y,t),t),\label{PTR}\\
  &&\partial_tv(y,t) = [\partial_tu+(u\cdot\nabla)u](\varphi(y,t),t),\label{PTV}\\
  &&\partial_t\vartheta(y,t)=(\partial_t\theta+u\cdot\nabla\theta)(\varphi(y,t),t)\label{PTT}.
\end{eqnarray*}
By Proposition \ref{cs} and recalling the regularities of $(\rho, u, \theta)$ stated in Proposition \ref{existence1},
one deduces by the H\"older and Gagliardo-Nirenberg inequalities that
\begin{eqnarray*}
  \sup_{0\leq t\leq T}\|\varrho_t\|&=&\sup_{0\leq t\leq{T_0}}\|(\rho_t+u\cdot\nabla\rho)(\varphi(y,t),t)\|\\
  &\leq&C\sup_{0\leq t\leq{T_0}}\| \rho_t+u\cdot\nabla\rho \| \leq C\sup_{0\leq t\leq{T_0}}(\| \rho_t\|+\|u\|_6\|\nabla\rho \|_3)\\
  &\leq& C\sup_{0\leq t\leq{T_0}}(\| \rho_t\|+\|\nabla u\|_2\|\nabla\rho \|_q)\leq C,\\
  \int_0^{T_0}\|\sqrt\varrho v_t\|^2dt&=&\int_0^{T_0}\|[\sqrt\rho(\partial_tu+(u\cdot\nabla)u)](\varphi(y,t),t)\|^2dt
  \nonumber\\
  &\leq&C\int_0^{T_0}\|\sqrt\rho(\partial_tu+(u\cdot\nabla)u)\|^2dt\nonumber\\
  &\leq& C\int_0^{T_0}\left(\|\sqrt\rho u_t\|^2+\|u\|_6^2\|\nabla u\|_3^2\right) dt\nonumber\\
  &\leq&C\int_0^{T_0}\left(\|\sqrt \rho u_t\|^2+\|\nabla u\|^3\|\nabla^2u\|\right)dt\leq C
\end{eqnarray*}
and
\begin{eqnarray*}
  \int_0^{T_0} t\|\nabla v_t\|^2dt&=&\int_0^{T_0}t\|\nabla [(\partial_tu+(u\cdot\nabla)u)(\varphi(y,t),t)]\|^2dt
  \nonumber\\
  &\leq&C\int_0^{T_0}t\|\nabla(\partial_tu+(u\cdot\nabla)u)\|^2dt \nonumber\\
  &\leq&C\int_0^{T_0} t\left(\|\nabla u_t\|^2+\|u\|_\infty^2\|\nabla^2u\|^2+\|\nabla u\|_4^4\right)dt \nonumber\\
  &\leq&C\int_0^{T_0} t\left(\|\nabla u_t\|^2+\|\nabla u\|_2\|\nabla^2u\|^3\right)dt\leq C.
\end{eqnarray*}
Similarly,
\begin{eqnarray*}
  \int_0^{T_0} t\|\sqrt\varrho\vartheta_t\|^2dt&=&\int_0^{T_0} t\|(\sqrt\rho\theta_t+\sqrt\rho u\cdot\nabla\theta)(\varphi(y,t), t)\|^2dt\nonumber \\
  &\leq&C\int_0^{T_0} t\|\sqrt\rho\theta_t+\sqrt\rho u\cdot\nabla\theta\|^2dt \nonumber\\
  &\leq& C\int_0^{T_0} t(\|\sqrt\rho\theta_t\|^2+\|\rho\|_\infty\|u\|_6^2\|\nabla\theta\|_3^2)dt\nonumber\\
  &\leq&C\int_0^{T_0} t(\|\sqrt\rho\theta_t\|^2+\|\nabla u\|^2\|\nabla\theta\|\|\nabla^2\theta\|)dt\leq C
\end{eqnarray*}
and
\begin{eqnarray*}
  \int_0^{T_0}t^2\|\nabla\vartheta_t\|^2dt&=&\int_0^{T_0} t^2\|\nabla[(\theta_t+u\cdot\nabla\theta)(\varphi(y,t),t)]\|^2dt\\
  &\leq&C\int_0^{T_0}t^2\|\nabla (\theta_t+u\cdot\nabla\theta) \|^2dt\\
  &\leq&C\int_0^{T_0}t^2\int_\Omega(|\nabla\theta_t|^2+|u|^2|\nabla^2\theta|^2+|\nabla u|^2|\nabla\theta|^2) dxdt \\
  &\leq&C\int_0^{T_0}t^2(\|\nabla\theta_t\|^2+\|u\|_\infty^2\|\nabla^2\theta\|^2+\|\nabla u\|_6^2\|\nabla\theta\|\|\nabla\theta\|_6)dt\\
  &\leq&C\int_0^{T_0}t^2(\|\nabla\theta_t\|^2+\|\nabla u\|\|\nabla^2u\|\|\nabla^2\theta\|^2)dt\\
  &&+C\int_0^{T_0} t^2\|\nabla^2 u\|^2\|\nabla\theta\|\|\nabla^2\theta\| dt\leq C.
\end{eqnarray*}
Therefore,
\begin{equation}
\label{REGRVT3}
\left.
\begin{aligned}
  \varrho_t\in L^\infty(0,{T_0}; L^2),\quad
  \sqrt\varrho v_t\in L^2(0,{T_0}; L^2), \quad
  \sqrt tv_t\in L^2(0,{T_0}; H^1),\\
 \sqrt t\sqrt\varrho\vartheta_t\in L^2(0,{T_0}; L^2),\quad
  t\vartheta_t\in L^2(0,{T_0}; H^1).
  \end{aligned}
  \right\}
\end{equation}

So, we have the following proposition:

\begin{proposition}
  \label{PROPREGRVT}
Given initial data $(\rho_0, u_0, \theta_0)$ satisfying the assumptions in Theorem \ref{u}. Let $(\rho, u, \theta)$
be the solution established in Proposition \ref{existence1} and $\varphi$ the corresponding mapping defined by
(\ref{LG1}). Then, $(\varrho, v, \vartheta)$ defined by (\ref{RVT}) satisfies (\ref{REGRVT1}) and (\ref{REGRVT3}).
\end{proposition}

Let $A$ be defined as before in the previous subsection and denote
 $$
\nabla_{A} f:= A \nabla f,\quad \dive_{A}v:=A:(\nabla v)^{T},\quad \nabla v=(\partial_iv_j)_{3\times3}.
 $$
Then, by direct computations, one can derive from (\ref{a})--(\ref{c}) and (\ref{LG5})--(\ref{LG6}) that
\begin{align}
 \label{bd'}\varrho_{t}+\dive_{A}v\varrho=0,\\
 \label{bg}J\varrho_{0} v_{t}-\mu\dive_{A}(\nabla_{A}v)-(\mu+\lambda)\nabla_{A}(\dive_{A}v)+R\nabla_{A}(J\varrho_{0}\vartheta)=0,\\
 \label{bh}c_{v}J\varrho_{0}\vartheta_{t}+RJ\varrho_{0}\vartheta\dive_{A}v-\kappa\dive_{A}(\nabla_{A}\vartheta)=\frac{\mu}{2}\left|\nabla_{A}v+(\nabla_{A}v)^{T} \right|^{2}+\lambda(\dive_{A}v)^{2},\\
 \label{bi} A_{t}+\nabla_{A} v A=0,\\
 \label{bj}J_{t}+\dive_{A}vJ=0.
\end{align}
System (\ref{bd'})--(\ref{bj}) are satisfied a.e.\,in $\Omega\times(0,{T_0})$.
Here in (\ref{bg}) and (\ref{bh}) we have used the fact that
\begin{equation}
\frac{\varrho}{J}=\frac{\rho_{0}}{J_{0}}=\rho_{0}\label{RHO-RHO0}
\end{equation}
to replace $\varrho$ with $J\rho_0$, as $\partial_t(\frac{\varrho}{J})=0$ guaranteed by (\ref{bd'}) and (\ref{bj}).
The component form of (\ref{bg}) reads as
\begin{align*}
J\rho_0\partial_{t}v_{i}-\mu a_{kl}\partial_{l}(a_{km}\partial_{m}v_{i})-(\mu+\lambda)a_{il}\partial_{l}(a_{km}\partial_{m}v_{k})
+Ra_{il}\partial_{l}(\varrho\vartheta)=0.
\end{align*}

The initial-boundary conditions read as
\begin{eqnarray}
 \label{bk}&(\varrho v,\varrho\vartheta, A, J)|_{t=0}=(\rho_0u_0, \rho_0\theta_0, I, 1),\\
 \label{bl}&v|_{\partial\Omega}=0,\quad\vartheta|_{\partial\Omega}=0.
\end{eqnarray}
Since $\varrho v\in C([0,{T_0}]; L^2)$ and $\varrho\vartheta\in C_w([0,{T_0}]; L^2)$,
guaranteed by (\ref{REGRVT1}), and $A, J\in C([0,{T_0}]; L^2)$, guaranteed by Proposition \ref{cs},
the initial condition (\ref{bk}) is well-defined.

Finally, we state and prove the continuities of $\sqrt{\rho_0}v$ and $\sqrt{\rho_0}\vartheta$.

\begin{proposition}\label{cn}
Under the same assumptions in Proposition \ref{PROPREGRVT}, it holds that
$$
\sqrt{\rho_0}v,\sqrt{\rho_0}\vartheta\in C([0,{T_0}];L^{2}),
$$
and
$$
(\sqrt{\rho_0}v,\sqrt{\rho_0}\vartheta)\rightarrow (\sqrt{\rho_0}u_0, \sqrt{\rho_0}\theta_0),
\quad \mbox{ in }L^2,\quad\mbox{as }t\rightarrow0.
$$
\end{proposition}

\begin{proof}
We only give the proof for $\sqrt{\rho_0}\vartheta$ while that for $\sqrt{\rho_0}v$ can be
done similarly.
Due to (\ref{REGRVT1}) and (\ref{REGRVT3}), one has
$\sqrt{\rho_0}\vartheta\in C((0,{T_0}];L^{2}).$
It remains to show
\begin{equation}\label{COV-SQRT-2}
\sqrt{\rho_0}\vartheta\rightarrow\sqrt{\rho_0}\theta_0,
\quad \mbox{ in }L^2,\quad\mbox{as }t\rightarrow0.
\end{equation}
Noticing that $\inf_{(y,t)\in\Omega\times[0,{T_0}]}J(y,t)>0$, $J\in L^\infty(0,{T_0}; W^{1,q})$, and
$J_t\in L^\infty(0,{T_0}; L^2)$, guaranteed by Proposition \ref{cs}, one can verify easily that $\frac 1J\in C([0,{T_0}]; C(\overline\Omega))$. Thanks to this and recalling (\ref{RHO-RHO0}), it follows from
(\ref{REGRVT1}) and (\ref{bk}) that
\begin{equation}
\label{VF4}
\rho_0\vartheta=\frac1J\varrho\vartheta\rightharpoonup\rho_0\theta_0 \quad\mbox{ in }L^2(\Omega)\quad \mbox{ as }t\rightarrow0.
\end{equation}

In order to show $\sqrt{\rho_0}\vartheta\rightarrow\sqrt{\rho_0}\theta_0$ in $L^2$ as $t\rightarrow0$, one needs to verify
\begin{eqnarray}
  &\label{VF1}\sqrt{\rho_0}\vartheta\rightharpoonup\sqrt{\rho_0}\theta_0\quad \mbox{ in }L^2\quad\mbox{ as }
  t\rightarrow0,\\
  &\label{VF2}\displaystyle\varlimsup_{t\rightarrow0}\|\sqrt{\rho_0}\vartheta\|^2\leq\|\sqrt{\rho_0}\theta_0\|^2.
\end{eqnarray}
To verify (\ref{VF1}), since $\sqrt{\rho_0}\vartheta=0$ on $\Omega_0\times(0,{T_0})$ and $\sqrt{\rho_0}
\theta_0=0$ on $\Omega_0$, it suffices to show
$$
\sqrt{\rho_0}\vartheta\rightharpoonup\sqrt{\rho_0}\theta_0\quad\mbox{in }L^2(\Omega_+)\quad\mbox{as }t\rightarrow0,
$$
where $\Omega_+=\{y\in\Omega|\rho_0(y)>0\}$ and $\Omega_0=\{y\in\Omega|\rho_0(y)=0\}.$
Recalling that $\sqrt{\rho_0}\vartheta\in L^\infty(0,{T_0}; L^2)$ and since $C_c^\infty(\Omega_+)$ is dense in $L^2(\Omega_+)$, one only needs to check
\begin{equation}
  \int_{\Omega_+}\sqrt{\rho_0}\vartheta\chi dy\rightarrow\int_{\Omega_+}\sqrt{\rho_0}\theta_0\chi dy \quad\mbox{ as }t
  \rightarrow0,\quad\forall\chi\in C_c^\infty(\Omega). \label{VF3}
\end{equation}
Take arbitrary $\chi\in C_c^\infty(\Omega_+)$ and denote $S:=\{y\in\Omega_+|\chi(y)\not=0\}$. Then $S\subset\overline S\subset\Omega_+$. Since $\rho_0\in C(\overline\Omega)$ and $\rho_0>0$ on $\Omega_+$, it follows
that $\min_{y\in\overline S}\rho_0(y)>0$ and, thus, $\frac{\chi}{\sqrt{\rho_0}}\in L^2(S)$. So, it follows from (\ref{VF4}) that
\begin{equation*}
  \int_{\Omega_+}\sqrt{\rho_0}\vartheta\chi dy=\int_S\rho_0\vartheta \frac{\chi}{\sqrt{\rho_0}}dy
  \rightarrow\int_S\rho_0\theta_0 \frac{\chi}{\sqrt{\rho_0}}dy=
  \int_{\Omega_+}\sqrt{\rho_0}\theta_0\chi dy \quad\mbox{ as }t
  \rightarrow0.
\end{equation*}
Therefore, (\ref{VF3}) and further (\ref{VF1}) hold.

We now verify (\ref{VF2}). First, noticing that
$$
\sqrt{\rho_0}\vartheta=\frac{1}{\sqrt J}\sqrt\varrho\vartheta,\quad
\frac 1J\in C([0,{T_0}];C(\overline\Omega)),\quad J|_{t=0}=1,
$$
one has $\varlimsup_{t\rightarrow0}\|\sqrt{\rho_0}\vartheta\|^2=\varlimsup_{t\rightarrow0}\|\sqrt{\varrho}\vartheta\|^2$. Therefore,
it suffices to show
\begin{equation}\label{VF5}
\varlimsup_{t\rightarrow0}\|\sqrt{\varrho}\vartheta\|^2\leq\|\sqrt{\rho_0}\theta_0\|^2.
\end{equation}
Since $\text{det}\nabla\psi(x,t)=J(\psi(x,t),t)>0$, one deduces by
direct calculations that
\begin{eqnarray*}
  \|\sqrt\varrho\vartheta\|^2(t) &=&\int_\Omega\rho(\varphi(y,t),t)\theta^2(\varphi(y,t),t)dy \\
  &=&\int_\Omega\rho(x,t)\theta^2(x,t)|\text{det}\nabla\psi(x,t)|dx=\int_\Omega\rho(x,t)\theta^2(x,t)J(\psi(x,t),t)dx\\
  &=&\int_\Omega\rho(x,t)\theta^2(x,t)dt+\int_\Omega\rho(x,t)\theta^2(x,t)[J(\psi(x,t),t)-1]dt\\
  &\leq&\|\sqrt\rho\theta\|^2(t)+\|J-1\|_\infty(t)\|\sqrt\rho\theta\|^2(t).
\end{eqnarray*}
Thanks to the above and since
$$
\sqrt\rho\theta\in L^\infty(0,{T_0}; L^2),\quad
J\in C([0,{T_0}]; C(\overline\Omega)),\quad J|_{t=0}=1,
$$
it follows from Proposition \ref{existence1} that
\begin{eqnarray*}
  \varlimsup_{t\rightarrow0}\|\sqrt\varrho\vartheta\|^2(t)&\leq&\varlimsup_{t\rightarrow0}
  \|\sqrt\rho\theta\|^2(t)+C\varlimsup_{t\rightarrow}\|J-1\|_\infty(t)\\
  &=&\varlimsup_{t\rightarrow0}
  \|\sqrt\rho\theta\|^2(t) \leq \|\sqrt{\rho_0}\theta_0\|^2.
\end{eqnarray*}
This verifies (\ref{VF5}) and further (\ref{VF2}).
 \end{proof}

\subsection{Proof of Theorem \ref{u}}
We are now ready to give the proof of Theorem \ref{u}.
\begin{proof}[\textbf{Proof of Theorem \ref{u}}]
{(i) \textbf{Existence.}} By virtue of Proposition \ref{existence1}, it remains to show that
$\sqrt\rho u,\sqrt\rho\theta\in C([0,{T_0}];L^2)$. Let $(\varrho, v, \vartheta)$ be given by (\ref{RVT}).
Then, it follows from Proposition \ref{cn} that $\sqrt{\rho_0}v, \sqrt{\rho_0}\vartheta\in C([0,{T_0}];
L^2)$, from which, recalling that $\varrho=J\rho_0$ and noticing that $J\in C([0,{T_0}];
C(\overline\Omega))$ guaranteed by Proposition \ref{cs}, one gets $\sqrt\varrho v, \sqrt\varrho\vartheta\in C([0,{T_0}]; L^2)$. Thanks to these, similarly to the proof
of (ii) of Proposition \ref{CONTINUITY} (since $\psi$ has the same properties as those of $\varphi$), one can then show that $\sqrt\rho u, \sqrt\rho\theta\in C([0,{T_0}]; L^2)$.

{\textbf{(ii) Uniqueness.}}
Let $(\hat\rho, \hat u,\hat\theta)$ and $\check\rho, \check u, \check\theta)$ be two solutions to system (\ref{a})--(\ref{c}), subject to (\ref{d})--(\ref{e}), in $\Omega\times(0,{T_0})$, with the same initial
data $(\rho_0, u_0, \theta_0)$. Let $(\hat\varphi,\hat\psi,\hat\varrho,\hat{v},\hat{\vartheta},\hat{A},\hat{J})$ and $(\check\varphi,\check\psi,\check\varrho,\check{v},\check{\vartheta},\check{A},\check{J})$ be the corresponding quantities defined as before and denote
 $$
 (v,\vartheta,A,J)=(\hat{v},\hat{\vartheta},\hat{A},\hat{J})-(\check{v},\check{\vartheta},\check{A},\check{J}).
 $$
Then, $(\hat{v},\hat{\vartheta},\hat{A},\hat{J})$ and $(\check{v},\check{\vartheta},\check{A},\check{J})$
have the regularities \eqref{REGRVT1} and \eqref{REGRVT3}, satisfy
system \eqref{bg}--\eqref{bj} a.e.\,in $\Omega\times(0,T_0)$, and fulfill the initial-boundary conditions
\eqref{bk}--\eqref{bl}. One can check by direct calculations that $(v,\vartheta,A,J)$ satisfies:
\begin{align}
\varrho_{0}\hat{J}v_{t}&-\mu\dive_{\hat{A}}(\nabla_{\hat{A}}v)-(\mu+\lambda)\nabla_{\hat{A}}(\dive_{\hat{A}}v)
 = -\varrho_{0}J\check{v}_{t}\nonumber\\
&+\mu\dive_{\hat{A}}(\nabla_{A}\check{v})+\mu\dive_{A}(\nabla_{\check{A}}\check{v})+(\mu+\lambda)\nabla_{\hat{A}}
(\dive_{A}\check{v})+(\mu+\lambda)\nabla_{A}(\dive_{\check{A}}\check{v})\nonumber\\
&-R\nabla_{\hat{A}}(\varrho_{0}\hat{J}\vartheta+\varrho_{0}J\check{\vartheta})-R\nabla_{A}(\varrho_{0}
 \check{J}\check{\vartheta}), \label{bm}\\
c_{v}\varrho_{0}\hat{J}\vartheta_{t}&-\kappa\dive_{\hat{A}}(\nabla_{\hat{A}}\vartheta)=
-c_{v}\varrho_{0}J\check{\vartheta}_{t}+\kappa\dive_{\hat{A}}(\nabla_{A}\check{\vartheta})+\kappa\dive_{A}
(\nabla_{\check{A}}\check{\vartheta})\nonumber \\
&-R\varrho_{0}\Big(\hat{J}\hat{\vartheta}\dive_{\hat{A}}v+\hat{J}\hat{\vartheta}\dive_{A}\check{v}+\hat{J}\vartheta
\dive_{\check{A}}\check{v}+J\check{\vartheta}\dive_{\check{A}}\check{v}\Big)\nonumber\\
&+\frac{\mu}{2}\left(\nabla_{\hat{A}}^{i}\hat{v}_{j}+\nabla_{\hat{A}}^{j}\hat{v}_{i}+\nabla_{\check{A}}^{i}
\check{v}_{j}+\nabla_{\check{A}}^{j}\check{v}_{i}\right)\left(\nabla_{\hat{A}}^{i}v_{j}+\nabla_{\hat{A}}^{j}v_{i}
+\nabla_{A}^{i}\check{v}_{j}+\nabla_{A}^{j}\check{v}_{i} \right) 	\nonumber  \\
&+\lambda\left( \dive_{\hat{A}}\hat{v}+\dive_{\check{A}}\check{v}\right)
\left(\dive_{\hat{A}}v+\dive_{A}\check{v} \right), 	\label{bn} \\
A_{t}&+\nabla_{\hat{A}}\hat{v}A+\nabla_{\hat{A}} v\check{A}+\nabla_{A}\check{v}\check{A}=0, 	\label{bo}\\
J_{t}&+\dive_{\hat{A}}\hat{v}J+\dive_{\hat{A}}v\check{J}+\dive_{A}\check{v}\check{J}=0. 	\label{bp}
\end{align}
For any vector field $W$ and function $f$ such that either $W|_{\partial\Omega}=0$ or $f|_{\partial\Omega}=0$,
by Lemma \ref{cv}, it follows from integration by parts that
\begin{eqnarray}
  &&\int_\Omega\frac{1}{\hat J}\nabla_{\hat A}f\cdot W dy=\int_\Omega\frac{1}{\hat J}
  \hat a_{il}\partial_lf W_i dy \nonumber\\
  &=&-\int_\Omega\left(\partial_l\left(\frac{\hat a_{il}}{\hat J}\right)W_i+\frac{1}{\hat J} \hat a_{il} \partial_lW_i\right) fdy=-\int_\Omega\frac{1}{\hat J}\text{div}_{\hat A}Wfdy.\label{INBYP}
\end{eqnarray}

\textbf{Step 1. Energy inequalities.}
Multiplying \eqref{bm} with $\frac{v}{\hj}$ and using (\ref{INBYP}), one gets
\begin{eqnarray}
&&\fr\fracc\norm{\sqrt{\varrho_{0}}v}^{2}+\mu\left\|\frac{\nabla_{\hat A}v}{\sqrt{\hat J}}\right\|
+(\mu+\lambda)\left\|\frac{\text{div}_{\hat A}v}{\sqrt{\hat J}}\right\|  	\nonumber\\
&=&-\inn\varrho_{0}J\cv_{t}\frac{v}{\hj}\dx-\inn\frac{1}{\hat J}\big[\mu\na\cv:\nha v+\lr{\mu+\lambda}
  {\da\cv}\text{div}_{\hat A}v\big]\dx\nonumber\\	
&&+\inn\frac{1}{\hat J}\big[\mu\da\lr{\nca\cv}\cdot v +\lr{\mu+\lambda}\na{\lr{\dca\cv}}\cdot v\big]\dx-R\inn\na\lr{\varrho_{0}\cj\ct}
  \frac{v}{\hj}\dx\nonumber\\
&& +R\inn\lr{\varrho_{0}\hj\vartheta+\varrho_{0}J\ct}\frac{\text{div}_{\hat A}v}{\hj}\dx=:\sum_{i=1}^5N_i. 	\label{cf}
\end{eqnarray}
By Proposition \ref{cs}, it follows
\begin{equation*}
  \sup_{0\leq t\leq{T_0}}\left(\left\|\frac{1}{\hat J}\right\|_\infty+\|(\hat J, \check J,\hat A)\|_\infty+\|(\nabla\check A, \nabla\check J)\|_3\right)\leq C.
\end{equation*}
Thanks to this, terms $N_i, i=1, 2, \cdots, 5,$ are estimated by the H\"older, Gagliardo-Nerenberg, Sobolev,
and Young inequalities
as follows:
\begin{eqnarray*}
N_1&\leq& C\norm{J}\norm{\varrho_{0}\cv_{t}}_{3}\norm{v}_{6} \leq C\norm{J}\norm{\sqrt{\varrho_{0}}\cv_{t}}^{\fr}\norm{\cv_{t}}_6^{\fr}\norm{\nabla v}\\
	& \leq& \frac{\mu}{8}\normm{\frac{\nha v}{\sqrt{\hj}}}^{2}+C\norm{\sqrt{\varrho_{0}}\cv_{t}}\norm{\nabla\cv_{t}}\norm{J}^{2},\\
N_2&=&-\inn\frac{1}{\hat J}\big[\mu\na\cv\cdot\nha v+\lr{\mu+\lambda}
  {\da\cv}\text{div}_{\hat A}v\big]\dx\\
  &\leq& \frac{\mu}{8}\normm{\frac{\nha v}{\sqrt{\hj}}}^{2}+C\norm{\nabla\cv}_{\infty}^{2}\norm{A}^{2}
  \leq \frac{\mu}{8}\normm{\frac{\nha v}{\sqrt{\hj}}}^{2}+C\norm{\nabla^2\cv}_q^{2}\norm{A}^{2},\\
N_3&=&\inn\frac{1}{\hat J}\big[\mu\da\lr{\nca\cv}\cdot v +\lr{\mu+\lambda}\na{\lr{\dca\cv}}\cdot v\big]\dx\\
	&=&\inn\frac{1}{\hat J}\big[\mu  a_{kl}\ptl\lr{\check{a}_{km}\ptm\cv}\cdot\ v+(\mu+\lambda){a}_{il}\ptl\lr{\check a_{km}\ptm\cv_{k}} v_{i}\big]\dx\\
&\leq&C\int_\Omega|A|(|\check A||\nabla^2\check v|+|\nabla\check A||\nabla\check v|)|v|dy \\
&\leq& C\|A\|(\|\nabla^2\check v\|_3+\|\nabla\check A\|_3\|\nabla\check v\|_\infty)\|v\|_6 \\
 &\leq& \frac{\mu}{8}\normm{\frac{\nha v}{\sqrt{\hj}}}^{2}+C \norm{\nabla^{2}\cv}^{2}_{q} \norm{A}^{2},\\
 N_4&=& -R\inn a_{il}\ptl\lr{\varrho_{0}\cj\ct}\frac{v_{i}}{\hj}\dx\\
	& \leq& C\inn\nn{A}\lr{\nn{\nabla\varrho_{0}}\nn{\cj}\nn{\ct}+\varrho_{0}\nn{\nabla\cj}\nn{\ct}+\varrho_{0}\nn{\cj}\nn{\nabla\ct}}\frac{\nn{v}}{\cj}\dx\\
	& \leq& C\norm{A}\lr{\norm{\nabla\varrho_{0}}_{3} \norm{\ct}_{\infty}+\norm{\nabla\cj}_{3}\norm{\ct}_{\infty}+ \norm{\nabla\ct}_{3}}\norm{v}_{6}\\
	& \leq&\frac{\mu}{8}\normm{\frac{\nha v}{\sqrt{\hj}}}^{2}+C\norm{\nabla\ct}\norm{\nabla^{2}\ct}\norm{A}^{2},\\	N_5&\leq&  \frac{\mu+\lambda}{4}\normm{\frac{\dha v}{\sqrt{\hj}}}^{2}+C\lr{ \norm{\sqrt{\varrho_{0}}\vartheta}^2
+ \norm{\ct}^{2}_{\infty}\norm{J}^{2}}\\
	& \leq& \frac{\mu+\lambda}{4}\normm{\frac{\dha v}{\sqrt{\hj}}}^{2}+C\lr{\norm{\sqrt{\varrho_{0}}\vartheta}^{2}+\norm{\nabla\ct}\norm{\nabla^{2}\ct}\norm{J}^{2}}.
\end{eqnarray*}
 Substituting these estimates into \eqref{cf} yields
 \begin{eqnarray}
 \nonumber
 &&\fracc\norm{\sqrt{\varrho_{0}}v}^{2}+\mu\normm{\frac{\nha v}{\sqrt{\hj}}}^{2}\\
 &\leq& C\norm{\sqrt{\varrho_{0}}\vartheta}^{2}
+C\lr{\norm{\sqrt{\varrho_{0}}\cv_{t}}\|\nabla\check v_t\|+\norm{\nabla\ct}\norm{\nabla^{2}\ct}+\norm{\nabla^{2}\cv}^{2}_{q}}\norm{(A,J)}^{2}. \label{cg}
 \end{eqnarray}
  Testing \eqref{bn} with $\frac{\vartheta}{\hj}$ and using (\ref{INBYP}), one gets
  \begin{align}
  \nonumber&\frac{c_{v}}{2}\fracc\normm{\sqrt{\varrho_{0}}\vartheta}^{2}+\kappa\normm{\frac{\nha\vartheta}{\sqrt{\hj}}}^{2}\\
  \nonumber=&-c_{v}\inn\varrho_{0}J\ct_{t}\frac{\vartheta}{\hj}\dx-\kappa\inn\nabla_{A}\ct\cdot \frac{\nha\vartheta}{\hj}\dx
  +\kappa\inn \dive_{A}\lr{\nca\ct}\frac{\vartheta}{\hj}\dx\\
  \nonumber&-R\inn\varrho_{0}\lr{\hj\htt\dha v+\hj\htt\da\cv+\hj\vartheta\dca\cv+J\ct\dca\cv}\frac{\vartheta}{\hj}\dx\\
  \nonumber&+\frac{\mu}{2}\inn \left(\nabla_{\hat{A}}^{i}\hat{v}_{j}+\nabla_{\hat{A}}^{j}\hat{v}_{i}
  +\nabla_{\check{A}}^{i}\check{v}_{j}+\nabla_{\check{A}}^{j}\check{v}_{i}\right)\left(\nabla_{\hat{A}}^{i}v_{j}
  +\nabla_{\hat{A}}^{j}v_{i}+\nabla_{A}^{i}\check{v}_{j}+\nabla_{A}^{j}\check{v}_{i} \right)\\
  \nonumber& +\lambda\int_\Omega\left( \dive_{\hat{A}}\hat{v}+\dive_{\check{A}}\check{v}\right)
  \left(\dive_{\hat{A}}v+\dive_{A}\check{v} \right) \frac{\vartheta}{\hj}\dx=:\sum_{i=1}^6O_i.
  \end{align}
Similar to $N_1, N_2,$ and $N_3$, one has the following estimates for $O_1, O_2,$ and $O_3$:
\begin{eqnarray*}
	O_1&\leq& \frac{\kappa}{8}\normm{\frac{\nha \vartheta}{\sqrt{\hj}}}^{2}+C\norm{\sqrt{\varrho_{0}}\ct_{t}}\norm{\nabla\ct_{t}}\norm{J}^{2},\\
	O_2&\leq& \frac{\kappa}{8}\normm{\frac{\nha \vartheta}{\sqrt{\hj}}}^{2}+C\norm{\nabla\ct}_{\infty}^{2}\norm{A}^{2},\\
	O_3&\leq& \frac{\kappa}{8}\normm{\frac{\nha \vartheta}{\sqrt{\hj}}}^{2}+C\lr{\norm{\nabla^{2}\ct}_{3}^{2}+\norm{\nabla\ct}_{\infty}^{2}}\norm{A}^{2}.
\end{eqnarray*}
By the H\"older, Sobolev, and Young inequalities, one deduces
\begin{eqnarray*}
O_4&\leq &C\int_\Omega \rho_0|\vartheta|\big(|\hat\vartheta||\text{div}_{\hat A}v|+|\hat\vartheta||A||\nabla\check v|
+|\vartheta||\nabla\check v|+|J||\check\vartheta||\nabla\check v|\big)dy\\
&\leq& C\norm{\sqrt{\varrho_{0}}\vartheta}\Bigg( || \htt||_{\infty}\normm{\frac{\nha v}{\sqrt{\hj}}} +\| \htt\|_{\infty}\norm{\nabla\cv}_{\infty}\norm{A}\\ &&+\norm{\sqrt{\varrho_{0}}\vartheta}\norm{\nabla\cv}_{\infty}+\norm{J}\norm{ \ct}_{\infty}
\norm{\nabla\cv}_{\infty}\Bigg)\\
&\leq&C\lr{\| \htt\|_{\infty}+\| \ct\|_{\infty}}\norm{\nabla\cv}_{\infty}
\norm{\sqrt{\varrho_{0}}\vartheta}\lr{\norm{A}+\norm{J}} \\
&&+C\| \htt\|_{\infty}\normm{\frac{\nha v}{\sqrt{\hj}}}\norm{\sqrt{\varrho_{0}}\vartheta}+C\norm{\nabla\cv}_{\infty}\norm{\sqrt{\varrho_{0}}\vartheta}^{2},\\
O_5+O_6&\leq& C\inn\lr{\nn{\nabla\hv}+\nn{\nabla\cv}}\lr{\nn{\nha v}+\nn{A}\nn{\nabla\cv}}|\vartheta|\dx\\
	&\leq& C\norm{\nabla(\hv,\cv)}_{3}\normm {\frac{\nha v}{\sqrt{\hj}}}\norm{\vartheta}_{6}+C\norm{\nabla(\hv,\cv)}_{6}^{2}\norm{A}\norm{\vartheta}_{6}\\
	&\leq& \frac{\kappa}{8}\normm{\frac{\nha \vartheta}{\sqrt{\hj}}}^{2}+C\lr{\norm{\nabla(\hv,\cv)}_{3}^{2}\normm {\frac{\nha v}{\sqrt{\hj}}}^{2}+\norm{\nabla^{2}(\hv,\cv)}^{4}\norm{A}^{2}}.
	\end{eqnarray*}
Thus, combining the above yields
\begin{eqnarray}
&&c_{v}\fracc\norm{\sqrt{\varrho_{0}}\vartheta}^{2}+\kappa\normm{\frac{\nha \vartheta}{\sqrt{\hj}}}^{2}	\nonumber \\
&\leq&C\norm{\sqrt{\varrho_{0}}\ct_{t}}\norm{\nabla\ct_{t}}\norm{J}^{2}+C\lr{\norm{\nabla\ct}_{\infty}^{2}+\norm{\nabla^{2}
\ct}_{3}^{2}+\norm{\nabla^{2}(\hv,\cv)}^{4}}\norm{A}^{2}	\nonumber\\
&&+C\norm{\nabla(\hv,\cv)}_{3}^{2}\normm{\frac{\nha
v}{\sqrt{\hj}}}^{2}+C|| \htt||_{\infty}\normm{\frac{\nha
v}{\sqrt{\hj}}}\norm{\sqrt{\varrho_{0}}\vartheta}+C\norm{\nabla\cv}_{\infty}\norm{\sqrt{\varrho_{0}}\vartheta}^{2}\nonumber\\
&&+C\norm{\nabla\cv}_{\infty}\lr{|| \htt||_{\infty}
+|| \ct||_{\infty}}\norm{\sqrt{\varrho_{0}}\vartheta}\lr{\norm{A}+\norm{J}}.
\label{ch}
\end{eqnarray}

\textbf{Step 2. Growth estimates.}
We proceed to consider the growth estimates of $A$ and $J$. Testing (\ref{bo}) and (\ref{bp}), respectively, with $A$ and $J$, and summing the resultant up, one obtains after some straightforward computations that
 \begin{align}
 	&\fracc(\norm{A}^{2}+\|J\|^2)\leq C\norm{\nabla(\hv,\cv)}_{\infty}(\norm{A}^{2}+\|J\|^2)
 +C\norm{\nabla v}(\norm{A}+\|J\|)\label{dh}
 \end{align}
and, thus,
\begin{align}
 	&\fracc\sqrt{\norm{A}^{2}+\|J\|^2}\leq C\norm{\nabla(\hv,\cv)}_{\infty}\sqrt{\norm{A}^{2}+\|J\|^2}
 +C\norm{\nabla v}.\label{dh-1}
 \end{align}
Since $\hat v, \check v\in L^1(0,{T_0}; W^{2,q})\cap L^\infty(0,{T_0}; H^1)$ and $W^{1,q}\hookrightarrow L^\infty$ for $q\in(3,6)$, one has
$\norm{\nabla(\hv,\cv)}_{\infty}\in L^{1}((0,{T_0}))$ and $\norm{\nabla v}\in L^{\infty}((0,{T_0}))$. Thanks
to these and applying the Gr\"onwall inequality to (\ref{dh-1}), one deduces
\begin{eqnarray}
  \sqrt{\norm{A}^{2}+\|J\|^2} \leq  Ce^{C\int_{0}^{t}\norm{\nabla(\hv,\cv)}_{\infty}\ms}\int_0^t\|\nabla v\| ds
  \leq Ct,\qquad\forall t\in(0,{T_0}).
\label{ci}
\end{eqnarray}

Recalling that
\begin{eqnarray*}
  \sqrt{\rho_0}\check v_t\in L^2(0,T_0; L^2),\quad \sqrt t \check v_t\in L^2(0,T_0; H^1),\quad \check\vartheta\in L^2(0,T_0; H^1),\\
\quad\sqrt t\check \vartheta\in L^2(0,{T_0}; H^2),\quad \sqrt t\check v\in L^2(0,{T_0}; W^{2,q}),
\end{eqnarray*}
one gets
$$
\omega_{1}(t)\triangleq t\lr{\norm{\sqrt{\varrho_{0}}\cv_{t}}\|\nabla\check v_t\|+\norm{\nabla\ct}\norm{\nabla^{2}\ct}+\norm{\nabla^{2}\cv}^{2}_{q}}\in L^{1}((0,{T_0})).
$$
Since $\sqrt{\rho_0}\vartheta,\sqrt{\varrho_0}v\in C([0,T_0]; L^2)$ (guaranteed by Proposition \ref{cn}) and $\sqrt{\rho_0}v|_{t=0}=0$,
integrating \eqref{cg} with respect to $t$ and using (\ref{ci}) yield
$$
\norm{\sqrt{\varrho_{0}}v}^{2}(t)+\mu\int_0^t\normm{\frac{\nha v}{\sqrt{\hj}}}^{2}ds\leq
 Ct+Ct\int_0^t \omega_{1}ds\leq Ct,\quad\forall t\in(0,T_0).
$$
Combining this with (\ref{ci}) leads to
\begin{equation}
 	\label{cj}
 \sup_{0\leq t\leq {T_0}}\lr{\norm{A}+\norm{J}+\norm{\sqrt{\varrho_{0}}v}^{2}}+\int_{0}^{t}\norm{\nabla v}^{2}\ms\leq Ct,\qquad\forall t\in(0,{T_0}).
\end{equation}

\textbf{Step 3. Singular $t$-weighted energy inequalities and uniqueness.}
 Multiplying \eqref{dh} by $t^{-\frac{3}{2}}$ yields
 \begin{eqnarray*}
 &&\fracc\left(\frac{\|A\|^2}{t^{\frac32}}+\frac{\|J\|^2}{t^\frac32}\right)
 +\frac32\left(\frac{\|A\|^2}{t^{\frac52}}+\frac{\|J\|^2}{t^\frac52}\right)\\
 &\leq& C\norm{\nabla(\hv,\cv)}_{\infty}\left(\frac{\|A\|^2}{t^{\frac32}}+\frac{\|J\|^2}{t^\frac32}\right)
 +C\frac{\norm{\nabla v}}{t^\frac14}\left(\frac{\|A\|}{t^{\frac54}}+\frac{\|J\|}{t^\frac54}\right)\\
 &\leq&
 \frac12\left(\frac{\|A\|^2}{t^{\frac52}}+\frac{\|J\|^2}{t^\frac52}\right)+\frac{C}{\sqrt t}\left\|\frac{\nabla_{\hat A} v}{\sqrt{\hat J}}\right\|^2\\
  &&+C \norm{\nabla(\hv,\cv)}_{\infty}\left(\frac{\|A\|^2}{t^{\frac32}}+\frac{\|J\|^2}{t^\frac32}\right)
 \end{eqnarray*}
and, thus,
 \begin{equation}
\fracc\left(\frac{\|A\|^2}{t^{\frac32}}+\frac{\|J\|^2}{t^\frac32}\right)
\leq
 C \norm{\nabla(\hv,\cv)}_{\infty}\left(\frac{\|A\|^2}{t^{\frac32}}+\frac{\|J\|^2}{t^\frac32}\right)
 +\frac{C}{\sqrt t}\left\|\frac{\nabla_{\hat A} v}{\sqrt{\hat J}}\right\|^2. \label{ck}
 \end{equation}
Multiplying (\ref{cg}) with $\frac{1}{\sqrt t}$ and recalling the definition of $\omega_1(t)$ yield
\begin{eqnarray}
&&\fracc\lr{\frac{\norm{\sqrt{\varrho_{0}}v}^{2}}{\sqrt t}}+ \frac{\norm{\sqrt{\varrho_{0}}v}^{2}}{2t^{\frac{3}{2}}}+\frac{\mu}{\sqrt t}\normm{\frac{\nha v}{\sqrt{\hj}}}^{2} 	\nonumber\\
&\leq& C\left(\frac{1}{\sqrt t}+\omega_1(t)\right)\lr{\norm{\sqrt{\varrho_{0}}\vartheta}^{2}
 +\frac{\|A\|^2}{t^{\frac32}}+\frac{\|J\|^2}{t^\frac32}}.	\label{cl}
 \end{eqnarray}
Denote
\begin{eqnarray*}
&&\omega_{21}(t)\triangleq t^{\frac{3}{2}}\norm{\sqrt{\varrho_{0}}\ct_{t}}\norm{\nabla\ct_{t}},\quad\omega_{22}(t)\triangleq t^{\frac{3}{2}}\lr{\norm{\nabla\ct}^{2}_{\infty}+\norm{\nabla^{2}\ct}_{3}^{2}+\norm{\nabla^{2}(\hv,\cv)}^{4}}\\
&&\omega_{23}(t)\triangleq \sqrt{t}\norm{\nabla(\hv,\cv)}_{3}^{2},\quad
\omega_{24}(t) \triangleq  \sqrt t\norm{ \htt}^{2}_{\infty},\quad
\omega_{25}(t) \triangleq  \norm{\nabla\cv}_{\infty},\\
&&\omega_{26}(t)\triangleq t^{\frac{3}{4}}\norm{\nabla\cv}_{\infty}||( \htt, \ct)||_{\infty}.
\end{eqnarray*}
Recalling the regularities of $(\hat v, \hat\vartheta)$ and $(\check v, \check\vartheta)$, we have
\begin{align}
&\omega_{21}(t)\leq  \norm{\sqrt t\sqrt{\varrho_{0}}\ct_{t}}\norm{t\nabla\ct_{t}}\in L^1((0,{T_0})), \nonumber\\
&\omega_{23}(t)  \leq
 C\norm{\nabla(\hv,\cv)}\|\sqrt{t}\nabla^{2}(\hv,\cv)\|\in L^{\infty}((0,{T_0})), \nonumber \\
&\omega_{24}(t) \leq C\|\nabla\htt\|\|\sqrt{t}\nabla^{2}\htt\|\in L^1((0,{T_0})), \quad \omega_{25}(t)=\norm{\nabla\cv}_{\infty}\in L^1((0,{T_0})),\nonumber\\
&\omega_{26}(t) \leq C\|\sqrt{t}\nabla^{2}\cv\|_{q}||\nabla(\htt,\ct)\|^{\fr}\|\sqrt{t}\nabla^{2}(\htt,\ct)\|^{\fr}\in L^1((0,{T_0})). \nonumber
\end{align}
For $\omega_{22}$, by Proposition \ref{cs}, it follows from the Gagliardo-Nirenberg inequality that
\begin{eqnarray*}
  \omega_{22}(t)&\leq &C t^{\frac{3}{2}}(\norm{\nabla\check\theta}^{2}_{\infty}+\norm{\nabla^{2}\check\theta}_{3}^{2})+ C\sqrt t \norm{\nabla^{2}(\hv,\cv)}^{2}\|\sqrt{t}\nabla^{2}(\hv,\cv)\|^{2} \\
  &\leq&C\|\sqrt t\nabla^2\check\theta\|\|t\nabla^2\check\theta\|_6+
  C\sqrt t \norm{\nabla^{2}(\hv,\cv)}^{2}||\sqrt{t}\nabla^{2}(\hv,\cv)||^{2}\in L^1((0,{T_0})).
\end{eqnarray*}
In terms of $\omega_{2i}, i=1,2,\cdots,6$, one gets from \eqref{ch} by the Young inequality that
\begin{eqnarray*}
&&c_{v}\fracc\norm{\sqrt{\varrho_{0}}\vartheta}^{2}+\kappa\normm{\frac{\nha\vartheta}{\sqrt{\hj}}}^{2}\\
&\leq& C\lr{\omega_{21}(t)\frac{\|J\|^2}{t^\frac32}+\omega_{22}(t)\frac{\|A\|^2}{t^{\frac32}}
+\frac{\omega_{23}(t)}{\sqrt t}\normm{\frac{\nha v}{\sqrt{\hj}}}^{2}}	\nonumber\\
&&+C\sqrt{\omega_{24}(t)}\left(\frac{1}{\sqrt t}\normm{\frac{\nha v}{\sqrt{\hj}}}^{2}\right)^\frac12||\sqrt{\varrho_{0}}\vartheta|| +C\omega_{25}(t)||\sqrt{\varrho_{0}}\vartheta||^{2}\nonumber\\
&&+C\omega_{26}(t)\lr{\norm{\sqrt{\varrho_{0}}\vartheta}^{2}
 +\frac{\|A\|^2}{t^{\frac32}}+\frac{\|J\|^2}{t^\frac32}},
\end{eqnarray*}
from which, by the Young inequality and recalling that $\omega_{23}\in L^\infty((0,{T_0}))$, one gets
\begin{eqnarray}
&&c_{v}\fracc\norm{\sqrt{\varrho_{0}}\vartheta}^{2}+\kappa\normm{\frac{\nha\vartheta}{\sqrt{\hj}}}^{2}\nonumber\\
&\leq&C(\omega_{23}(t)+1)\frac{1}{\sqrt t}\normm{\frac{\nha v}{\sqrt{\hj}}}^{2}+C \omega_2(t)\lr{\norm{\sqrt{\varrho_{0}}\vartheta}^{2}
 +\frac{\|A\|^2}{t^{\frac32}}+\frac{\|J\|^2}{t^\frac32}}\nonumber\\
 &\leq& \frac{ C}{\sqrt t}\normm{\frac{\nha v}{\sqrt{\hj}}}^{2}+C \omega_2(t)\lr{\norm{\sqrt{\varrho_{0}}\vartheta}^{2}
 +\frac{\|A\|^2}{t^{\frac32}}+\frac{\|J\|^2}{t^\frac32}}, 	\label{cm}
\end{eqnarray}
 where
 $$
 \omega_2:=\omega_{21}+\omega_{22}+\omega_{24}+\omega_{25}+\omega_{26}\in L^1((0,{T_0})).
 $$

Multiplying \eqref{ck} and \eqref{cm} with a small positive number $\zeta$
and adding the resultants with \eqref{cl}, one obtains
 \begin{eqnarray*}
 \frac{d}{dt}\left[\frac{\norm{\sqrt{\varrho_{0}}v}^{2}}{\sqrt t}+\zeta\left(\frac{\|A\|^2}{t^{\frac32}}+ \frac{\|J\|^2}{t^\frac32}+ c_{v}\norm{\sqrt{\varrho_{0}}\vartheta}^{2}\right)\right] +\frac{\mu}{2\sqrt t}\normm{\frac{\nha v}{\sqrt{\hj}}}^{2}+\kappa\zeta\normm{\frac{\nha\vartheta}{\sqrt{\hj}}}^{2}\\
 	\leq C\left(\frac{1}{\sqrt t}+\omega_1(t)+\omega_2(t)+\|\nabla(\hat v,\check v)\|_\infty\right)\left(\frac{\|A\|^2}{t^{\frac32}}+ \frac{\|J\|^2}{t^\frac32}+ \norm{\sqrt{\varrho_{0}}\vartheta}^{2}\right).
 \end{eqnarray*}
By Proposition \ref{cn} and recalling (\ref{cj}), it follows that
\begin{eqnarray*}
  \lim_{t\rightarrow0}\left[\frac{\norm{\sqrt{\varrho_{0}}v}^{2}}{\sqrt t}+\zeta\left(\frac{\|A\|^2}{t^{\frac32}}+ \frac{\|J\|^2}{t^\frac32}+ c_{v}\norm{\sqrt{\varrho_{0}}\vartheta}^{2}\right)\right](t)=0.
\end{eqnarray*}
Thanks to this and by the Gr\"onwall inequality, one gets
 \begin{align*}
 &
 \left(\frac{\norm{\sqrt{\varrho_{0}}v}^{2}}{\sqrt t}+\frac{\|A\|^2}{t^{\frac32}}+ \frac{\|J\|^2}{t^\frac32}+ \norm{\sqrt{\varrho_{0}}\vartheta}^{2}\right)(t)+\int_{0}^{t}\lr{\frac{1}{\sqrt t}\normm{\frac{\nha v}{\sqrt{\hj}}}^{2}+ \normm{\frac{\nha\vartheta}{\sqrt{\hj}}}^{2}}\mathrm{d}\tau=0,
 \end{align*}
which implies $A=J=v=\vartheta=0$.

Recalling that $\hat\varrho=\hat J\rho_0$ and $\check\varrho=\check J\rho_0$,
it follows $\hat\varrho=\check\varrho$. Noticing that $\partial_t\hat\varphi(y,t)=\hat v(y,t)$ and $\partial_t\check\varphi(y,t)=\check v(y,t)$,
one has $\hat\varphi=\check\varphi$ and further that $\hat\psi=\check\psi$. Then, it follows
$$\hat u(x,t)=\hat v(\hat\psi(x,t),t)=\check v(\check\psi(x,t),t)=\check u(x,t),$$
that is $\hat u\equiv\check u$. Similarly, one has $\hat\theta\equiv\check\theta$ and $\hat\rho\equiv\check\rho$. This proves the uniqueness.
\end{proof}

\section{Appendix}

In this appendix, we give the proof of Proposition \ref{cs} and Proposition \ref{CONTINUITY}, as well as a lemma used in
(\ref{INBYP}).

\begin{proof}[\textbf{Proof of Proposition \ref{cs}}]
  Solving (\ref{LG6}) yields
$$
J(y,t)=e^{-\int_0^t\text{div} u(\varphi(y,s),s)ds}
$$
and, thus,
\begin{equation}
  \frac1{C_*}\leq J(y,t)\leq C_*,\quad\forall (y,t)\in\Omega\times[0,{T_0}], \label{LG7}
\end{equation}
where $C_*:=e^{\int_0^{{T_0}}\|\text{div} u\|_\infty dt}$. Since $\text{det}\nabla\psi(x,t)=J(\psi(x,t),t)$, it
holds that
\begin{eqnarray*}
\|g(\varphi(\cdot,t))\|_\alpha& =&\left(\int_\Omega|g(\varphi(y,t))|^\alpha dy\right)^\frac1\alpha   =\left(\int_\Omega|g(x)|^\alpha|\text{det}\nabla\psi(x,t)|dx\right)^\frac1\alpha\\
&=&\left(\int_\Omega|g(x)|^\alpha J(\psi(x,t),t)dx\right)^\frac1\alpha,
\end{eqnarray*}
from which, by (\ref{LG7}), one gets
$$
  C_*^{-\frac1\alpha}\|g\|_\alpha\leq
  \|g(\varphi(\cdot,t))\|_\alpha\leq C_*^\frac1\alpha\|g\|_\alpha,\quad \alpha\in[1,\infty).
$$
Letting $\alpha\rightarrow\infty$ in the above leads to the estimate for $\alpha=\infty$. Therefore,
\begin{equation}
  \|g\|_\alpha\simeq
  \|g(\varphi(\cdot,t))\|_\alpha,\quad \alpha\in[1,\infty].\label{LG14}
\end{equation}
 Applying the above to $g(\psi(x,t))$ leads to
\begin{equation}
 \|g\|_\alpha\simeq
  \|g(\psi(\cdot,t))\|_\alpha,\quad \alpha\in[1,\infty].\label{LG14'}
\end{equation}
Combining (\ref{LG14}) with (\ref{LG14'}) leads to (\ref{CS5}).

It follows from (\ref{LG4}) that
$$
\partial_t|B|^2=2B:(B\nabla u(\varphi,t))\leq C\|\nabla u\|_\infty|B|^2
$$
and, thus,
$$
|B|^2(y,t)\leq 3e^{C\int_0^t\|\nabla u\|_\infty ds},\quad\forall(y,t)\in\Omega\times[0,{T_0}].
$$
Therefore
\begin{equation}
  \label{LG8}
  \sup_{0\leq t\leq{T_0}}\|B\|_\infty(t)\leq Ce^{C\int_0^{T_0}\|\nabla u\|_\infty dt}\leq C.
\end{equation}
Similarly, one gets from (\ref{LG5}) that
\begin{equation}
  \sup_{0\leq t\leq T}\|A\|_\infty\leq C. \label{LG9}
\end{equation}
Thanks to (\ref{LG7}), (\ref{LG8}), and (\ref{LG9}), it follows from (\ref{LG5})--(\ref{LG6}) and (\ref{LG14}) that
\begin{eqnarray}
  \sup_{0\leq t\leq{T_0}}\|(J_t, A_t, B_t)\|\leq C\sup_{0\leq t\leq{T_0}}\|\nabla u(\varphi(\cdot,t),t)\|\leq
  C\sup_{0\leq t\leq{T_0}}\|\nabla u\|\leq C.\label{LG9-1}
\end{eqnarray}

One gets from (\ref{LG4}) that
\begin{equation*}
  \partial_t\partial_iB=\partial_iB\nabla u(\varphi,t)+B\nabla\partial_lu(\varphi,t)\partial_i\varphi_l
  =\partial_iB\nabla u(\varphi,t)+B\nabla\partial_lu(\varphi,t)b_{il}
\end{equation*}
and, thus, by (\ref{LG8}), it follows that
\begin{eqnarray*}
  \partial_t|\nabla B|^2&=&2\partial_iB:(\partial_iB\nabla u(\varphi,t)+B\nabla\partial_lu(\varphi,t)b_{il})\\
  &\leq&C(\|\nabla u\|_\infty|\nabla B|^2+|\nabla B||B|^2|\nabla^2u(\varphi,t)|)\\
  &\leq&C(\|\nabla u\|_\infty|\nabla B|^2+|\nabla^2u(\varphi,t)||\nabla B|),
\end{eqnarray*}
and further that
$$
\partial_t|\nabla B|^q\leq C(\|\nabla u\|_\infty|\nabla B|^q+|\nabla^2u(\varphi,t)|\nabla B|^{q-1}).
$$
Integrating the above over $\Omega$, it follows from the H\"older inequality and (\ref{LG14}) that
\begin{eqnarray*}
  \frac{d}{dt}\|\nabla B\|_q^q&\leq&C\left(\|\nabla u\|_\infty\|\nabla B\|_q^q+\|\nabla^2u(\varphi,t)\|_q\|\nabla B\|_q^{q-1}\right)\\
  &\leq&C\left(\|\nabla u\|_\infty\|\nabla B\|_q^q+\|\nabla^2u\|_q\|\nabla B\|_q^{q-1}\right)\\
\end{eqnarray*}
and, thus,
$$
\frac{d}{dt}\|\nabla B\|_q\leq C\left(\|\nabla u\|_\infty\|\nabla B\|_q+\|\nabla^2u\|_q\right).
$$
Applying the Gr\"onwall inequality to the above yields
\begin{equation}
  \sup_{0\leq t\leq{T_0}}\|\nabla B\|_q\leq Ce^{C\int_0^{T_0}\|\nabla u\|_\infty dt}\int_0^{T_0}\|\nabla^2u\|_q dt\leq C. \label{LG10}
\end{equation}
Similarly, one derives from (\ref{LG5}) and (\ref{LG6}) that
\begin{equation}
  \sup_{0\leq t\leq {T_0}}(\|\nabla A\|_q+\|\nabla J\|_q)\leq C. \label{LG11}
\end{equation}
Conclusion (\ref{CS2}) follows from (\ref{LG7}) and (\ref{LG8})--(\ref{LG11}).

Fix $t_0\in[0,{T_0}]$ and denote
$$
G(y):=g(\varphi(y,t_0)),\quad \forall y\in\Omega.
$$
Then, it is clear that
$$
g(x)=G(\psi(x,t_0)),\quad \forall x\in\Omega.
$$
By direct calculations and recalling the definitions of $a_{ij}(y,t)$ and $b_{ij}(y,t)$, one has
\begin{eqnarray*}
  &&\partial_iG(y)=b_{il}(y,t_0)\partial_lg(\varphi(y,t_0)),\quad\partial_ig(x)
  =a_{il}(\psi(x,t_0),t_0)\partial_lG(\psi(x,t_0)).
\end{eqnarray*}
Therefore, it follows from (\ref{LG8}) and (\ref{LG9}) that
\begin{eqnarray}
  &&|\nabla G(y)|\leq C|\nabla g(\varphi(y,t_0))|,\quad |\nabla g(x)|\leq C|\nabla G(\psi(x,t_0))|.\label{Gg3}
\end{eqnarray}
Thanks to (\ref{Gg3}), it follows from (\ref{LG14}) and (\ref{LG14'}) that
\begin{eqnarray}
  \|\nabla G\|_\alpha \leq C\|\nabla g(\varphi(\cdot,t_0))\|_\alpha\leq C\|\nabla g\|_\alpha,\quad 1\leq\alpha\leq\infty,
  \label{Gg5} \\
  \|\nabla g\|_\alpha\leq C\|\nabla G(\psi(\cdot,t_0))\|_\alpha\leq C\|\nabla G\|_\alpha,\quad 1\leq\alpha\leq\infty.
  \label{Gg6}
\end{eqnarray}
As a result, recalling the definition of $G$, one gets
$$
\|\nabla[g(\varphi(\cdot,t_0))]\|_\alpha \simeq \|\nabla g\|_\alpha,\quad 1\leq\alpha\leq\infty,
$$
which applied to $g(\psi(x,t_0))$ yields
$$
\|\nabla g\|_\alpha\simeq\|\nabla[ g(\psi(\cdot,t_0))]\|_\alpha,\quad 1\leq\alpha\leq\infty.
$$
Therefore (\ref{CS4}) holds.

By direct calculations and recalling the definitions of $a_{ij}(y,t)$ and $b_{ij}(y,t)$, one has
\begin{eqnarray*}
\partial_{ij}^2G(y)&=&\partial_ib_{jl}(y,t_0)\partial_lg(\varphi(y,t_0))+b_{il}(y,t_0)b_{jm}(y,t_0)
\partial_{lm}^2g(\varphi(y,t_0)),\\
\partial_{ij}^2g(x)&=&a_{jm}(\psi(x,t_0),t_0)\partial_ma_{il}(\psi(x,t_0),t_0)\partial_lG(\psi(x,t_0))\\
&&+a_{il}(\psi(x,t_0),t_0)a_{jm}(\psi(x,t_0),t_0)
\partial_{lm}^2G(\psi(x,t_0)).
\end{eqnarray*}
Then, it follows from (\ref{LG8}) and (\ref{LG9}) that
\begin{eqnarray*}
  &&|\nabla^2G(y)|\leq C|\nabla B(y,t_0)||\nabla g(\varphi(y,t_0))|+C|\nabla^2g(\varphi(y,t_0))|,\label{Gg2}\\
  &&|\nabla^2g(x)|\leq C|\nabla A(\psi(x,t_0),t_0)||\nabla G(\psi(x,t_0))|+C|\nabla^2G(\psi(x,t_0))|. \label{Gg4}
\end{eqnarray*}
Thanks to these, it follows from the H\"older and Sobolev inequalities, (\ref{LG10})--(\ref{LG11}), and (\ref{LG14})--(\ref{LG14'}) that: for $1\leq\alpha<3$,
\begin{eqnarray*}
  \|\nabla^2G\|_\alpha&\leq&C\left(\|\nabla B\|_3\|\nabla g(\varphi(\cdot,t_0))\|_{\frac{3\alpha}{3-\alpha}}+\|\nabla^2g(\varphi(
  \cdot, t_0))\|_\alpha\right)\nonumber\\
  &\leq&C\left(\|\nabla B\|_q\|\nabla g\|_{\frac{3\alpha}{3-\alpha}}+\|\nabla^2g\|_\alpha\right)\leq C\|\nabla g\|_{W^{1,\alpha}}, \label{Gg7}\\
  \|\nabla^2g\|_\alpha&\leq&C\left(\|\nabla A(\psi(\cdot,t_0),t_0)\|_3\|\nabla G(\psi(\cdot,t_0))\|_{\frac{3\alpha}{3-\alpha}}+\|\nabla^2G(\psi(\cdot,t_0))\|_\alpha\right)\nonumber\\
  &\leq&C\left(\|\nabla A\|_3\|\nabla G\|_{\frac{3\alpha}{3-\alpha}}+\|\nabla^2G\|_\alpha\right)\leq C\|\nabla G\|_{W^{1,\alpha}};
\end{eqnarray*}
for $\alpha=3$,
\begin{eqnarray*}
  \|\nabla^2G\|_3&\leq&C\left(\|\nabla B\|_q\|\nabla g(\varphi(\cdot,t_0))\|_{\frac{3q}{q-3}}+\|\nabla^2g(\varphi(
  \cdot, t_0))\|_3\right)\nonumber\\
  &\leq&C\left( \|\nabla g\|_{\frac{3q}{q-3}}+\|\nabla^2g\|_3\right)\leq C\|\nabla g\|_{W^{1,3}},\label{Gg8}\\
  \|\nabla^2g\|_3&\leq&C\left(\|\nabla A(\psi(\cdot,t_0),t_0)\|_q\|\nabla G(\psi(\cdot,t_0))\|_{\frac{3q}{q-3}}+\|\nabla^2G(\psi(\cdot,t_0))\|_3\right)\nonumber\\
  &\leq&C\left( \|\nabla A \|_q\|\nabla G\|_{\frac{3q}{q-3}}+\|\nabla^2G\|_3\right)\leq C\|\nabla G\|_{W^{1,3}};
\end{eqnarray*}
and for $3<\alpha\leq q$,
\begin{eqnarray*}
  \|\nabla^2G\|_\alpha&\leq&C\left(\|\nabla B\|_q\|\nabla g(\varphi(\cdot,t_0))\|_\infty+\|\nabla^2g(\varphi(
  \cdot, t_0))\|_\alpha\right)\nonumber\\
  &\leq&C\left(\|\nabla g\|_\infty+\|\nabla^2g\|_\alpha\right)\leq C\|\nabla g\|_{W^{1,\alpha}},\label{Gg9}\\
  \|\nabla^2g\|_\alpha&\leq&C\left(\|\nabla A(\psi(\cdot,t_0),t_0)\|_q\|\nabla G(\psi(\cdot,t_0))\|_\infty+\|\nabla^2G(\psi(\cdot,t_0))\|_\alpha\right)\nonumber\\
  &\leq&C\left(\|\nabla A\|_q\|\nabla G\|_\infty+\|\nabla^2G\|_\alpha\right)\leq C\|\nabla G\|_{W^{1,\alpha}}.
\end{eqnarray*}
Therefore, for any $1\leq\alpha\leq q$, it holds that
$$
\|\nabla^2G\|_\alpha\leq C\|\nabla g\|_{W^{1,\alpha}},\quad \|\nabla^2g\|_\alpha\leq C\|\nabla G\|_{W^{1,\alpha}}.
$$
Thanks to these, by (\ref{Gg5})--(\ref{Gg6}), and recalling the definition of $G$, it follows that
\begin{equation*}
  \|\nabla[g(\varphi(\cdot,t_0))]\|_{W^{1,\alpha}}\leq C\|\nabla g\|_{W^{1,\alpha}},\quad \|\nabla g\|_{W^{1,\alpha}}\leq C \|\nabla[g(\varphi(\cdot,t_0))]\|_{W^{1,\alpha}},
\end{equation*}
for any $1\leq\alpha\leq q$. This proves
$$
\|\nabla[g(\varphi(\cdot,t_0))]\|_{W^{1,\alpha}}\simeq\|\nabla g\|_{W^{1,\alpha}},\quad  \forall\  1\leq\alpha\leq q,
$$
which, applied to $g(\psi(x,t))$, yields further
\begin{equation*}
  \|\nabla g\|_{W^{1,\alpha}}\simeq \|\nabla[g(\psi(x,t_0))]\|_{W^{1,\alpha}},\quad \forall\ 1\leq\alpha\leq q.\label{Gg11}
\end{equation*}
Therefore, (\ref{CS3}) holds.
\end{proof}

\begin{proof}[\textbf{Proof of Proposition \ref{CONTINUITY}}]
By Proposition \ref{cs}, it holds that $\|h(\varphi(\cdot,t),t)\|\leq C\|h\|$ for any $t\in[0,{T_0}]$ and, thus,
$$
\|h(\varphi(\cdot,t),t)\|_{L^\infty(0,{T_0}; L^2)}\leq C\|h\|_{L^\infty(0,{T_0}; L^2)}.
$$
Fix $t_0\in[0,{T_0}]$ and take arbitrary $\varepsilon>0$.
Choose $\xi\in C_c^\infty(\Omega)$ such that
$$
\|\xi-h(\cdot,t_0)\|_{L^2}\leq\varepsilon.
$$

(i) By the H\"older inequality and by Proposition \ref{cs}, one deduces
\begin{eqnarray}
   && \|h(\varphi(\cdot,t),t)-h(\varphi(\cdot,t_0),t_0)\| \nonumber\\
  &\leq&\|h(\varphi(\cdot,t),t)-h(\varphi(\cdot,t),t_0)\|+ \|h(\varphi(\cdot,t),t_0)-\xi(\varphi(\cdot,t))\| \nonumber\\
  &&+\|\xi(\varphi(y,t))-\xi(\varphi(y,t_0))\|+\|\xi(\varphi(y,t_0))-h(\varphi(\cdot,t_0),t_0)\| \nonumber\\
  &\leq&C (\|h(\cdot,t)-h(\cdot,t_0)\|+ \|h(\cdot,t_0)-\xi\| +\|\xi(\varphi(\cdot,t))
  -\xi(\varphi(\cdot,t_0))\|) \nonumber\\
  &\leq&C(\|h(\cdot,t)-h(\cdot,t_0)\|+\varepsilon+\|\xi(\varphi(\cdot,t))-\xi(\varphi(\cdot,t_0))\|) \label{LGG5}
\end{eqnarray}
for any $t\in[0,{T_0}]$.
Recalling (\ref{LG1}) and by Proposition \ref{cs}, it follows from the Gagliardo-Nirenberg and H\"older inequalities
that
\begin{eqnarray}
  &&\left\|\xi(\varphi(\cdot,t))-\xi(\varphi(\cdot,t_0))\right\|
   = \left\|\int_{t_0}^t\nabla\xi(\varphi(\cdot,s))\cdot\partial_t\varphi(\cdot,s)ds\right\|\nonumber\\
   &=&\left\|\int_{t_0}^t
  \nabla\xi(\varphi(\cdot,s))\cdot u(\varphi(\cdot,s),s)ds\right\|
  \leq \left|\int_{t_0}^t\|\nabla\xi(\varphi(\cdot,s))\|\|u\|_\infty ds\right|\nonumber\\
  &\leq& C\left|\int_{t_0}^t\|\nabla\xi\|\|\nabla u\|^\frac12\|\nabla^2u\|^\frac12ds\right|\leq C |t-t_0|^\frac34,\quad\forall
  t\in[0,{T_0}].\label{LGG7}
\end{eqnarray}
Plugging this estimate into (\ref{LGG5}) leads to
\begin{equation}
  \|h(\varphi(\cdot,t),t)-h(\varphi(\cdot,t_0),t_0)\|\leq C(\|h(\cdot,t)-h(\cdot,t_0)\|+\varepsilon+ |t-t_0|^\frac34) \label{LGG6}
\end{equation}
which implies $\|h(\varphi(\cdot,t),t)-h(\varphi(\cdot,t_0),t_0)\|\rightarrow0$ as $t\rightarrow t_0$, for any
$t_0\in[0,{T_0}]$. Therefore, $h(\varphi(\cdot,t),t)\in C([0,{T_0}]; L^2)$.

(ii) Note that
\begin{eqnarray*}
  &&\int_\Omega[h(\varphi(y,t),t)-h(\varphi(y,t_0),t_0)]\chi(y)dy\\
  &=&\int_\Omega[h(\varphi(y,t),t)-h(\varphi(y,t),t_0)]\chi(y)dy+\int_\Omega[h(\varphi(y,t),t_0)-\xi(\varphi(y,t))]\chi(y)dy\\
  &&+\int_\Omega[\xi(\varphi(y,t))-\xi(\varphi(y,t_0))]\chi(y)dy+\int_\Omega[\xi(\varphi(y,t_0))-h(\varphi(y,t_0),t_0)]\chi(y)dy\\
  &=:&R_1+R_2+R_3+R_4.
\end{eqnarray*}
Since $\text{det}\nabla\psi(x,t)=J(\psi(x,t),t)>0$, one deduces
\begin{eqnarray*}
  R_1&=&\int_\Omega[h(\varphi(y,t),t)-h(\varphi(y,t),t_0)]\chi(y)dy\\
  &=&\int_\Omega[h(x,t)-h(x,t_0)]\chi(\psi(x,t))|\text{det}\nabla\psi(x,t)|dx\\
  &=&\int_\Omega[h(x,t)-h(x,t_0)]\chi(\psi(x,t))J(\psi(x,t),t)dx\\
  &=&\int_\Omega[h(x,t)-h(x,t_0)][\chi(\psi(x,t))J(\psi(x,t),t)-\chi(\psi(x,t_0))J(\psi(x,t_0),t_0)]dx\\
  &&+\int_\Omega[h(x,t)-h(x,t_0)]\chi(\psi(x,t_0))J(\psi(x,t_0),t_0)dx\\
  &=:&R_{11}+R_{12}.
\end{eqnarray*}
Since $\chi(\psi(x,t_0))J(\psi(x,t_0),t_0)\in L^2$, guaranteed by Proposition \ref{cs}, one has
\begin{equation}
  R_{12}\rightarrow0,\quad\mbox{as }t\rightarrow t_0. \label{LGG1}
\end{equation}
Recalling (\ref{LG1-1}) and (\ref{LG6}), one has
\begin{eqnarray*}
  &&\chi(\psi(x,t))J(\psi(x,t),t)-\chi(\psi(x,t_0))J(\psi(x,t_0),t_0)\\
  &=&\int_{t_0}^t[\chi(\psi(x,s))(\partial_tJ(\psi(x,s),s)+\nabla J(\psi(x,s),s)\cdot\partial_t\psi(x,s))\\
  &&+\nabla\chi(\psi(x,s))\cdot\partial_t\psi(x,s)J(\psi(x,s),s)]ds \\
  &=&-\int_{t_0}^t(u(x,s)\cdot\nabla)\psi(x,s)\cdot[\chi(\psi(x,s)) \nabla J(\psi(x,s),s)\\
  &&+\nabla\chi(\psi(x,s))J(\psi(x,s),s)]ds-\int_{t_0}^t\chi(\psi(x,s))\text{div} u(x,s)J(\psi(x,s),s)ds\\
  &=&-\int_{t_0}^tu(x,s)\cdot A(\psi(x,s),s)\cdot[\chi(\psi(x,s)) \nabla J(\psi(x,s),s)\\
  &&+\nabla\chi(\psi(x,s))J(\psi(x,s),s)]ds-\int_{t_0}^t\chi(\psi(x,s))\text{div} u(x,s)J(\psi(x,s),s)ds
\end{eqnarray*}
and, thus, by Proposition \ref{cs}, it follows
\begin{eqnarray*}
  &&\|\chi(\psi(\cdot,t))J(\psi(\cdot,t),t)-\chi(\psi(\cdot,t_0))J(\psi(\cdot,t_0),t_0)\|\\
  &\leq&\left|\int_{t_0}^t\|u\|_\infty\|A(\psi(\cdot,s),s)\|\|\nabla\chi\|_\infty\|J\|_\infty+
  \|u\|_\infty\|A\|_\infty\|\chi\|_\infty\|\nabla J(\psi(\cdot,s),s)\| ds\right|\\
  &&+\left|\int_{t_0}^t\|\chi\|_\infty\|\text{div} u\|\|J\|_\infty ds\right|\\
  &\leq&C\left|\int_{t_0}^t\left(\|\nabla u\|^\frac12\|\nabla^2u\|^\frac12\right)ds\right|+C|t-t_0|\leq C|t-t_0|^\frac34.
\end{eqnarray*}
Thanks to this, one has
\begin{eqnarray}
  |R_{11}|
   &\leq& \|h(\cdot,t)-h(\cdot,t_0)\|
  \|\chi(\psi(\cdot,t))J(\psi(\cdot,t),t)-\chi(\psi(\cdot,t_0))J(\psi(\cdot,t_0),t_0)\|\nonumber\\
   &\leq& C|t-t_0|^\frac34.  \label{LGG2}
\end{eqnarray}
Using the H\"older inequality and by Proposition \ref{cs}, it follows
\begin{eqnarray}
  |R_2+R_4|&\leq&C\|\chi\|(\|h(\varphi(\cdot,t),t_0)-\xi(\varphi(\cdot,t))\|
  +\|h(\varphi(\cdot,t_0),t_0)-\xi(\varphi(\cdot,t_0))\|)\nonumber\\
  &\leq&C\|\chi\| \|h(\cdot, t_0)-\xi\| \leq C \varepsilon. \label{LGG3}
\end{eqnarray}
For $R_3$, it follows from the H\"older inequality and (\ref{LGG7}) that
\begin{eqnarray}
  |R_3|\leq\|\chi\|\|\xi(\varphi(y,t))-\xi(\varphi(y,t_0))\|\leq
   C |t-t_0|^\frac34.\label{LGG4}
\end{eqnarray}
Combining (\ref{LGG2})--(\ref{LGG4}), one gets
$$
\left|\int_\Omega[h(\varphi(y,t),t)-h(\varphi(y,t_0),t_0)]\chi(y)dy\right|\leq C |t-t_0|^\frac34+C\varepsilon+|R_{12}|.
$$
With the aid of this and recalling (\ref{LGG1}), one derives
$$
\int_\Omega[h(\varphi(y,t),t)-h(\varphi(y,t_0),t_0)]\chi(y)dy\rightarrow0,\quad\mbox{as }t\rightarrow t_0,
$$
which implies that $h(\varphi(\cdot,t),t)$ is weakly continuous in $L^2(\Omega)$ at any $t_0\in[0,{T_0}]$. Therefore, $h(\varphi(\cdot,t),t)\in C_w([0,{T_0}]; L^2)$.
\end{proof}

Finally, we prove the following lemma which is used in (\ref{INBYP}) during the proof of the uniqueness in the previous section.

\begin{lemma}
\label{cv}
Given a bounded domain $\Omega$ in $\mathbb R^3$. Let $\Phi\in W^{2,q}$ with $q\in(3,6)$ be a bijective mapping
on $\Omega$. Denote $x=\Psi(y)$ and $y=\Psi(x)$.
Then, it holds that
$$
\text{div}_y\left(\frac{\partial_{x_i}y}{\text{det}\nabla_xy}\right)=0 \quad\mbox{and}\quad \text{div}_x\left(\frac{\partial_{y_i}x}{\text{det}\nabla_yx}\right)=0.
$$
\end{lemma}

\begin{proof}
We only give the proof of the first identity while the second one can be proved in the same way.
In the proof of this lemma, for a matrix $A=(a_{ij})_{3\times3}$, we use $R_i(A)$, $M_{ij}(A)$, and
$A_{adj}$, respectively, to denote
the $i$-th row of $A$, the minor of the entry $a_{ij}$, and the classical adjoint of $A$
that is, $A_{adj}=((-1)^{i+j}M_{ij})^T$.
Denote
$\nabla_xy=(\partial_{x_i}y_j)_{3\times3}$ and $\nabla_yx=(\partial_{y_i}x_j)_{3\times3}$. Then, the chain rule
gives $\nabla_xy\nabla_yx=I$. Thanks to this, one deduces
$$
\frac{\nabla_xy}{\text{det} \nabla_xy}=\frac{(\nabla_yx)^{-1}}{{\text{det} \nabla_xy}}=\frac{(\nabla_yx)_{adj}}{\text{det}\nabla_yx{\text{det} \nabla_xy}}=(\nabla_yx)_{adj}.
$$
Therefore,
$$
\text{div}_y\left(\frac{\partial_{x_i}y}{\text{det}\nabla_xy}\right)=\text{div}_y\Big(R_i((\nabla_yx)_{adj})\Big).
$$
It remains to show that $\text{div}_y\Big(R_i((\nabla_yx)_{adj})\Big)=0$ for $i=1,2,3.$ We only prove the case $i=1$,
the proofs for $i=2,3$ are the same. By definition, one has
$$
\text{div}_y\Big(R_1((\nabla_yx)_{adj})\Big)=\left|\begin{array}{ccc}
                                                     \partial_{y_1} & \partial_{y_1}x_2 & \partial_{y_1}x_3 \\
                                                     \partial_{y_2} & \partial_{y_2}x_2 & \partial_{y_2}x_3 \\
                                                     \partial_{y_3} & \partial_{y_3}x_2 & \partial_{y_3}x_3
                                                   \end{array}
\right|,
$$
where the determinant is understood by expanding along the first column. By direct calculations, one can verify that the above determinant is identically zero and, thus, the conclusion holds.
\end{proof}

\section*{ACKNOWLEDGMENTS}
{This work was supported in part by the National
Natural Science Foundation of China (11971009 and 11871005), by the Key Project of National Natural Science Foundation of China (12131010), and by
the Guangdong Basic and Applied Basic Research Foundation (2019A1515011621,
2020B1515310005, 2020B1515310002, and 2021A1515010247).}

\end{document}